\documentclass[10pt]{article}

\usepackage[utf8]{inputenc}
\usepackage[dvipsnames]{xcolor}
\usepackage{float}
\usepackage{cancel}
\usepackage{leftidx}
\usepackage{amsmath}
\usepackage{amscd}
\usepackage{amssymb}
\usepackage{rotating}
\usepackage{amsfonts}
\usepackage{mathrsfs}
\usepackage{amsthm}
\usepackage{verbatim}
\usepackage{simpler-wick}
\usepackage{diagbox}
\usepackage[overload]{empheq}
\usepackage{cases}
\usepackage[top=1in, bottom=1in, left=1.25in, right=1.25in]{geometry}

\usepackage[colorlinks,linkcolor=Cyan, anchorcolor=Cyan, urlcolor=blue, citecolor=blue]{hyperref}
\usepackage{cleveref} 
 \usepackage[hyperpageref]{backref}

\usepackage{titlesec}
       \titleformat{\chapter}[display]
             {\normalfont\Large\bfseries}{\thechapter}{11pt}{\Large}
       \titleformat{\section}
             {\normalfont\large\bfseries}{\thesection}{11pt}{\large}
       \titlespacing*{\chapter}{0pt}{0pt}{15pt} 
       \titlespacing*{\section}{0pt}{3.5ex plus 1ex minus .2ex}{2.3ex plus .2ex}

    \usepackage[titletoc]{appendix}

\input xy
\xyoption{all}

\allowdisplaybreaks[3]

\newcommand{\pqed}{\hfill\qedsymbol\\}
\newcommand{\nn}{\nonumber}
\newcommand{\sfh}{\mathsf{h}}
\newcommand{\tsfh}{\tilde{\mathsf{h}}}
\newcommand{\sfF}{\mathsf{F}}
\newcommand{\hypergeometricF}{\mathfrak{F}}
\newcommand{\PZc}{\mathrm{c}}
\newcommand{\Mbar}{\overline{\mathcal{M}}}
\newcommand{\Mtd}{\widetilde{\mathcal{M}}}
\newcommand{\sqp}{\diamond}

\newcommand{\Nsim}{\stackrel{N}{\sim}}
\newcommand{\Nlsim}{\stackrel{(N,l)}{\sim}}
\newcommand{\tauG}{\mathrm{g}}
\newcommand{\sfbd}{\mathbf{d}^{\mathbf{d}}}

\newtheorem{theorem}{Theorem}[section]
\newtheorem{proposition}[theorem]{Proposition}
\newtheorem{lemma}[theorem]{Lemma}
\newtheorem{corollary}[theorem]{Corollary}
\newtheorem{conjecture}[theorem]{Conjecture}
\newtheorem{theorem/definition}[theorem]{Theorem/Definition}
\newtheorem{question}[theorem]{Question}

\theoremstyle{remark}
\newtheorem{remark}[theorem]{Remark}

\theoremstyle{definition}
 
\newtheorem{definition}[theorem]{Definition}

\begin{document}

\title
{\large{\textbf{ON GENUS 1 GROMOV-WITTEN INVARIANTS OF FANO COMPLETE INTERSECTIONS}}}
\author{\normalsize Xiaowen Hu\\
\date{}
\\
}
\maketitle

\begin{abstract}
We study  genus 1 Gromov-Witten invariants of Fano  complete intersections in the projective spaces. Among other things, we show a reconstruction theorem for  genus 1 invariants with only ambient insertions, and compute the genus 1 invariants with 1 marked point. For cubic hypersurfaces of dimension $\neq 4$ and odd dimensional intersections of two quadrics, we obtain a complete reconstruction theorem for genus 1 Gromov-Witten invariants.
\end{abstract}

\tableofcontents

\section{Introduction}

The genus 1 Gromov-Witten invariants of Calabi-Yau  hypersurfaces in projective spaces are first computed by Zinger (\cite{Zin08}, \cite{Zin09b}). The result is then extended to Calabi-Yau complete intersections by Popa (\cite{Pop12}). Afterwards, there are different new proofs (\cite{KL18} for CY complete intersections \cite{GR19} and \cite{CGLZ20} for quintic 3-folds).

In this paper, we study the genus 1 Gromov-Witten invariants of Fano complete intersections of dimension $\geq 3$.  One difficulty in this case is that there exist nontrivial invariants of arbitrary lengths, contrary to the Calabi-Yau case. Thus it is natural to first use the tautological relations to reduce the computation to invariants as simple as possible. Even if one wants to compute the genus 1 invariants with only ambient insertions (cup products of hyperplane classes), such a process definitely needs information of the genus 0 invariants with primitive insertions. This gives one motivation for our previous paper \cite{Hu15}, where we study the full quantum cohomology  (i.e. allowing arbitrary primitive insertions)  of \emph{non-exceptional} complete intersections. For the latter notion, see Definition \ref{def-exceptional}; by \cite[\S 3]{BaM04} and \cite{Hu21}, the exceptional complete intersections  are exactly the ones that have semi-simple quantum cohomology.

The approach of \cite{Hu15} is using the monodromy invariance and thus the symmetric reduction of the WDVV equation, and  geometry in some special cases (= the cases in Theorem \ref{thm-reconstruction-cubic-(2,2)-intro}). We obtained a conditional algorithm\footnote{The algorithm is implemented as a Macaulay2 package\\ \url{https://github.com/huxw06/Quantum-cohomology-of-Fano-complete-intersections}}; \emph{conditional} means that the algorithm will automatically check some conditions in the running for the recursion to work. That the conditions are always true, or equivalently  the effectiveness of the algorithm, is the \emph{double root recursion conjecture}. See \cite[\S 1]{Hu15} for a full account.

There is a master space localization method \cite{FL19} which in genus 1 produces an algorithm that computes at least the invariants with ambient insertions. Recently there is another approach \cite{ABPZ21}, which shows a new degeneration formula, and gives  an effective algorithm in all genera. However, at present, this algorithm seems not easy to give explicit computations. 

In this paper, as an application and a continuation of \cite{Zin09b} and \cite{Hu15},  we compute the symmetric reduction, which arises from the big monodromy groups of smooth projective complete intersections, of the tautological relations for genus 1 Gromov-Witten invariants induced by the Getzler's relation in $H^*(\Mbar_{1,4})$ \cite[Theorem 1.8]{Get97}, and compute the genus 1 invariants with 1 marked point. The main results are the following theorems.

\begin{theorem}[= Proposition \ref{prop-reconstruction-genus1-fromLength1}]\label{thm-reconstruction-genus1-fromLength1-intro}
For any non-exceptional Fano complete intersection $X$ of dimension $\geq 3$ in projective spaces, and any $k\geq 0$, one can reconstruct the  genus 1 Gromov-Witten invariants of $X$ with $k$ primitive insertions, from the following:
\begin{enumerate}
       \item[(i)] the genus 0 Gromov-Witten invariants of $X$;
       \item[(ii)] the genus 1 Gromov-Witten invariants of $X$ with less than $k$ primitive insertions;
       \item[(iii)] the genus 1 Gromov-Witten invariants of $X$ with  $k$ primitive insertions and at most one ambient insertions.
 \end{enumerate} 
In particular, one can reconstruct the  genus 1 Gromov-Witten invariants of $X$ with only \emph{ambient} insertions, from the genus 0 Gromov-Witten invariants of $X$, and genus 1 Gromov-Witten invariants of $X$ with only one marked point, using the monodromy invariance.
\end{theorem}

\begin{theorem}[= Proposition \ref{prop-reconstruction-cubic-(2,2)}]\label{thm-reconstruction-cubic-(2,2)-intro}
In the following two cases of   Fano complete intersections $X$ in projective spaces, one can reconstruct the  genus 1 Gromov-Witten invariants of $X$ with arbitrary insertions, from the genus 0 Gromov-Witten invariants of $X$, and genus 1 Gromov-Witten invariants of $X$ with only one marked point, using the monodromy invariance:
\begin{enumerate}
      \item[(i)] cubic hypersurfaces of dimension 3 and dimension $\geq 5$;
      \item[(ii)] complete intersections of odd dimensions $\geq 3$, of two quadrics.
\end{enumerate}
\end{theorem}
Here Gromov-Witten invariants are by default \emph{primary} ones, i.e. without inserting $\psi$-classes. By the topological recursion relations in genus 0 and 1 (e.g. \cite[(3)]{Get98}), the non-primary invariants in genus 0 and 1 can be reconstructed from the primary ones.

\begin{theorem}[= Theorem \ref{thm-genus1-GWInv-Fano}]\label{thm-genus1-GWInv-Fano-intro}
Let $X$ be a smooth complete intersection of multidegree $\mathbf{d}=(d_1,\dots,d_r)$ in $\mathbb{P}^{n-1}$, with Fano index $\nu_{\mathbf{d}}=n-|\mathbf{d}|\geq 1$. Suppose $d_i\geq 2$ for $1\leq i\leq r$.
 Then for $0\leq b\leq \lfloor\frac{n-1}{\nu_{\mathbf{d}}}\rfloor$,
\begin{eqnarray}\label{eq-genus1-GWInv-intro}
&& \langle \sfh_{1+\nu_{\mathbf{d}}b}\rangle_{1,b}\nn\\
&=& \frac{1}{2}
\mathrm{Coeff}_{q^{b}}\Big\{
\frac{\Theta_{1+\nu_{\mathbf{d}}b}^{(0)}(q)\big(\sum_{p=0}^{n-1-r}      \Theta_{p}^{(1)}(q)\Theta_{n-1-r-p}^{(0)}(q)
+\sum_{p=1}^{r}   \Theta_{n-p}^{(1)}(q)\Theta_{n-1-r+p}^{(0)}(q)\big)}{\Phi_0(q)}\Big\}\nn\\
&&+\frac{n}{24}
\mathrm{Coeff}_{q^{b}}\bigg\{ \big(\frac{n-1}{2}-\sum_{k=1}^r \frac{1}{d_k}\big)\Big(1-\sum_{\beta=0}^{\infty}\tilde{\PZc}_{1+\nu_{\mathbf{d}}b,1+\nu_{\mathbf{d}}b- \nu_{\mathbf{d}}\beta}^{(\beta)}q^{\beta}\big(L(q)^{1+\nu_{\mathbf{d}}b- \nu_{\mathbf{d}}\beta}-1\big)\Big)\nn\\
&&-L(q)'\sum_{\beta=0}^{\infty}\tilde{\PZc}_{1+\nu_{\mathbf{d}}b,1+\nu_{\mathbf{d}}b- \nu_{\mathbf{d}}\beta}^{(\beta)}q^{\beta+1}
\binom{1+\nu_{\mathbf{d}}b- \nu_{\mathbf{d}}\beta}{2} L(q)^{1+\nu_{\mathbf{d}}b- \nu_{\mathbf{d}}\beta-2}\nn\\
&&-\frac{\Phi'_0(q)}{\Phi_0(q)}\sum_{\beta=0}^{\infty}\tilde{\PZc}_{1+\nu_{\mathbf{d}}b,1+\nu_{\mathbf{d}}b- \nu_{\mathbf{d}}\beta}^{(\beta)}q^{\beta+1}
(1+\nu_{\mathbf{d}}b- \nu_{\mathbf{d}}\beta)  L(q)^{1+\nu_{\mathbf{d}}b- \nu_{\mathbf{d}}\beta-1} \nn\\
&&-\sum_{\beta=0}^{\infty}\tilde{\PZc}_{1+\nu_{\mathbf{d}}b,1+\nu_{\mathbf{d}}b- \nu_{\mathbf{d}}\beta-1}^{(\beta)}q^{\beta}\big(L(q)^{1+\nu_{\mathbf{d}}b- \nu_{\mathbf{d}}\beta-1}-1\big)\bigg\}\nn\\
&&-\frac{\prod_{k=1}^r d_k}{24}
\mathrm{Res}_{w=0 }\Big\{
\frac{(1+w)^n (\tilde{\PZc}_{1+\nu_{\mathbf{d}}b,0}^{(b)}+\tilde{\PZc}_{1+\nu_{\mathbf{d}}b,1}^{(b)}w)}{w^{n-r}\prod_{k=1}^r(d_k w+1)}
\Big\}.
\end{eqnarray}
\end{theorem} 
Here $\tilde{\PZc}_{i,j}^{(b)}\in \mathbb{Q}$ are constants introduced  by Popa-Zinger \cite{PoZ14}, which are defined recursively:
\begin{equation}\label{eq-def-numbers-c-intro}
      \sum_{\beta=0}^{\infty}\sum_{l=0}^\infty \PZc_{p,l}^{(\beta)} w^l q^{\beta}
      =\sum_{\beta=0}^{\infty}q^{\beta}\frac{(w+\beta)^p \prod_{k=1}^{r}\prod_{i=1}^{d_k \beta}(d_k w+i)}
      {\prod_{j=1}^{\beta}(w+j)^n},
\end{equation}
\begin{equation}\label{eq-def-numbers-tildec-intro}
      \sum_{\begin{subarray}{c}\beta_1+\beta_2=\beta\\ \beta_1,\beta_2\geq 0\end{subarray}}
      \sum_{k=0}^{p- \nu_{\mathbf{d}}\beta_1}\tilde{\PZc}_{p,k}^{(\beta_1)}\PZc_{k,l}^{(\beta_2)}=\delta_{\beta,0}
      \delta_{p,l},\ \mbox{for}\ \beta,l\in \mathbb{Z}_{\geq 0},\ l\leq p- \nu_{\mathbf{d}}\beta.
\end{equation}
The series $\Theta^{(0)}_p(q)$ and $\Theta^{(1)}_p(q)$ are Fano counterparts to the series with the same notations in  the cases of Calabi-Yau hypersurfaces \cite[(3.42), (3,43)]{Zin09b}. We have (see Lemma \ref{lem-Theta(0)-formula-0} and \ref{lem-Theta(1)-formula-0})
\begin{eqnarray}\label{eq-Theta(0)-formula-0-intro}
\Theta_{p}^{(0)}(q)
= \Phi_0(q)\sum_{\beta=0}^{\infty}\tilde{\PZc}_{p,p- \nu_{\mathbf{d}}\beta}^{(\beta)}q^{\beta}
L(q)^{p- \nu_{\mathbf{d}}\beta},
\end{eqnarray}
and
\begin{eqnarray}\label{eq-Theta(1)-formula-0-intro}
\Theta_{p}^{(1)}(q)
&=& \Phi_0(q)\sum_{\beta=0}^{\infty}\tilde{\PZc}_{p,p- \nu_{\mathbf{d}}\beta-1}^{(\beta)}q^{\beta}
L(q)^{p- \nu_{\mathbf{d}}\beta-1}\nn\\
&&+ \Phi_1(q)\sum_{\beta=0}^{\infty}\tilde{\PZc}_{p,p- \nu_{\mathbf{d}}\beta}^{(\beta)}q^{\beta}
L(q)^{p- \nu_{\mathbf{d}}\beta}\nn\\
&&+\Phi'_0(q)\sum_{\beta=0}^{\infty}\tilde{\PZc}_{p,p- \nu_{\mathbf{d}}\beta}^{(\beta)}q^{\beta+1}
(p- \nu_{\mathbf{d}}\beta)
L(q)^{p- \nu_{\mathbf{d}}\beta-1}\nn\\
&&+L(q)'\Phi_0(q)\sum_{\beta=0}^{\infty}\tilde{\PZc}_{p,p- \nu_{\mathbf{d}}\beta}^{(\beta)}q^{\beta+1}
\binom{p- \nu_{\mathbf{d}}\beta}{2}
L(q)^{p- \nu_{\mathbf{d}}\beta-2}.
\end{eqnarray}
The series of $q$, $L(q)$, $\Phi_0(q)$ and $\Phi_1(q)$ in (\ref{eq-genus1-GWInv-intro}), (\ref{eq-Theta(0)-formula-0-intro}) and (\ref{eq-Theta(1)-formula-0-intro}) in turn have the following formulae (Proposition \ref{prop-formula-L0} and \ref{prop-formula-Phi0-Phi1})
\begin{equation}\label{eq-formula-L0(q)-intro}
      L(q)=\sum_{k=0}^{\infty}\frac{\prod_{i=1}^{k-1}(k|\mathbf{d}|+1-in)}{k!}(\frac{\mathbf{d}^{\mathbf{d}}q}{n})^k,
\end{equation}
\begin{equation}\label{eq-formula-Phi0-intro}
      \Phi_0(q)=L(q)^{\frac{r+1}{2}}\cdot \big(1+ \mathbf{d}^{\mathbf{d}}(1-\frac{|\mathbf{d}|}{n})q L(q)^{|\mathbf{d}|}\big)^{-\frac{1}{2}},
      \end{equation}
\begin{eqnarray}\label{eq-formula-Phi1-intro}
\Phi_1(q)&=&\frac{3 r^2-2  |\mathbf{d}|\sum_{k=1}^{r}\frac{1}{d_k} -1}{24|\mathbf{d}|}L(q)^{\frac{r-1}{2}}\big(L(q)-1\big) \big(1+ \mathbf{d}^{\mathbf{d}}(1-\frac{|\mathbf{d}|}{n})q L(q)^{|\mathbf{d}|}\big)^{-\frac{1}{2}}\nn\\
&&+\frac{L(q)^{\frac{r-1}{2}}\cdot \big(1+ \mathbf{d}^{\mathbf{d}}(1-\frac{|\mathbf{d}|}{n})q L(q)^{|\mathbf{d}|}\big)^{-\frac{7}{2}} }{24|\mathbf{d}|n^3}
\times\Big(|\mathbf{d}| ^3   \left(|\mathbf{d}|  n-|\mathbf{d}| -3 r^2+1\right)L(q)\nn\\
&&+|\mathbf{d}| ^2 n  \left(2 |\mathbf{d}| ^2-6 |\mathbf{d}|  n-6 |\mathbf{d}|  r+3 n^2+6 n r+n+3 r^2-1\right)L(q)^n\nn\\
&&+3 |\mathbf{d}| ^2  (n-|\mathbf{d}| ) \left(|\mathbf{d}|  n-|\mathbf{d}| -3 r^2+1\right) L(q)^{n+1}\nn\\
&&+|\mathbf{d}|  n (n-|\mathbf{d}| ) \left(4 |\mathbf{d}| ^2-5 |\mathbf{d}|  n-12 |\mathbf{d}|  r-2 n^2+6 n r+n+6 r^2-2\right)L(q)^{2 n} \nn\\
&&+3 |\mathbf{d}|   (n-|\mathbf{d}| )^2  \left(|\mathbf{d}|  n-|\mathbf{d}| -3 r^2+1\right)L(q)^{2 n+1}\nn\\
&&+n (n-|\mathbf{d}| )^2 \left(2 |\mathbf{d}| ^2+|\mathbf{d}|  n-6 |\mathbf{d}|  r+3 r^2-1\right)L(q)^{3 n}\nn \\
&&+ (n-|\mathbf{d}| )^3  \left(|\mathbf{d}|  n-|\mathbf{d}| -3 r^2+1\right)L(q)^{3 n+1}\Big),
\end{eqnarray}
where (following the notations in \cite{PoZ14})
\begin{equation}\label{eq-constants-related-to-d}
      |\mathbf{d}|=\sum_{i=1}^r d_i,\ \mathbf{d}^{\mathbf{d}}=\prod_{i=1}^r d_i^{d_i},\ \mathbf{d}!=\prod_{i=1}^{r}d_i!.
\end{equation}

The proof of Theorem \ref{thm-genus1-GWInv-Fano-intro}  follows the line of the approach invented by Zinger in the Calabi-Yau cases (CY hypersurfaces in \cite{Zin09b}, and CY complete intersections in \cite{Pop12}). Namely, we compute the reduced genus 1 GW invariants by the Atiyah-Bott localization theorem \cite{AB84}, and then compute the difference between the reduced invariants and the standard ones. In the localization computation, we find a way that greatly simplifies the computations in \cite{Zin09b} and \cite{Pop12}: we evaluate the equivariant generating function by equating  the equivariant parameters to the $n$-th roots of unity (see Lemma \ref{lem-calY-SpecializedAtRootsOfUnity} and Corollary \ref{cor-equivariant-Z-SpecializedAtRootsOfUnity}). This avoids  elaborate use of symmetric functions in loc. cit.

The formula (\ref{eq-genus1-GWInv-intro}) is  not as neat as that for the Calabi-Yau complete intersections  (compare to \cite[Theorem 2]{Zin08} and \cite[Theorem 1]{Pop12}). This is due to the  positive Fano index $\nu_{\mathbf{d}}$. In Section \ref{sec:simplifications-A(q)} we propose some conjectural formulae related to the first row of the right hand side of (\ref{eq-genus1-GWInv-intro}), based on numerical examples and residue computations in lower degrees.

Theorem \ref{thm-genus1-GWInv-Fano-intro} together with Theorem \ref{thm-reconstruction-genus1-fromLength1-intro} gives a complete reconstruction result for genus 1 GW invariants of non-exceptional Fano complete intersections with only ambient insertions, assuming the double root recursion conjecture \cite[Conjecture 1.14]{Hu15}. 
We regard the two cases in Theorem \ref{thm-reconstruction-cubic-(2,2)-intro} as the \emph{border cases} of the non-exceptional Fano complete intersections. 
By \cite[Theorem 1.10, 1.11, 1.13]{Hu15}, the genus GW invariants in the border cases can be reconstructed by essentially linear recursions (and do not need the double root recursion). Combined with
Theorem \ref{thm-genus1-GWInv-Fano-intro} and Theorem \ref{thm-reconstruction-cubic-(2,2)-intro}, we obtain  a complete reconstruction result for genus 1 GW invariants in the border cases. 

However, in the cases other than these two border cases,
our approach does not achieve a complete reconstruction result for genus 1 GW invariants with an arbitrary number of primitive insertions.
We have computed the equations  obtained by the symmetric reduction of Getzler's relation, and only 3 of them are presented in this paper (see (\ref{eq-symmetridReduction-0<=a,b,c,d<=n}), (\ref{eq-symmetricReduction-0<=a,b,c<=n,n+1<=d<=n+m}), and (\ref{eq-symReduction-n+1<=a,b,c,d<=n+m})). Then numerical experiments suggest that  the undetermined initial values in Theorem \ref{thm-reconstruction-genus1-fromLength1-intro} cannot be determined by these equations.
This is reminiscent of the case of the genus 0 initial values $F^{(k)}(0)$ (defined in (\ref{eq-F&G-s-derivatives})) in \cite{Hu15}, which turn out to be determined by the double root recursion. Our experience from \cite{Hu15} then suggests the following question:
\begin{question}
Are the undetermined initial values in Theorem \ref{thm-reconstruction-genus1-fromLength1-intro} reconstructible from symmetric reductions of tautological relations in higher genera? More generally, 
to what extent can the cohomological field theory associated with a Fano complete intersection be determined by its monodromy invariance?
\end{question}

The paper is organized as follows. In Section \ref{sec:recap} we collect the results in \cite{Hu15} on genus 0 Gromov-Witten invariants of Fano complete intersections that will be used in this paper.  In Section \ref{sec:symmetricReduction} we compute the symmetric reduction of the genus 1 Getzler relation of Gromov-Witten invariants, and prove Theorem \ref{thm-reconstruction-genus1-fromLength1-intro} and \ref{thm-reconstruction-cubic-(2,2)-intro}. In Section \ref{sec:ReducedGenus1GW-Fano} and Section \ref{sec:SvR}, we compute the genus 1 invariants with 1 marked point. Sections \ref{sec:formulae-hypergeometricSeries-Fano} and \ref{sec:simplifications-A(q)} contain formulae  related to hypergeometric series associated with Fano complete intersections that are used in Section \ref{sec:ReducedGenus1GW-Fano}.\\

\emph{Acknowledgment}\quad 
 I am grateful to Huai-Liang Chang, Honglu Fan,  Hua-Zhong Ke, and Yijie Lin for helpful discussions. 
  This work is supported by  NSFC 11701579 and NSFC 11831017.

\section{Recap of quantum cohomology of Fano complete intersections}\label{sec:recap}

Let $X$ be a smooth Fano complete intersection of dimension $\geq 3$ in a projective space $\mathbb{P}$.  Denote the Fano index of $X$ by $\mathsf{a}_X$. Denote the subgroup of $H^*(X;\mathbb{C})$ generated by the image of $H^*(\mathbb{P};\mathbb{C})$, via the imbedding, by 
$H_{\mathrm{amb}}^*(X;\mathbb{C})$. Then $H^*(X;\mathbb{C})=H_{\mathrm{amb}}^*(X;\mathbb{C})\oplus H_{\mathrm{prim}}^*(X;\mathbb{C})$.
Let $\gamma_0,\dots,\gamma^M$ be a basis of $H^*(X;\mathbb{C})$. 
Let $T^0,\dots,T^M$ be the dual basis with respect to $\gamma_0,\dots,\gamma^M$. Then the  generating functions of genus $0$ and $1$ Gromov-Witten invariants are  defined respectively as
 \begin{equation}\label{eq-def-generatingFunction-genus0}
      F(T^0,\dots,T^{M})=\sum_{k\geq 0} \sum_{\beta=0}^{\infty} \frac{1}{k!}\big\langle \sum_{i=0}^M \gamma_i T^i,\dots,\sum_{i=0}^M\gamma_i T^i\big\rangle_{0,k,\beta},
 \end{equation}
\begin{equation}\label{eq-def-generatingFunction-genus1}
      G(T^0,\dots,T^{M})=\sum_{k\geq 0} \sum_{\beta=0}^{\infty} \frac{1}{k!}\big\langle \sum_{i=0}^M \gamma_i T^i,\dots,\sum_{i=0}^M\gamma_i T^i\big\rangle_{1,k,\beta},
 \end{equation}
where the GW invariants outside of the stable range are defined to be zero, as a convention. 

The evaluation of a function $h(T^0,\dots,T^M)$ at $T^0=\dots=T^M=0$ is simply denoted by $h(0)$. For example, if each $\gamma_i$ has pure complex degree $b_i$, then the dimension formula of the virtual fundamental class yields
\begin{equation}\label{eq-coefficientsOfF}
     \frac{\partial^k F}{\partial T^{i_1}\dots \partial T^{i_k}}(0)=
      \begin{cases}
      \langle \gamma_{i_1},\dots,\gamma_{i_k}\rangle_{1,k,\beta},& \mbox{if}\ \beta=\frac{\dim X-3+\sum_{p=1}^k b_{i_p}-k}{\mathsf{a}_X}\in \mathbb{Z}_{\geq 0};\\
      0,& \mbox{otherwise},
      \end{cases}
\end{equation}
and
\begin{equation}\label{eq-coefficientsOfG}
      \frac{\partial^k G}{\partial T^{i_1}\dots \partial T^{i_k}}(0)=
      \begin{cases}
      \langle \gamma_{i_1},\dots,\gamma_{i_k}\rangle_{1,k,\beta},& \mbox{if}\ \beta=\frac{\sum_{p=1}^k b_{i_p}-k}{\mathsf{a}_X}\in \mathbb{Z}_{\geq 0};\\
      0,& \mbox{otherwise}.
      \end{cases}
\end{equation}

We recall the following result \cite[Theorem 3.7]{Hu15}, which is a consequence of the monodromy invariance of Gromov-Witten invariants (e.g. \cite[Theorem 4.2]{LT98}) and the bigness of the monodromy groups of projective complete intersections (\cite{Del73}).
\begin{definition}\label{def-exceptional}
$X_n(\mathbf{d})$ is called an \emph{exceptional complete intersection} (in the projective spaces), if it is one of the following cases:
(i) $X_{n}(2)$;
(ii) $X_{n}(2,2)$, $n$ is even;
(iii) $X_2(3)$, i.e., a cubic surface.
Otherwise $X_n(\mathbf{d})$ is called \emph{non-exceptional}.
\end{definition}

\begin{theorem}\phantomsection\label{thm-monodromythm}
\begin{itemize}
\item[(i)] Let $X$ be an even-dimensional non-exceptional complete intersection of dimension $n\geq 4$. Suppose $\gamma_0,\cdots,\gamma_{n+m}$ be a basis of $H^{*}(X,\mathbb{C})$, where $\gamma_{0},\cdots,\gamma_{n}$ is a basis of $H_{\mathrm{amb}}^{*}(X,\mathbb{C})$, and 
$\gamma_{n+1},\cdots,\gamma_{n+m}$ is an orthonormal basis of $H^{n}_{\mathrm{prim}}(X,\mathbb{C})$, in the sense 
$Q(\alpha_{i},\alpha_{j})=\delta_{ij}$. Let 
$t^{0},\cdots, t^{n+m}$ be the dual basis of $\gamma_{0},\cdots,\gamma_{n+m}$. Then for  $g\geq 0$, the genus $g$ generating function $\mathcal{F}_{g}$ can be written in a unique way as a series in $t^{0},\cdots,t^{n}$ and $s$, where
\begin{eqnarray}\label{eq-sumS-evenDim}
s=\sum_{\mu=n+1}^{n+m}\frac{(t^{\mu})^2}{2}.
\end{eqnarray}
\item[(ii)] Let $X$ be an odd-dimensional complete intersection of dimension $n\geq 3$. Suppose $\gamma_0,\cdots,\gamma_{n+m}$ be a basis of $H^{*}(X,\mathbb{C})$, 
where $\gamma_{0},\cdots,\gamma_{n}$ is a basis of $H_{\mathrm{amb}}^{*}(X,\mathbb{C})$, and $\gamma_{n+1},\cdots,\gamma_{n+m}$ is a symplectic basis of $H^{n}_{\mathrm{prim}}(X,\mathbb{C})$. Let 
$t^{0},\cdots, t^{n+m}$ be the dual basis of $\gamma_{0},\cdots,\gamma_{n+m}$. Then for any $g\geq 0$, the genus $g$ generating function $\mathcal{F}_{g}$ can be written in a unique way as a series in $t^{0},\cdots,t^{n}$ and $s$ with the degree of $s$ not greater than $\frac{m}{2}$, where
\begin{eqnarray}\label{eq-sumS-oddDim}
s=-\sum_{\mu=n+1}^{n+\frac{m}{2}}t^{\mu}t^{\mu+\frac{m}{2}}.
\end{eqnarray}
\end{itemize}
\end{theorem}

Now let $X$ be a smooth non-exceptional Fano complete intersection of multidegree $\mathbf{d}=(d_1,\dots,d_r)$ in the projective space $\mathbb{P}^{n+r}$, and $n\geq 3$. Let $\sfh$ be the hyperplane class in $H^*(X;\mathbb{C})$. 
 We define
\begin{equation}
      \sfh_i=\underbrace{\sfh\cup\cdots\cup\sfh}_{\mbox{$i$ factors}},
\end{equation}
\begin{eqnarray}\label{eq-ss5-intro}
\tilde{\sfh}=\left\{
\begin{array}{cc}
\sfh, & \mathsf{a}_X\geq 2,\\ 
\sfh+\mathbf{d}!, & \mathsf{a}_X=1,
\end{array}\right.
\end{eqnarray}
and
\begin{equation}
      \tsfh_i=\underbrace{\tsfh\sqp\cdots\sqp\tsfh}_{\mbox{$i$ factors}}
\end{equation}
where $\sqp$ stands for the small quantum product. Let $M$ and $W$ be the transition matrices between $\sfh_i$ and $\tsfh_i$:
\begin{equation}\label{eq-h-tildeh-transform}
      \sfh_i=\sum_{j=0}^{n}M_i^j \tsfh_j,\ \tsfh_i=\sum_{j=0}^{n}W_i^j \sfh_j,\ \mbox{for } 0\leq i\leq n.
\end{equation}

\begin{lemma}[{\cite[Lemma 6.1]{Hu15}}]\label{lem-Matrix-MW-triangular}
The matrices $M$ and $W$ are unipotent and upper triangular. Namely, $M_{i}^i=W_{i}^i=1$ for $0\leq i\leq n$, and $M_i^j=W_i^j=0$ for $0\leq i<j\leq n$.
\end{lemma}

The Poincaré pairing matrix, with respect to the basis $\tsfh_i$'s, $(\tauG_{e,f})_{0\leq e,f\leq n}$ and its inverse matrix $(\tauG^{e,f})_{0\leq e,f\leq n}$ are given by (\cite[Lemma 6.2]{Hu15})
\begin{eqnarray}\label{eq-pairing1}
\tauG_{e,f}=\int_X\tsfh_e\cup\tsfh_f=
\left\{
\begin{array}{cc}
(\mathbf{d}^{\mathbf{d}})^{\frac{e+f-n}{\mathsf{a}_X}}\prod_{i=1}^{r}d_i & 
\mathrm{if}\ \frac{e+f-n}{\mathsf{a}_X}\in \mathbb{Z}_{\geq 0};\\
0, & \mathrm{otherwise},
\end{array}
\right.
\end{eqnarray}
\begin{eqnarray}\label{eq-pairing2}
\tauG^{e,f}=\frac{1}{\prod_{i=1}^{r}d_i}\cdot
\left\{
\begin{array}{cc}
- \mathbf{d}^{\mathbf{d}}, & \mathrm{if}\ e+f=n-\mathsf{a}_X;\\
1, & \mathrm{if}\ e+f=n;\\
0, & \mathrm{otherwise}.
\end{array}
\right.
\end{eqnarray}
We extend $1,\tsfh,\dots,\tsfh_n$ to a basis as in Theorem \ref{thm-monodromythm}. Namely, $\gamma_i=\tsfh_i$ for $0\leq i\leq n$, and if $n$ is even (resp. $n$ is odd), $\gamma_i$ for $n+1\leq i\leq n+m$ is an orthonormal (resp. symplectic) basis of $H^*_{\mathrm{prim}}(X)$. 
Let $\tau^0,\dots,\tau^{n+m}$ be the corresponding dual basis, and let $s$ be defined as (\ref{eq-sumS-evenDim}) (resp. as (\ref{eq-sumS-oddDim})). We extend the Poincaré matrix $(\tauG_{i,j})$  to allow $0\leq i,j\leq n+m$, and similarly the matrix $(\tauG^{i,j})$. For example, if $n$ is odd,
\begin{eqnarray}\label{eq-inverse-symplecticPairing}
\tauG^{i,i+\frac{m}{2}}=-1,\ \tauG^{i+\frac{m}{2},i}=1,\ \mbox{for}\ n+1\leq i\leq n+\frac{m}{2},& \tauG^{i,j}=0\ \mbox{for}\ |i-j|\neq \frac{m}{2}.
\end{eqnarray}
By Theorem \ref{thm-monodromythm}, the genus 0 and 1 generating functions $F$ and $G$ are the functions of $\tau^0,\dots,\tau^n$ and $s$.
 Then we define
\begin{eqnarray}\label{eq-F&G-s-derivatives}
F^{(k)}(\tau^0,\cdots,\tau^{n})=\Big(\frac{\partial ^{k}}{\partial s^{k}}F\Big) \Big |_{s=0},\
G^{(k)}(\tau^0,\cdots,\tau^{n})=\Big(\frac{\partial ^{k}}{\partial s^{k}}G\Big) \Big |_{s=0}.
\end{eqnarray}
When the above basis is chosen, we often use  compressed notations of the following form
\begin{equation}
      F_{i_1,\dots,i_k}:=\frac{\partial^k F}{\partial \tau^{i_1}\dots \partial \tau^{i_k}},\
      F_{ss}=\frac{\partial^2 F}{\partial s^2},
\end{equation}
and similarly for $G$. Commas in subscripts or superscripts are sometimes omitted to save space, when no confusion arises. 

The advantage of the coordinates $\tau_i$'s is that some  coefficients of $F$, which in view of (\ref{eq-coefficientsOfF}) are combinations of Gromov-Witten invariants with pure degree classes $\sfh_i$'s, have concise formulae.
Let us see some examples that will be used in this paper. By \cite[Lemma 6.5]{Hu15},
\begin{eqnarray}\label{eq-qp1.5}
F_{a,b,c}^{(0)}(0)=\left\{
\begin{array}{cc}
(\mathbf{d}^{\mathbf{d}})^{\frac{a+b+c-n}{\mathsf{a}_X}}\prod_{i=1}^{r}d_i, & \mathrm{if}\ \frac{a+b+c-n}{\mathsf{a}_X}\in \mathbb{Z}_{\geq 0}; \\
0, & \mathrm{otherwise}.
\end{array}
\right.
\end{eqnarray}
By \cite[Theorem 6.7]{Hu15},
\begin{eqnarray}\label{eq-qp4}
\sum_{e=0}^{n}\sfF_{a,b,c,e}^{(0)}(0)\tauG^{e,0}=\left\{
\begin{array}{cc}
\mathsf{c}(n,l,\mathbf{d})(\sfbd)^{l} \mathsf{q}^l, & \mathrm{if}\ a,b,c\geq 1\ \mathrm{and}\  l=\frac{a+b+c-1}{\mathsf{a}_X}\in \mathbb{Z}_{\geq 0};\\
0, &  \mathrm{otherwise},
\end{array}
\right.
\end{eqnarray}
where
\begin{eqnarray}\label{eq-def-functionC}
\mathsf{c}(n,l,\mathbf{d})&:=&
1+\sum_{i=n-l\mathsf{a}_X}^{n}\sum_{j=0}^n 
\frac{j-i}{\mathsf{a}_X} M_{j}^{i}W_{n}^{j}(\mathbf{d}^{\mathbf{d}})^{\frac{i-n}{\mathsf{a}_X}}  \nn\\
&&-\sum_{i=n-l\mathsf{a}_X}^{n}\sum_{j=0}^n 
\frac{j-i}{\mathsf{a}_X} M_{j}^{i}W_{n- \mathsf{a}_X}^{j}(\mathbf{d}^{\mathbf{d}})^{1+\frac{i-n}{\mathsf{a}_X}}.
\end{eqnarray}
By \cite[Theorem 6.8 and 6.10]{Hu15},
\begin{equation}\label{eq-F(1)i(0)-Kronecker}
     F_i^{(1)}(0)=\delta_{i,0}.
\end{equation}
The formulae (\ref{eq-qp1.5}) and (\ref{eq-qp4}) suggest the following convention: if $i_p>n$, then we set
\begin{equation}\label{eq-convention-subscript-outOfRange}
      F_{i_1,\dots,i_{p-1},i_p,i_{p+1},\dots,i_k}= \mathbf{d}^{\mathbf{d}}F_{i_1,\dots,i_{p-1},i_p- \mathsf{a}_X,i_{p+1},\dots,i_k}.
\end{equation}
We apply the same convention to $G_{i_1,\dots,i_k}$. Then  (\ref{eq-qp1.5}) and (\ref{eq-pairing2}) yield
\begin{equation}\label{eq-giventalRelation}
     \sum_{e=0}^n\sum_{f=0}^n F_{abe}^{(0)}(0)\tauG^{ei}F_{cdf}(0)=F_{a+b,c,d}^{(0)}(0).
\end{equation}
In fact, this is a direct consequence of (\cite[Corollary 9.3 and 10.9]{Giv96}) that in the small quantum cohomology ring
\begin{equation}
      \tsfh^{ n+1}= \mathbf{d}^{\mathbf{d}}\tsfh^{(n+1- \mathsf{a}_X)}.
\end{equation}

The following result on the symmetric reduction of WDVV equation for $F$ (\cite[Theorem 4.1 and 4.2]{Hu15}) will be used in the symmetric reduction of  Getzler's relation in genus 1.
\begin{theorem}\label{thm-mainthm1}
If $\dim X$ is even, the  WDVV equation for $F$ is equivalent to the WDVV equation for $F^{(0)}$ together with
\begin{multline}\label{eq-wdvv230}
\sum_{e=0}^{n}\sum_{f=0}^{n}F_{abe}\tauG^{ef}F_{sf}+2sF_{sab}F_{ss}=F_{sa}F_{sb},\ (resp.\ \mbox{mod}\ s^{m/2}\ \mbox{if}\ \dim X\ \mbox{is odd})\\
\mbox{for}\ 0\leq a,b\leq n,
\end{multline}
and
\begin{eqnarray}\label{eq-wdvv240}
\sum_{e=0}^{n}\sum_{f=0}^{n}F_{se}\tauG^{e,f}F_{sf}+2s(F_{ss})^2= 0.\ (resp.\ \mbox{mod}\ s^{m/2}\ \mbox{if}\ \dim X\ \mbox{is odd})
\end{eqnarray}
\end{theorem}

\section{Symmetric reduction of Getzler's relations}\label{sec:symmetricReduction}
Throughout this section, we suppose  that $X$ is a smooth non-exceptional Fano complete intersection of multidegree $\mathbf{d}=(d_1,\dots,d_r)$ in the projective space $\mathbb{P}^{n+r}$, $n=\dim X\geq 3$, and $m=\mathrm{rank} H^*_{\mathrm{prim}}(X)$. We  follow the notations introduced  in Section \ref{sec:recap}. First, we recall the effect of the Euler vector field on the genus 1 generating function $G$.

\begin{lemma}\label{lem-EulerField-genus1}
\begin{eqnarray}\label{eq-EulerField-genus1} 
     && \big(\sum_{i=0}^{n}\sum_{j=0}^n\sum_{k=0}^n
(1-i)W_j^i  M_{i}^k \tau^j\frac{\partial}{\partial \tau^k}
+\mathsf{a}_X\frac{\partial}{\partial \tau^1}
-\delta_{1,\mathsf{a}_X} \mathbf{d}! \frac{\partial}{\partial \tau^0}\big)G\nn\\
&=&-\frac{1}{24}\int_X c_{n-1}(T_X)c_1(T_X).
\end{eqnarray}
\end{lemma}
\begin{proof}
Let $t^0,\dots,t^n$ be the basis of $H^*_{\mathrm{amb}}(X;\mathbb{C})^{\vee}$ dual to $1,\sfh,\dots,\sfh_n$. The Euler vector field in terms of $t^i$ and $s$ is equal to
\begin{eqnarray*}\label{eq-EV-symmetricReducted}
E&=&\sum_{i=0}^{n}(1-i)t^{i}\frac{\partial}{\partial t^i}+(2-n)s\frac{\partial }{\partial s}+\mathsf{a}_X\frac{\partial}{\partial t^1}.
\end{eqnarray*}
By Lemma \ref{lem-Matrix-MW-triangular} and (\ref{eq-ss5-intro}) have
\[
M_1^k=\begin{cases}
- \mathbf{d}! \delta_{1,\mathsf{a}_X},& \mbox{if}\ k=0;\\
1,& \mbox{if}\ k=1;\\
0,& \mbox{otherwise}
\end{cases}
\]
It follows that  in the coordinates $\tau_i$'s we have
\begin{eqnarray*}
E&=&\sum_{i=0}^{n}\sum_{j=0}^n\sum_{k=0}^n
(1-i)W_j^i  M_{i}^k \tau^j\frac{\partial}{\partial \tau^k}
+(2-n)s\frac{\partial }{\partial s}
+\mathsf{a}_X\frac{\partial}{\partial \tau^1}
-\delta_{1,\mathsf{a}_X} \mathbf{d}! \frac{\partial}{\partial \tau^0}.
\end{eqnarray*}
Then the divisor equation, the dimension formula of virtual fundamental classes, and the formula of degree 0  GW invariants (\cite[\S VI.6.3]{Man99}) yield (\ref{eq-EulerField-genus1}).
\end{proof}

Then we recall Getzler's relation.
We use $\sum_{\mathscr{P}(a,b,c,\dots)}$ to denote the summation over all the permutations of $a,b,c,\dots$. 
We use $\sum_{\mathscr{P}(a,b,c,\dots)}(\pm)$ to denote the \emph{signed} summation over all the permutations of $a,b,c,\dots$. For example, 
\begin{equation}
      \sum_{\mathscr{P}(a,b,c,d)}(\pm)h(\tau^a,\tau^b,\tau^c,\tau^d)=\sum_{\sigma\in S_4}\mathrm{sgn}\big(\sigma(a,b,c,d)\big)h\big(\tau^{\sigma(a)},\tau^{\sigma(b)},\tau^{\sigma(c)},\tau^{\sigma(d)}\big)
\end{equation}
where $\mathrm{sgn}\big(\sigma(a,b,c,d)\big)=\pm1$ is determined by
\begin{equation}
\tau^{\sigma(a)}\tau^{\sigma(b)}\tau^{\sigma(c)}\tau^{\sigma(d)}
     = \mathrm{sgn}\big(\sigma(a,b,c,d)\big)\cdot \tau^a \tau^b \tau^c \tau^d.
\end{equation}
So there are nontrivial signs only if the dimension $n$ is odd, and some of $a,b,c,\dots$ lie in the range $[n+1,n+m]$, i.e. corresponding to primitive cohomology classes. 

 Getzler's relation \cite[Theorem 1.8]{Get97} for genus 1 GW invariants is
\begin{multline}\label{Grelation-F&G}
\sum_{\mathscr{P}(a,b,c,d)}(\pm)\sum_{\mu=0}^{n+m} \sum_{\nu=0}^{n+m} \sum_{\lambda=0}^{n+m} \sum_{\rho=0}^{n+m}\big(
3F_{ab\lambda}\tauG^{\lambda\mu}G_{\mu\nu}\tauG^{\nu\rho}F_{\rho cd}-4F_{ab\lambda}\tauG^{\lambda\mu}F_{\mu c \rho}\tauG^{\rho\nu}G_{\nu d}\\
-F_{ab\lambda}\tauG^{\lambda\mu}F_{\mu cd\rho}\tauG^{\rho\nu}G_{\nu} +2F_{abc\lambda}\tauG^{\lambda\mu}F_{\mu d\rho}\tauG^{\rho\nu}G_{\nu} \\
+\frac{1}{6}F_{abc\lambda}\tauG^{\lambda\mu}F_{\mu d \nu \rho}\tauG^{\nu\rho}+\frac{1}{24}F_{abcd\lambda}\tauG^{\lambda\mu}F_{\mu\nu\rho}\tauG^{\nu\rho}
-\frac{1}{4}F_{ab\nu \lambda}\tauG^{\lambda\mu}F_{\mu cd\rho}\tauG^{\nu\rho}\big)= 0.
\end{multline}
Here we have  arranged the order in  the products, the superscripts, and the subscripts so that (\ref{Grelation-F&G}) holds in both even and odd dimensions.

\begin{definition}
We use $\sim$ to express that the difference between both sides are combinations (i.e. sums of products) of terms involving only genus 0 GW invariants. 
We use $\Nlsim$ to express that the difference between both sides are combinations of terms involving only the following type of GW invariants:
\begin{enumerate}
      \item[(i)] genus 0 invariants;
      \item[(ii)] genus 1  invariants with at most $2N-2$ primitive insertions;
      \item[(iii)] genus 1  invariants with  $2N$ primitive insertions and at most $l-1$ ambient insertions.
\end{enumerate}
In particular, we denote $\stackrel{(N,0)}{\sim}$ by $\Nsim$.
\end{definition}

Our strategy to compute the symmetric reduction of Getzler's relation is the same as that of the WDVV equation (\cite[\S 4]{Hu15}). Namely, for $0\leq k\leq 4$ we assume $k$ of $a,b,c,d$ lies in $[n+1,n+m]$, and assume they are equal or not, and then compute  (\ref{Grelation-F&G}) in the coordinates $\tau_i$, $0\leq i\leq n$, and $s$. 
We present only 3 equations in Section \ref{sec:0<=a,b,c,d<=n} to \ref{sec:n+1<=a,b,c,d<=n+m} which will be used in this paper.
We adopt Einstein's summation convention in the range $[0,n]\cap \mathbb{Z}$. For example, 
\[
F_{abe}\tauG^{ei}F_{cdf}\tauG^{fj}G_{ij}:=\sum_{e=0}^n\sum_{f=0}^n\sum_{i=0}^n\sum_{j=0}^n F_{abe}\tauG^{ei}F_{cdf}\tauG^{fj}G_{ij}.
\]
and also
\[
F_{abce}^{(0)}(0)\tauG^{ei}y_{i+d}:=\sum_{i=0}^{n}\sum_{e=0}^nF_{abce}^{(0)}(0)\tauG^{ei}y_{i+d}.
\]

\subsection{\texorpdfstring{$0\leq a,b,c,d\leq n$}{0<=a,b,c,d<=n}}\label{sec:0<=a,b,c,d<=n}
\begin{lemma}\label{lem-symmetridReduction-0<=a,b,c,d<=n}
Suppose $0\leq a,b,c,d\leq n$. Then
\begin{eqnarray}\label{eq-symmetridReduction-0<=a,b,c,d<=n}
&&\sum_{\mathscr{P}(a,b,c,d)} \big(3F_{abe}\tauG^{ei}F_{cdf}\tauG^{fj}G_{ij}+6sF_{abs}F_{cdf}\tauG^{fj}G_{sj}+6sF_{abe}\tauG^{ei}F_{cds}G_{i s}\nn\\
&&+4sF_{abs}F_{cds}G_{s}+12s^2F_{abs}F_{cds}G_{ss}\nn\\
&&-4F_{abe}\tauG^{ei}F_{cfi}\tauG^{fj}G_{dj}-8sF_{abs}F_{cfs}\tauG^{fj}G_{dj}
-8sF_{abe}\tauG^{ei}F_{csi}G_{ds}\nn\\
&&-8sF_{abs}F_{cs}G_{ds}-16s^2F_{abs}F_{css}G_{ds}\nn\\
&&-F_{abe}\tauG^{ei}F_{icdf}\tauG^{fj}G_{j}-2sF_{abs}F_{scdf}\tauG^{fj}G_{j}-2sF_{abe}\tauG^{ei}F_{icds}G_{s}-4s^2F_{abs}F_{cdss}G_{s}\nn\\
&&+2F_{abce}\tauG^{ei}F_{idf}\tauG^{fj}G_{j}+4sF_{abcs}F_{sdf}\tauG^{fj}G_{j}+4sF_{abce}\tauG^{ei}F_{ids}G_{s}\nn\\
&&+4sF_{abcs}F_{ds}G_{s}+8s^2F_{abcs}F_{dss}G_{s}\big)\nn\\
&=& \sum_{\mathscr{P}(a,b,c,d)} \Big(- \frac{1}{6}\big(F_{abce}\tauG^{ei}F_{idjk}\tauG^{kj}+2s F_{abcs}F_{sdjk}\tauG^{kj}+2s F_{abce}\tauG^{ei}F_{idss}
+(-1)^n m F_{abce}\tauG^{ei}F_{ids}\nn\\
&&+\big(2(-1)^n m+4\big)s F_{abcs}F_{dss}+4s^2 F_{abcs}F_{dsss}\big)\nn\\
&&-\frac{1}{24}\big(F_{abcde}\tauG^{ei}F_{ijk}\tauG^{kj}+2s F_{abcds}F_{sjk}\tauG^{kj}+2sF_{abcde}\tauG^{ei}F_{iss}
+(-1)^n m F_{abcde}\tauG^{ei}F_{is}\nn\\
&&+\big(2(-1)^n m+4\big)s F_{abcds}F_{ss}+4s^2 F_{abcds}F_{sss}\big)\nn\\
&&+\frac{1}{4}( F_{abij}F_{cdkl}\tauG^{jk}\tauG^{il}+4s F_{abis}F_{cdsl}\tauG^{il}\nn\\
&&+ (-1)^n m F_{abs}F_{cds}+2s F_{abs}F_{cdss}+2s F_{abss}F_{cds}+4s^2 F_{abss}F_{cdss})\Big).
\end{eqnarray}
\end{lemma}
\begin{proof} 
When the dimension $n$ of $X$ is even, a direct computation gives
\begin{eqnarray}\label{eq-symmetridReduction-0<=a,b,c,d<=n-LHS-1-evenDim}
 &&\sum_{\mu=0}^{n+m} \sum_{\nu=0}^{n+m}\sum_{\lambda=0}^{n+m} \sum_{\rho=0}^{n+m}  F_{ab\lambda}\tauG^{\lambda\mu}G_{\mu\nu}\tauG^{\nu\rho}F_{\rho cd}\nn\\
&=& F_{abe}\tauG^{ei}F_{cdf}\tauG^{fj}G_{ij}+2sF_{abs}F_{cdf}\tauG^{fj}G_{sj}+2sF_{abe}\tauG^{ei}F_{cds}G_{i s}\nn\\
&&+2sF_{abs}F_{cds}G_{s}+4s^2F_{abs}F_{cds}G_{ss}.
\end{eqnarray}
When $n$ is odd, as suggested by (\ref{eq-inverse-symplecticPairing}), we extend the range of the superscripts of the primitive variables $\tau^i$ and the corresponding inverse pairings $\tauG^{i,j}$ from $[n+1,n+m]\cap \mathbb{Z}$ to $[n+1,\infty)\cap \mathbb{Z}$, and define
\begin{equation}
      \tau^{i+m}=-\tau^i,\ \mbox{for}\ i\geq n+1
\end{equation}
and
\begin{equation}
      \tauG^{i+m,j}=\tauG^{i,j+m}=-\tauG^{i,j},\ \mbox{for}\ i,j\geq n+1.
\end{equation}
Moreover, we set $\mu'=\mu+\frac{m}{2}$ for $\mu\geq n+1$. These conventions make calculations of Grassmann variables $\tau^i$ uniformly for $i\geq n+1$. For example, we have
\begin{equation*}
      \tauG^{i',i}=1,\ \frac{\partial F}{\partial \tau^i}=-\tau^{i'}   \frac{\partial F}{\partial s}. 
\end{equation*}
We adopt these conventions in the odd dimensional cases in this proof. Then
\begin{eqnarray*}\label{eq-symmetridReduction-0<=a,b,c,d<=n-LHS-1-oddDim}
 &&\sum_{\mu=0}^{n+m} \sum_{\nu=0}^{n+m}\sum_{\lambda=0}^{n+m} \sum_{\rho=0}^{n+m}  F_{ab\lambda}\tauG^{\lambda\mu}G_{\mu\nu}\tauG^{\nu\rho}F_{\rho cd}\\
&=& F_{abe}\tauG^{ei}F_{cdf}\tauG^{fj}G_{ij}
+\sum_{\mu=n+1}^{n+m} \tau^{\mu}F_{abs}(-\tau^{\mu'}G_{si})\tauG^{ij}F_{jcd}\\
&&+\sum_{\nu=n+1}^{n+m} F_{abi}\tauG^{ij}(-\tau^{\nu'}G_{sj})(-\tau^{\nu}F_{scd})\\
&&+\sum_{\mu=n+1}^{n+m} \sum_{\nu=n+1}^{n+m} \tau^{\mu}F_{abs} (-\delta_{\mu,\nu'}G_s+\tau^{\mu'}\tau^{\nu'}G_{ss})(-\tau^{\nu}F_{scd})\\
&=& F_{abe}\tauG^{ei}F_{cdf}\tauG^{fj}G_{ij}+2sF_{abs}F_{cdf}\tauG^{fj}G_{sj}+2sF_{abe}\tauG^{ei}F_{cds}G_{i s}\\
&&+4s^2F_{abs}F_{cds}G_{ss}+2sF_{abs}F_{cds}G_{s}.
\end{eqnarray*}
Namely, (\ref{eq-symmetridReduction-0<=a,b,c,d<=n-LHS-1-evenDim}) remains true when $n$ is odd. Similar computations show that for both even and odd $n$ we have
\begin{eqnarray*}
&&\sum_{\mu=0}^{n+m} \sum_{\nu=0}^{n+m} \sum_{\lambda=0}^{n+m} \sum_{\rho=0}^{n+m} F_{ab\lambda}\tauG^{\lambda\mu}F_{\mu c \rho}\tauG^{\rho\nu}G_{\nu d}\\
&=& F_{abe}\tauG^{ei}F_{cfi}\tauG^{fj}G_{dj}+2sF_{abs}F_{cfs}\tauG^{fj}G_{dj}
+2sF_{abe}\tauG^{ei}F_{csi}G_{ds}\\
&&+2sF_{abs}F_{cs}G_{ds}+4s^2F_{abs}F_{css}G_{ds},
\end{eqnarray*}
\begin{eqnarray*}
&&\sum_{\mu=0}^{n+m} \sum_{\nu=0}^{n+m} \sum_{\lambda=0}^{n+m} \sum_{\rho=0}^{n+m} F_{ab\lambda}\tauG^{\lambda\mu}F_{\mu cd\rho}\tauG^{\rho\nu}G_{\nu}\\
&=&F_{abe}\tauG^{ei}F_{icdf}\tauG^{fj}G_{j}+2sF_{abs}F_{scdf}\tauG^{fj}G_{j}+2sF_{abe}\tauG^{ei}F_{icds}G_{s}\nn\\
&&+2sF_{abs}F_{cds}G_{s}+4s^2F_{abs}F_{cdss}G_{s},
\end{eqnarray*}
and
\begin{eqnarray*}
&&\sum_{\mu=0}^{n+m} \sum_{\nu=0}^{n+m} \sum_{\lambda=0}^{n+m} \sum_{\rho=0}^{n+m}F_{abc\lambda}\tauG^{\lambda\mu}F_{\mu d\rho}\tauG^{\rho\nu}G_{\nu}\\
&=&F_{abce}\tauG^{ei}F_{idf}\tauG^{fj}G_{j}+2sF_{abcs}F_{sdf}\tauG^{fj}G_{j}+2sF_{abce}\tauG^{ei}F_{ids}G_{s}\\
&&+2sF_{abcs}F_{ds}G_{s}+4s^2F_{abcs}F_{dss}G_{s}.
\end{eqnarray*}
Combining these yields that the sum of the first two rows of (\ref{Grelation-F&G}) is equal to LHS of 
(\ref{eq-symmetridReduction-0<=a,b,c,d<=n}). The  RHS of (\ref{eq-symmetridReduction-0<=a,b,c,d<=n}) is obtained in the same way:
\begin{eqnarray*}
&&\sum_{\mu=0}^{n+m} \sum_{\nu=0}^{n+m} \sum_{\lambda=0}^{n+m} \sum_{\rho=0}^{n+m}F_{abc\lambda}\tauG^{\lambda\mu}F_{\mu d \nu \rho}\tauG^{\nu\rho}\\
&=&F_{abce}\tauG^{ei}F_{idjk}\tauG^{jk}+2s F_{abcs}F_{sdjk}\tauG^{jk}-m F_{abce}\tauG^{ei}F_{ids}+2s F_{abce}\tauG^{ei}F_{idss}\nn\\
&&+\big(2(-1)^n m+4\big)s F_{abcs}F_{dss}+4s^2 F_{abcs}F_{dsss},
\end{eqnarray*}
\begin{eqnarray*}
&&\sum_{\mu=0}^{n+m} \sum_{\nu=0}^{n+m} \sum_{\lambda=0}^{n+m} \sum_{\rho=0}^{n+m}F_{abcd\lambda}\tauG^{\lambda\mu}F_{\mu\nu\rho}\tauG^{\nu\rho}\\
&=&F_{abcde}\tauG^{ei}F_{ijk}\tauG^{jk}+2s F_{abcds}F_{sjk}\tauG^{jk}-m F_{abcde}\tauG^{ei}F_{is}+2s F_{abcde}\tauG^{ei}F_{iss}\nn\\
&&+\big(2(-1)^n m+4\big)s F_{abcds}F_{ss}+4s^2 F_{abcds}F_{sss},
\end{eqnarray*}
\begin{eqnarray*}
&&\sum_{\mu=0}^{n+m} \sum_{\nu=0}^{n+m} \sum_{\lambda=0}^{n+m} \sum_{\rho=0}^{n+m}F_{ab\nu \lambda}\tauG^{\lambda\mu}F_{\mu cd\rho}\tauG^{\nu\rho}\\
&=&F_{abij}F_{cdkl}\tauG^{jk}\tauG^{il}+4s F_{abis}F_{cdsj}\tauG^{ij}\nn\\
&& +(-1)^n m F_{abs}F_{cds}+2s F_{abs}F_{cdss}+2s F_{abss}F_{cds}+4s^2 F_{abss}F_{cdss}.
\end{eqnarray*}
\end{proof}

Fix $N\in \mathbb{Z}_{\geq 0}$ ($0\leq N\leq \frac{m}{2}$ if $n$ is odd), and $I=(i_1,\dots,i_l)$, where $0\leq i_1,\dots,i_l\leq n$. Let
\begin{equation}\label{eq-def-y-yi-yij} 
y=G_{i_1,\dots,i_l}^{(N)}(0),\
y_{i}=G_{i,i_1,\dots,i_l}^{(N)}(0),\ 
y_{i,j}=G_{i,j,i_1,\dots,i_l}^{(N)}(0).
\end{equation}
Note that by the convention (\ref{eq-convention-subscript-outOfRange}), any subscripts $i,j\geq 0$ are allowed.
For $1\leq a,b,c,d\leq n$ let
\begin{eqnarray}\label{eq-def-zabcd}
&&w_{a,b,c,d}\nn\\
&=& \sum_{\mathscr{P}(a,b,c,d)}\Big(\frac{N}{2}F_{ab}^{(1)}(0)y_{c+d}
-\frac{1}{24}F_{a+b,c,d,f}^{(0)}(0)\tauG^{fi}y_{i}
+\frac{1}{12}F_{abce}^{(0)}(0)\tauG^{ei}y_{i+d}
-\frac{1}{3}N F_{a+b,c}^{(1)}(0) y_d\nn\\
&&+\big(\frac{N(3N-2)}{6}F_{ab}^{(1)}(0)F_{cd}^{(1)}(0)
-\frac{N}{12} F_{a+b,c,d}^{(1)}(0)
+\frac{N}{6}F_{abce}^{(0)}(0)\tauG^{ei}F_{id}^{(1)}(0)\big)y\Big).
\end{eqnarray}

\begin{proposition}\label{prop-reconstruction-genus1-fromLength1}
Using the Euler field and Getzler's relation (\ref{Grelation-F&G}), one can reconstruct $G^{(0)}$ from $G^{(0)}_i(0)$ for $1\leq i\leq n$ and genus zero GW invariants. More generally, for $k\geq 1$, one can reconstruct $G^{(k)}$ from  $G^{(k)}(0)$ and $G^{(k)}_j(0)$ for $1\leq j\leq n$, and $G^{(j)}$ for $0\leq j\leq k-1$, and genus zero GW invariants.
\end{proposition}
\begin{proof}
For $0\leq a,b,c,d\leq n$, we compute the coefficients of $\tau^{i_1}\cdots \tau^{i_l} s^N$ in (\ref{eq-symmetridReduction-0<=a,b,c,d<=n}). Using (\ref{eq-giventalRelation}) we get
\begin{multline*}
 \sum_{\mathscr{P}(a,b,c,d)}(3y_{a+b,c+d}-4y_{a,b+c+d})+24w_{a,b,c,d}\\
 +N\sum_{\mathscr{P}(a,b,c,d)}\big(-2F_{ab}^{(1)}(0)F_{c}^{(1)}(0)y_{d}
+4F_{abc}^{(1)}(0)F_{d}^{(1)}(0) y\big)\
 \Nlsim 0.
\end{multline*}
Suppose $1\leq a,b,c,d\leq n$. Then by (\ref{eq-F(1)i(0)-Kronecker}) we obtain
\begin{eqnarray}\label{eq-recursionGenusOne-1}
 (y_{a+b,c+d}+y_{a+c,b+d}+y_{a+d,b+c})\nn\\
-(y_{a,b+c+d}+y_{b,a+c+d}+y_{c,a+b+d}+y_{d,a+b+c})
+w_{a,b,c,d}&\Nlsim& 0.
\end{eqnarray}
Taking $c=d=1$ in (\ref{eq-recursionGenusOne-1}), we get
\begin{equation}\label{eq-recursionGenusOne-1-c=d=1}
      y_{2,a+b}+2y_{a+1,b+1}-y_{a,b+2}-y_{b,a+2}-2y_{1,a+b+1}+w_{a,b,1,1}\Nlsim0.
\end{equation}
So
\begin{equation*}
      (y_{a+1,b+1}-y_{a,b+2})-(y_{a+2,b}-y_{a+1,b+1})\Nlsim 2y_{1,a+b+1}-y_{2,a+b}-w_{a,b,1,1}.
\end{equation*}
Let $k=a+b$. Then for $1\leq i\leq k-1$,
\begin{equation*}
      (y_{i+1,k-i+1}-y_{i,k-i+2})-(y_{i+2,k-i}-y_{i+1,k-i+1})\Nlsim2y_{1,k+1}-y_{2,k}-w_{i,k-i,1,1}.
\end{equation*}
Summing over $i=1,\dots,k-1$, we get
\begin{equation}\label{eq-recursionGenusOne-y2k}
      y_{2,k}\Nlsim \frac{2k}{k+1}y_{1,k+1}-\frac{1}{k+1}\sum_{i=1}^{k-1}w_{i,k-i,1,1}.
\end{equation}
By Lemma \ref{lem-EulerField-genus1}, the Euler vector field yields
\begin{equation}\label{eq-recursionGenusOne-y1i}
      y_{1,i}\Nlsim y_i.
\end{equation}
So for all $1\leq i,j\leq n$ we have $y_{i,j}\Nlsim$ a linear combination of $y_1,\dots,y_n$ and $y$, with coefficients the genus 0 GW invariants of $X$. This holds also for $i$ or $j=0$ by the fundamental class axiom. Hence the conclusion follows from  induction on $l$ and $N$.
\end{proof}

\subsection{\texorpdfstring{$0\leq a,b,c\leq n$ and $n+1\leq d\leq n+m$}{0<=a,b,c<=n,n+1<=d<=n+m}}\label{sec:0<=a,b,c<=n,n+1<=d<=n+m}
\begin{lemma}\label{lem-symmetricReduction-0<=a,b,c<=n,n+1<=d<=n+m}
Suppose that $n$ is even (resp. $n$ is odd).
Suppose $0\leq a,b,c\leq n$. Then
\begin{eqnarray}\label{eq-symmetricReduction-0<=a,b,c<=n,n+1<=d<=n+m}
&&  12\sum_{\mathscr{P}(a,b,c)}F_{abe}\tauG^{ei}F_{scf}\tauG^{fj}G_{ij}\nn\\
&&+\sum_{\mathscr{P}(a,b,c)}(12F_{abf}F_{sc}-4F_{abe}\tauG^{ei}F_{cif}+16sF_{sab}F_{scf}
+24sF_{abf}F_{ssc})\tauG^{fj}G_{sj}\nn\\
&&-12\sum_{\mathscr{P}(a,b,c)}(F_{sbc}F_{sf}+F_{sc}F_{sbf}+F_{sb}F_{scf}
-F_{bcfe}\tauG^{ei}F_{si}-2sF_{sbcf}F_{ss})\tauG^{fj}G_{aj}\nn\\
&&-12\sum_{\mathscr{P}(a,b,c)}(F_{sb}F_{sc}+2sF_{sb}F_{ssc}+2sF_{ssb}F_{sc}+4s^2 F_{ssb}F_{ssc}+2s F_{sbe}\tauG^{ef} F_{scf})G_{sa}\nn\\
&&+12\big(\sum_{\mathscr{P}(a,b,c)}(F_{sbc}F_{saf}-F_{bcfe}\tauG^{ei}F_{sia}-2sF_{sbcf}F_{ssa})\nn\\
&&-2F_{abcfe}\tauG^{ei}F_{si}+2F_{abce}\tauG^{ei} F_{sif}+4F_{sabc}F_{sf}
-4s F_{sabcf}F_{ss}+4s F_{sabc}F_{ssf}\big)\tauG^{fj}G_j\nn\\
&&+24\big(s\sum_{\mathscr{P}(a,b,c)}F_{sbc}F_{ssa}+ 2F_{abce}\tauG^{ei}F_{si}
+2sF_{abce}\tauG^{ei}F_{ssi}+6sF_{sabc}F_{ss}+4s^2 F_{sabc}F_{sss}\big)G_s\nn\\
&&+48s(s\sum_{\mathscr{P}(a,b,c)}F_{sbc}F_{ssa}+F_{abce}\tauG^{ei}F_{si}+2sF_{sabc}F_{ss})G_{ss}\nn\\
&=&-\Big(F_{abce}\tauG^{ei}F_{sifj}\tauG^{fj}
+((-1)^n m+2)F_{abce}\tauG^{ei}F_{ssi}+2sF_{abce}\tauG^{ei}F_{sssi}\nn\\
&&+\big((2(-1)^n m+4)F_{ss}+2((-1)^n m+6)sF_{sss}+4s^2F_{ssss}+2F_{sfj}\tauG^{fj}+2sF_{ssfj}\tauG^{fj}\big)F_{sabc}\nn\\
&&+F_{sabce}\tauG^{ei}F_{ifj}\tauG^{fj}+2sF_{ssabc}F_{sfj}\tauG^{fj}+(-1)^n mF_{sabce}\tauG^{ei}F_{si}+2sF_{sabce}\tauG^{ei}F_{ssi}\nn\\
&&+2((-1)^n m+2)sF_{ssabc}F_{ss}+4s^2F_{ssabc}F_{sss}\nn\\
&&+\frac{1}{2}\sum_{\mathscr{P}(a,b,c)}\big(F_{sabe}\tauG^{ei}F_{cifj}\tauG^{fj}
+F_{sab}F_{scfj}\tauG^{fj}+2sF_{ssab}F_{scfj}\tauG^{fj}
+((-1)^n m-4)F_{sabe}\tauG^{ei}F_{sci}\nn\\
&&-6sF_{sabe}\tauG^{ei}F_{ssci}-((-1)^n m+2)F_{sab}F_{ssc}-2sF_{sab}F_{sssc}+2((-1)^n m-4)sF_{ssab}F_{ssc}\nn\\
&&-4s^2F_{ssab}F_{sssc}-2F_{abef}\tauG^{ei}\tauG^{fj}F_{scij}\big)\Big).
\end{eqnarray}
(resp. $\mod s^{\frac{m}{2}}$ if $n$ is odd).
\end{lemma}
The proof is given in Appendix \ref{sec:proof-lem-symmetricReduction-0<=a,b,c<=n,n+1<=d<=n+m}.

\begin{proposition}\label{prop-elimination-ambient-SOrderAtLeastOne}
With the notations defined in (\ref{eq-def-y-yi-yij}), we have
\begin{enumerate}
      \item[(i)] if $\mathsf{a}_X=n-1$, then for $2\leq k\leq n$,
      \begin{equation}\label{eq-elimination-ambient-SOrderAtLeastOne-fanoIndex=n-1}
            y_{k}\Nlsim 2Ny\times \begin{cases}
            \mathsf{c}(n,\frac{k-1}{\mathsf{a}_X},\mathbf{d})\mathsf{b}(\mathbf{d})^{\frac{k-1}{\mathsf{a}_X}}, & \mbox{if}\  \frac{k-1}{\mathsf{a}_X}\in \mathbb{Z}_{\geq 0};\\
            0, &  \mbox{otherwise}.
            \end{cases}
      \end{equation}
      \item[(ii)]  if $\mathsf{a}_X<n-1$, then for $3\leq k\leq n$,
      \begin{equation}\label{eq-elimination-ambient-SOrderAtLeastOne-fanoIndex<n-1}
            y_{k}\Nlsim 2Ny\times \begin{cases}
            \mathsf{c}(n,\frac{k-1}{\mathsf{a}_X},\mathbf{d})\mathsf{b}(\mathbf{d})^{\frac{k-1}{\mathsf{a}_X}}, & \mbox{if}\  \frac{k-1}{\mathsf{a}_X}\in \mathbb{Z}_{\geq 0};\\
            0, &  \mbox{otherwise}.
            \end{cases}
      \end{equation}      
\end{enumerate}
\end{proposition}
\begin{proof}
For $0\leq a,b,c\leq n$, computing the coefficients of $\tau^{i_1}\cdots \tau^{i_l} s^N$ of (\ref{eq-symmetricReduction-0<=a,b,c<=n,n+1<=d<=n+m}) yields
\begin{eqnarray*}
&3F_{c}^{(1)}(0)y_{a+b}+3F_{b}^{(1)}(0)y_{a+c}+3F_{a}^{(1)}(0)y_{b+c}-3y_{a+b+c}&\\
&- 3\big(F_{a}^{(1)}(0)F_{b}^{(1)}(0)y_{c}+ F_{a}^{(1)}(0)F_{c}^{(1)}(0)y_{b}+ F_{b}^{(1)}(0)F_{c}^{(1)}(0)y_{a}\big)
+6NF_{abce}^{(0)}(0)\tauG^{ei}F_{i}^{(1)}(0)y &\Nlsim 0.
\end{eqnarray*}
Suppose $1\leq a,b,c\leq n$. Then by (\ref{eq-F(1)i(0)-Kronecker}) we get
\begin{eqnarray}\label{eq-recursionGenusOne-3}
      -3 y_{a+b+c}+6NF_{abce}^{(0)}(0)\tauG^{e0}y &\Nlsim& 0.
\end{eqnarray}
Applying (\ref{eq-qp4}), we obtain (\ref{eq-elimination-ambient-SOrderAtLeastOne-fanoIndex<n-1}), and (\ref{eq-elimination-ambient-SOrderAtLeastOne-fanoIndex=n-1}) in the range $3\leq k\leq n$. When $\mathsf{a}_X=n-1$, by the definition (\ref{eq-def-y-yi-yij}) of $y_k$ and the convention (\ref{eq-convention-subscript-outOfRange}), we have $y_{n+1}= \mathbf{d}^{\mathbf{d}}y_2$. Hence taking $a=b=1$ and $c=n-1$ in (\ref{eq-recursionGenusOne-3}) yields (\ref{eq-elimination-ambient-SOrderAtLeastOne-fanoIndex=n-1}).
\end{proof}

\subsection{\texorpdfstring{$n+1\leq a,b,c,d\leq n+m$}{n+1<=a,b,c,d<=n+m}}\label{sec:n+1<=a,b,c,d<=n+m}
\begin{lemma}\label{lem-symReduction-n+1<=a,b,c,d<=n+m}
Suppose $(n,\mathbf{d})\neq \big(3,(2,2)\big)$. Then
\begin{eqnarray}\label{eq-symReduction-n+1<=a,b,c,d<=n+m}
&&3(F_{se}F_{sf}\tauG^{ei}\tauG^{fj}G_{ij}+4sF_{ss}F_{se}\tauG^{ei}G_{si}\nn\\
&&+2F_{ss}F_{sf}\tauG^{fj}G_j+6sF_{ss}F_{ss}G_s+4s^2F_{ss}F_{ss}G_{ss})\nn\\
&=&-\frac{1}{2}F_{ss}F_{sfj}\tauG^{fj}-\frac{1}{8}F_{sse}\tauG^{ei}F_{ifj}\tauG^{fj}
-\frac{1}{4}sF_{sss}F_{sfj}\tauG^{fj}+\frac{3}{4}sF_{sse}\tauG^{ei}F_{ssi}\nn\\
&&+\frac{1}{4}F_{sef}\tauG^{ei}\tauG^{fj}F_{sij}
-(\frac{(-1)^n m}{8}-1)F_{ss}F_{ss}-\frac{1}{2}sF_{sss}F_{ss}
+\frac{1}{2}s^2F_{sss}F_{sss}.
\end{eqnarray}
\end{lemma}
\begin{proof}
We first treat the case that the  dimension $n$ is even. We take $a=b=c=d$ in (\ref{Grelation-F&G}). A direct computation gives
\begin{eqnarray}\label{eq-symReduction-a=b=c=d-LHS}
&&\sum_{\mu=0}^{n+m} \sum_{\nu=0}^{n+m} \sum_{\lambda=0}^{n+m} \sum_{\rho=0}^{n+m}\big(
3F_{aa\lambda}\tauG^{\lambda\mu}G_{\mu\nu}\tauG^{\nu\rho}F_{\rho aa}-4F_{aa\lambda}\tauG^{\lambda\mu}F_{\mu a\rho}\tauG^{\rho\nu}G_{\nu a}\nn\\
&&-F_{aa\lambda}\tauG^{\lambda\mu}F_{\mu aa\rho}\tauG^{\rho\nu}G_{\nu} +2F_{aaa\lambda}\tauG^{\lambda\mu}F_{\mu a\rho}\tauG^{\rho\nu}G_{\nu}\big)\nn\\
&=&  3F_{se}F_{sf}\tauG^{ei}\tauG^{fj}G_{ij}+12sF_{ss}F_{se}\tauG^{ei}G_{si}\nn\\
&&+(-F_{se}\tauG^{ei}F_{sif}-2sF_{ss}F_{ssf}+6F_{ss}F_{sf})\tauG^{fj}G_j\nn\\
&&+(8sF_{ss}F_{ss}-4F_{se}\tauG^{ei}F_{si}-2s F_{se}\tauG^{ei}F_{ssi}-4s^2 F_{ss}F_{sss})G_s\nn\\
&&+12s^2F_{ss}F_{ss}G_{ss}+(\tau^a)^2(\dots)+(\tau^a)^4(\dots),
\end{eqnarray}
where $(\dots)$ means functions of $\tau^0,\dots,\tau^n$ and $s$. From (\ref{eq-wdvv230}) and (\ref{eq-wdvv240}) we have
\begin{eqnarray*}
 -F_{se}\tauG^{ei}F_{sif}-2sF_{ss}F_{ssf}+6F_{ss}F_{sf}= 6F_{ss}F_{sf},
\end{eqnarray*}
and
\begin{eqnarray*}
8sF_{ss}F_{ss}-4F_{se}\tauG^{ei}F_{si}-2s F_{se}\tauG^{ei}F_{ssi}-4s^2 F_{ss}F_{sss}
= 18sF_{ss}F_{ss}.
\end{eqnarray*}
It follows that the sum of the terms without common factor $(\tau^a)^2$ or $(\tau^a)^4$ in  (\ref{eq-symReduction-a=b=c=d-LHS}) is equal to LHS of 
(\ref{eq-symReduction-n+1<=a,b,c,d<=n+m}). On the other hand, 
\begin{eqnarray}\label{eq-symReduction-a=b=c=d-RHS}
&&\sum_{\mu=0}^{n+m} \sum_{\nu=0}^{n+m} \sum_{\lambda=0}^{n+m} \sum_{\rho=0}^{n+m}\big(
\frac{1}{6}F_{abc\lambda}\tauG^{\lambda\mu}F_{\mu d \nu \rho}\tauG^{\nu\rho}+\frac{1}{24}F_{abcd\lambda}\tauG^{\lambda\mu}F_{\mu\nu\rho}\tauG^{\nu\rho}
-\frac{1}{4}F_{ab\nu \lambda}\tauG^{\lambda\mu}F_{\mu cd\rho}\tauG^{\nu\rho}\big)\nn\\
&=&\frac{1}{2}F_{ss}F_{sfj}\tauG^{fj}+\frac{1}{8}F_{sse}\tauG^{ei}F_{ifj}\tauG^{fj}
+\frac{1}{4}sF_{sss}F_{sfj}\tauG^{fj}-\frac{3}{4}sF_{sse}\tauG^{ei}F_{ssi}\nn\\
&&-\frac{1}{4}F_{sef}\tauG^{ei}\tauG^{fj}F_{sij}+\frac{m}{8}F_{sse}\tauG^{ei}F_{si}
+\frac{m-4}{4}F_{ss}F_{ss}+\frac{m+2}{4}sF_{sss}F_{ss}
-\frac{1}{2}s^2F_{sss}F_{sss}\nn\\
&&+(\tau^a)^2(\dots)+(\tau^a)^4(\dots).
\end{eqnarray}

By (\ref{eq-wdvv240}),
\begin{equation*}
            F_{se}\tauG^{ef}F_{ssf}+F_{ss}F_{ss}+2sF_{ss}F_{sss}=0.
\end{equation*}
It follows that the sum of the terms without common factor $(\tau^a)^2$ or $(\tau^a)^4$ in  (\ref{eq-symReduction-a=b=c=d-RHS}) gives RHS of (\ref{eq-symReduction-n+1<=a,b,c,d<=n+m}). Recall that (\cite[Corollary A.3]{Hu15}) $m=\mathrm{rank} H^*_{\mathrm{prim}}(X)\geq 3$. So (\ref{eq-symReduction-a=b=c=d-LHS})=$-$(\ref{eq-symReduction-a=b=c=d-RHS}) implies that the sums of the terms without common factor $(\tau^a)^2$ or $(\tau^a)^4$ of both sides are equal. Hence we have (\ref{eq-symReduction-n+1<=a,b,c,d<=n+m}).

If $n$ is odd, to obtain a non-trivial equation, we take $n+1\leq a<c\leq n+\frac{m}{2}$, $b=a'=a+\frac{m}{2}$ and $d=c'=c+\frac{m}{2}$ in (\ref{Grelation-F&G}). 
We adopt the conventions introduced in the odd dimensional cases of the proof of Lemma \ref{lem-symmetridReduction-0<=a,b,c,d<=n}. Then a calculation of Grassmann variables yields 
\begin{eqnarray}\label{eq-symReduction-a=b=c=d-LHS-oddDim}
&&\sum_{\mathscr{P}(a,a',c,c')}(\pm)
\sum_{\mu=0}^{n+m} \sum_{\nu=0}^{n+m} \sum_{\lambda=0}^{n+m} \sum_{\rho=0}^{n+m}\big(
3F_{aa'\lambda}\tauG^{\lambda\mu}G_{\mu\nu}\tauG^{\nu\rho}F_{\rho cc'}-4F_{aa'\lambda}\tauG^{\lambda\mu}F_{\mu c\rho}\tauG^{\rho\nu}G_{\nu c'}\nn\\
&&-F_{aa'\lambda}\tauG^{\lambda\mu}F_{\mu cc'\rho}\tauG^{\rho\nu}G_{\nu} +2F_{aa'c\lambda}\tauG^{\lambda\mu}F_{\mu c'\rho}\tauG^{\rho\nu}G_{\nu}\big)\nn\\
&=&8\times\big(3F_{se}F_{sf}\tauG^{ei}\tauG^{fj}G_{ij}+12sF_{ss}F_{se}\tauG^{ei}G_{si}
+(-F_{se}\tauG^{ei}F_{sif}-2sF_{ss}F_{ssf}+6F_{ss}F_{sf})\tauG^{fj}G_j\nn\\
&&+(8sF_{ss}F_{ss}-4F_{se}\tauG^{ei}F_{si}-2s F_{se}\tauG^{ei}F_{ssi}-4s^2 F_{ss}F_{sss})G_s+12s^2F_{ss}F_{ss}G_{ss}\big)\nn\\
&&+\tau^a \tau^{a'}(\dots)+\tau^c\tau^{c'}(\dots)(\dots)+\tau^a \tau^{a'}\tau^c \tau^{c'}(\dots),
\end{eqnarray}
and
\begin{eqnarray}\label{eq-symReduction-a=b=c=d-RHS-oddDim}
&&\sum_{\mathscr{P}(a,a',c,c')}(\pm)
\sum_{\mu=0}^{n+m} \sum_{\nu=0}^{n+m} \sum_{\lambda=0}^{n+m} \sum_{\rho=0}^{n+m}\big(
\frac{1}{6}F_{aa'c\lambda}\tauG^{\lambda\mu}F_{\mu c' \nu \rho}\tauG^{\nu\rho}\nn\\
&&+\frac{1}{24}F_{aa'cc'\lambda}\tauG^{\lambda\mu}F_{\mu\nu\rho}\tauG^{\nu\rho}
-\frac{1}{4}F_{aa'\nu \lambda}\tauG^{\lambda\mu}F_{\mu cc'\rho}\tauG^{\nu\rho}\big)\nn\\
&=&8\times\big(\frac{1}{2}F_{ss}F_{sfj}\tauG^{fj}+\frac{1}{8}F_{sse}\tauG^{ei}F_{ifj}\tauG^{fj}
+\frac{1}{4}sF_{sss}F_{sfj}\tauG^{fj}-\frac{3}{4}sF_{sse}\tauG^{ei}F_{ssi}\nn\\
&&-\frac{1}{4}F_{sef}\tauG^{ei}\tauG^{fj}F_{sij}+\frac{m}{8}F_{sse}\tauG^{ei}F_{si}
+\frac{m-4}{4}F_{ss}F_{ss}+\frac{m+2}{4}sF_{sss}F_{ss}
-\frac{1}{2}s^2F_{sss}F_{sss}\big)\nn\\
&&+\tau^a \tau^{a'}(\dots)+\tau^c\tau^{c'}(\dots)(\dots)+\tau^a \tau^{a'}\tau^c \tau^{c'}(\dots),
\end{eqnarray}
where $(\dots)$ means functions of $\tau^0,\dots,\tau^n$ and $s$.
By \cite[Theorem A.2]{Hu15}, if $n\geq 3$ or $\mathbf{d}\neq (2,2)$, we have $m=\mathrm{rank} H^*_{\mathrm{prim}}(X)\geq 6$. So (\ref{eq-symReduction-a=b=c=d-LHS-oddDim})=$-$(\ref{eq-symReduction-a=b=c=d-RHS-oddDim}) yields (\ref{eq-symReduction-n+1<=a,b,c,d<=n+m}).
\end{proof}

\begin{corollary}
Let $y,y_i,y_{i,j}$ be defined as in (\ref{eq-def-y-yi-yij}). Then
\begin{eqnarray}\label{eq-recursionEq00}
 F_{e}^{(1)}(0)F_{f}^{(1)}(0)\tauG^{ei}\tauG^{fj}y_{i,j}+(4N+2)F^{(2)}(0)F_{e}^{(1)}(0)\tauG^{ei}y_{i}
+(4N^2+2N)F^{(2)}(0)^2 y&\Nlsim& 0.
\end{eqnarray}
\end{corollary}
\begin{proof}
If $(n,\mathbf{d})\neq \big(3,(2,2)\big)$, by computing the coefficients of $\tau^{i_1}\cdots \tau^{i_l} s^N$ of (\ref{eq-symReduction-n+1<=a,b,c,d<=n+m}) we get (\ref{eq-recursionEq00}). Now suppose $(n,\mathbf{d})=\big(3,(2,2)\big)$. Then by Riemann-Roch one can compute $m=\mathrm{rank} H^*_{\mathrm{prim}}(X)=4$. With the notations in the proof of Lemma \ref{lem-symReduction-n+1<=a,b,c,d<=n+m}, we have
\begin{equation*}
      s=-\tau^a \tau^{a'}-\tau^{c}\tau^{c'},\ s^2=2\tau^a \tau^{a'}\tau^c \tau^{c'}.
\end{equation*}
Thus (\ref{eq-symReduction-a=b=c=d-LHS-oddDim})=$-$(\ref{eq-symReduction-a=b=c=d-RHS-oddDim}) yields a more complicated identity involving $G$ and $F$. But one can check that in the abbreviated terms $\tau^a \tau^{a'}(\dots)$ and $\tau^c \tau^{c'}(\dots)$ in (\ref{eq-symReduction-a=b=c=d-LHS-oddDim}) and (\ref{eq-symReduction-a=b=c=d-RHS-oddDim}), there is at least one additional factor $s$ for $G_s$ and $G_{si}$ ($0\leq i\leq n$), and an additional factor $s^2$ for $G_{ss}$. Similarly, in the abbreviated term $\tau^a \tau^{a'}\tau^c \tau^{c'}(\dots)$ in (\ref{eq-symReduction-a=b=c=d-LHS-oddDim}) and (\ref{eq-symReduction-a=b=c=d-RHS-oddDim}), there is  one additional factor $s$ for $G_{ss}$. Hence computing the coefficients of $\tau^{i_1}\cdots \tau^{i_l} s^N$ of (\ref{eq-symReduction-a=b=c=d-LHS-oddDim})=$-$(\ref{eq-symReduction-a=b=c=d-RHS-oddDim}), we get again (\ref{eq-recursionEq00}).
\end{proof}

\subsection{Cubic hypersurfaces and odd dimensional complete intersections of two quadrics}
We use the results in Section \ref{sec:0<=a,b,c,d<=n}-\ref{sec:n+1<=a,b,c,d<=n+m} to show Theorem \ref{thm-reconstruction-cubic-(2,2)-intro}. First we recall that (\cite[Theorem 1.8 and 1.10]{Hu15}) for $\mathbf{d}=(3)$ or $(2,2)$,
\begin{equation}\label{eq-F(2)(0)-cubic-(2,2)}
      F^{(2)}(0)=1.
\end{equation}
\begin{lemma}\label{lem-CubicHypersurfaces-(2,2)-sumZ}
Let $y,y_i,y_{i,j}$ be defined as in (\ref{eq-def-y-yi-yij}), where $I=\emptyset$. Then
\begin{eqnarray}\label{eq-CubicHypersurfaces-(2,2)-sumZ}
&& \frac{\mathbf{d}^{\mathbf{d}}}{2n}\sum_{i=1}^{n-2}w_{i,n-1-i,1,1}+\frac{1}{2}w_{n-1,n-1,1,1}\nn\\
&=&(-4\mathbf{d}!N+\mathbf{d}^{\mathbf{d}}-2\mathbf{d}!)y_n
+\big(\frac{(4n-2)\mathbf{d}!}{n}N+2 \mathbf{d}!- \mathbf{d}^{\mathbf{d}}+\frac{\mathbf{d}^{\mathbf{d}}-\mathbf{d}!}{n}\big)\mathbf{d}^{\mathbf{d}}y_1\nn\\
&&+\big(4\mathbf{d}!N+2 \mathbf{d}!-\frac{2 \mathbf{d}^{\mathbf{d}}}{n})N\mathbf{d}!y.
\end{eqnarray}
\end{lemma}
The proof is given in the appendix \ref{sec:proof-eq-CubicHypersurfaces-(2,2)-sumZ}.

\begin{proposition}\label{prop-reconstruction-cubic-(2,2)}
Suppose (i) $\mathbf{d}=(3)$ and $n\geq 3$ and $n\neq 4$ or (ii) $\mathbf{d}=(2,2)$ and $n\geq 3$ and $n$ is odd. Then the genus 1 Gromov-Witten invariants of $X$ can be reconstructed from the genus 0 invariants, and the genus 1 invariants with 1 marked point.
\end{proposition}
\begin{proof}
By (\ref{eq-pairing2}), (\ref{eq-F(1)i(0)-Kronecker}) and (\ref{eq-F(2)(0)-cubic-(2,2)}), (\ref{eq-recursionEq00}) reads
\begin{equation}\label{eq-recursionEq00-cubic-(2,2)}
 \frac{1}{(\prod_{i=1}^r d_i)^2}\big(y_{n,n}-2 \mathbf{d}^{\mathbf{d}}y_{1,n}+(\mathbf{d}^{\mathbf{d}})^2 y_{1,1}\big)
+\frac{4N+2}{\prod_{i=1}^r d_i}(y_n- \mathbf{d}^{\mathbf{d}}y_1)+(4N^2+2N)y\Nsim 0.
\end{equation}
Taking $a=b=n-1$ in (\ref{eq-recursionGenusOne-1-c=d=1}) we get
\begin{equation}\label{eq-recursionGenusOne-1-c=d=1-cubic-(2,2)}
      \mathbf{d}^{\mathbf{d}}y_{2,n-1}+2y_{n,n}-2\mathbf{d}^{\mathbf{d}}y_{2,n-1}-2 \mathbf{d}^{\mathbf{d}}y_{1,n}+w_{n-1,n-1,1,1}\Nsim 0.
\end{equation}
Then by (\ref{eq-recursionGenusOne-y2k}) we get
\begin{eqnarray}\label{eq-ynn-cubic-(2,2)}
y_{n,n}&\Nsim &   \frac{(2n-1)\mathbf{d}^{\mathbf{d}}}{n}y_{1,n}
-\frac{\mathbf{d}^{\mathbf{d}}}{2n}\sum_{i=1}^{n-2}w_{i,n-1-i,1,1}-\frac{1}{2}w_{n-1,n-1,1,1}.
\end{eqnarray}
By (\ref{eq-EulerField-genus1}) or directly using the divisor equation, we have
\begin{eqnarray}\label{eq-y1n-y11-y1-cubic-(2,2)}
y_{1,n}=(\frac{N(n-2)}{n-1}+1)y_n-(\mathbf{d}^{\mathbf{d}}- \mathbf{d}!)y_1,\
y_{1,1}= \frac{N(n-2)}{n-1} y_1,\
y_{1}= \frac{N(n-2)}{n-1} y.
\end{eqnarray}
Recall \cite[(203)]{Hu15},
\[
\mathsf{c}(n,1,\mathbf{d})=\frac{\mathbf{d}!}{\mathbf{d}^{\mathbf{d}}}.
\]
So by Proposition \ref{prop-elimination-ambient-SOrderAtLeastOne} (i),
\begin{equation}\label{eq-yn-cubic-(2,2)}
      y_n\Nsim 2N \mathsf{c}(n,1,\mathbf{d})\cdot \mathbf{d}^{\mathbf{d}} y=2 \mathbf{d}! N y.
\end{equation}
Applying (\ref{eq-CubicHypersurfaces-(2,2)-sumZ}) to (\ref{eq-ynn-cubic-(2,2)}), and then substituting (\ref{eq-ynn-cubic-(2,2)}), (\ref{eq-y1n-y11-y1-cubic-(2,2)}) and (\ref{eq-yn-cubic-(2,2)}) into (\ref{eq-recursionEq00-cubic-(2,2)}), after simplification we get, (i) if $\mathbf{d}=3$,
\begin{equation}
      \frac{9 (n-4) N (n N+n-4 N-1)}{(n-1)^2} y\Nsim 0.
\end{equation}
If $n\geq 5$, $n N+n-4 N-1>0$. So $y\sim 0$ for $n\geq 5$. If $n=3$, $\frac{9N(N-2)}{4}y\sim 0$. But when $N=2$, $y=0$ for the degree reason.

(ii) if $\mathbf{d}=(2,2)$,  
\begin{equation}
   \frac{16 N^2}{(n-1)^2}y\Nsim 0.
\end{equation}
So we are done.
\end{proof}

\section{Reduced genus 1 Gromov-Witten invariants of Fano complete intersections}\label{sec:ReducedGenus1GW-Fano}
Let $X$ be a smooth Fano hypersurface of $\mathbb{P}^{n-1}$ of multidegree $\mathbf{d}=(d_1,\dots,d_r)$. Without loss of generality, we assume that $d_i\geq 2$ for $1\leq i\leq r$. We denote the Fano index of $X$ by
\begin{equation}
      \nu_{\mathbf{d}}=n-|\mathbf{d}|=n-\sum_{i=1}^r d_i.
\end{equation}
 We are going to compute the reduced genus 1 Gromov-Witten invariant
$\langle \sfh_{a}\rangle_{1}^{0}$
of $X$, defined by Zinger \cite{Zin09a}. We follow the strategy of Zinger \cite{Zin09b}. 
By the dimension constraint this invariant vanishes unless
\begin{equation}
      \beta=\frac{a-1}{\nu_{\mathbf{d}}}\in \mathbb{Z}_{\geq 0}.
\end{equation}
Let $\mathcal{M}^0_{1,1}(\mathbb{P}^{n-1},\beta)$ be the open substack of $\Mbar_{1,1}(\mathbb{P}^{n-1},\beta)$ parametrizing genus 1 stable maps with smooth domain curves. Let  $\Mbar^0_{1,1}(\mathbb{P}^{n-1},\beta)$ be the Zariski closure of $\mathcal{M}^0_{1,1}(\mathbb{P}^{n-1},\beta)$.
Let $\pi:\mathfrak{U}\rightarrow \Mbar_{1,1}(\mathbb{P}^{n-1},\beta)$ be the universal curve, $\mathrm{ev}:\mathfrak{U}\rightarrow \mathbb{P}^{n-1}$ be the universal morphsim. The pushforward $\pi_* \mathrm{ev}^* \bigoplus_{k=1}^r \mathcal{O}(d_k)$ is locally free on $\mathcal{M}^0_{1,1}(\mathbb{P}^{n-1},\beta)$, but not locally free on $\Mbar_{1,1}^0(\mathbb{P}^{n-1},\beta)$. Nevertheless, there exists a well-defined Euler class of $\pi_* \mathrm{ev}^* \bigoplus_{k=1}^r \mathcal{O}(d_k)$ (\cite{Zin07a}), and we have the hyperplane property (\cite{LZ09}) of reduced genus 1 Gromov-Witten invariants 
\begin{equation}
      \langle \sfh_{a}\rangle_{1}^{0}=\int_{[\Mbar_{1,1}^0(\mathbb{P}^{n-1},\beta)]} \mathrm{ev}^* \sfh_a \cup \mathbf{e}\big(\pi_* \mathrm{ev}^* \bigoplus_{k=1}^r \mathcal{O}(d_k)\big).
\end{equation}

Vakil and Zinger \cite{VZ08} constructed  (see also \cite{HL10} for a simplified construction)  a modular desingularization $\rho:\Mtd_{1,1}^0(\mathbb{P}^{n-1},\beta)\rightarrow\Mbar_{1,1}^0(\mathbb{P}^{n-1},\beta)$, and a locally free sheaf $\mathcal{V}_1$, such that $\rho_* \mathcal{V}_1\cong \pi_* \mathrm{ev}^* \bigoplus_{k=1}^r \mathcal{O}(d_k)$, and then
\begin{equation}
      \langle \sfh_{a}\rangle_{1}^{0}
      =\int_{[\Mtd_{1,1}^0(\mathbb{P}^{n-1},\frac{a-1}{\nu_{\mathbf{d}}})]}\mathrm{ev}_1^* \sfh_a \cup \mathbf{e}(\mathcal{V}_1).
\end{equation}

Then we compute $\langle \sfh_{a}\rangle_{1}^{0}$ via Atiyah-Bott localization theorem. We follow the set-up of \cite{Zin09b}. We equip $\mathbb{P}^{n-1}$ a $\mathbb{T}=(\mathbb{C}^*)^n$-action, such that 
\begin{enumerate}
      \item let $\gamma_{n-1}$ be the induced equivariant tautological line bundle, and $x=c_1(\gamma_{n-1}^*)$, the equivariant first Chern class;
      \item let $P_1,\dots,P_n$ be the $n$ fixed points of $\mathbb{P}^n$, then $H_{\mathbb{T}}^*(P_i)=\mathbb{Q}[\alpha_1,\dots,\alpha_n]$, and $x|_{P_i}=\alpha_i$;
      \item let $\phi_i=\prod_{k\neq i}(x-\alpha_k)\in H_{\mathbb{T}}^*(\mathbb{P}^{n-1})$, then $\mathbf{e}(T \mathbb{P}^{n-1})|_{P_i}=\prod_{k\neq i}(\alpha_i-\alpha_k)=\phi_i|_{P_i}$.
\end{enumerate}

 We linearize $\mathcal{O}(1)$ such that  $e\big(\mathcal{O}(1)\big)=\alpha_i$ at $P_i$, for $1\leq i\leq n$, and linearize $T \mathbb{P}^{n-1}$ such that $e(T \mathbb{P}^{n-1})=\prod_{k\neq i}(\alpha_i- \alpha_k)$ at $P_i$.
 Let 
 \begin{equation}
      \mathcal{F}(\alpha,x,q)=\sum_{\beta=1}^{\infty}q^\beta \big(\mathrm{ev}_{1,\beta*}\mathbf{e}(\mathcal{V}_1)\big)
      \in H^*_{\mathbb{T}}(\mathbb{P}^{n-1})[q].
 \end{equation}

\begin{lemma}\label{lem-reducedGenus1GW-generatingFunc}
We have
\begin{equation}
      \mathcal{F}(\alpha,\alpha_i,q)
       =\sum_{\beta=1}^{\infty}q^\beta 
      \int_{\Mtd_{1,1}^{0}(\mathbb{P}^{n-1},\beta)}\mathbf{e}(\mathcal{V}_1)\mathrm{ev}_1^*\phi_i,
\end{equation}
and 
\begin{equation}
\langle \sfh_{1+\nu_{\mathbf{d}}\beta}\rangle_{1,\beta}^0=
     \sum_{i=1}^n \frac{\alpha_i^{1+\nu_{\mathbf{d}}\beta}\mathrm{Coeff}_{q^\beta}\big\{ \mathcal{F}(\alpha,\alpha_i,q)\big\}}{\prod_{j\neq i}(\alpha_i- \alpha_j)}.
\end{equation}
\end{lemma}
\begin{proof}
From the definition we have
 \begin{equation}
      \mathcal{F}(\alpha,x,q)=\sum_{\beta=1}^{\infty}q^{\beta}\big(
      \mathcal{F}_{ \beta,0}x^{n-2-\nu_{\mathbf{d}}\beta}+
      \sum_{i=1}^{n-2-\nu_{\mathbf{d}}\beta}\mathcal{F}_{\beta,i}(\alpha)x^{n-2-\nu_{\mathbf{d}}\beta-i}\big),
 \end{equation}
 where  $\mathcal{F}_{\beta,i}(\alpha)\in \mathbb{Q}[\alpha_1,\dots,\alpha_n]$ is a symmetric polynomial homogeneous  of degree $i$. 
In particular,
 \begin{equation}
 \mathcal{F}_{\beta,0}=\langle \sfh_{1+\nu_{\mathbf{d}}\beta}\rangle_{1,\beta}^0.
 \end{equation}
By the projection formula we compute the localization of $\mathcal{F}(\alpha,x,q)$ at $P_i$:
\begin{equation}
      \mathcal{F}(\alpha,\alpha_i,q)=\mathcal{F}(\alpha,x,q)|_{P_i}=\sum_{\beta=1}^{\infty}q^\beta 
      \int_{\mathbb{P}^{n-1}}\big(\mathrm{ev}_{1,\beta*}\mathbf{e}(\mathcal{V}_1)\big)\phi_i
      =\sum_{\beta=1}^{\infty}q^\beta 
      \int_{\Mtd_{1,1}^{0}(\mathbb{P}^{n-1},\beta)}\mathbf{e}(\mathcal{V}_1)\mathrm{ev}_1^*\phi_i.
\end{equation}
One notices that $\mathcal{F}(\alpha,x,q)$ is completely determined by $\mathcal{F}(\alpha,\alpha_1,q)$,\dots, $\mathcal{F}(\alpha,\alpha_n,q)$. Moreover
\begin{equation}
      \mathcal{F}(\alpha,\alpha_i,q \alpha_i^{\nu_{\mathbf{d}}} )=\sum_{\beta=1}^{\infty}q^{\beta}\big(
      \mathcal{F}_{ \beta,0}\alpha_i^{n-2}+
      \sum_{i=1}^{n-2-\nu_{\mathbf{d}}\beta}\mathcal{F}_{\beta,i}(\alpha)\alpha_i^{n-2-i}\big).
\end{equation}
Hence
\begin{equation}\label{eq-reducedGenus1Inv-inTermsOfEquivInv}
      \mathcal{F}_{\beta,0}=
     \sum_{i=1}^n \frac{\alpha_i^{1+\nu_{\mathbf{d}}\beta}\mathrm{Coeff}_{q^\beta}\big\{ \mathcal{F}(\alpha,\alpha_i,q)\big\}}{\prod_{j\neq i}(\alpha_i- \alpha_j)}.
\end{equation}
\end{proof}

\begin{remark}\label{rem:comparision-Zinger'sMethod-ours}
Our task is to extract the coefficient of $q^{\beta}\alpha_i^{n-2}$ in $\mathcal{F}(\alpha,\alpha_i,q \alpha_i^{\nu_{\mathbf{d}}} )$ to get  $\mathcal{F}_{\beta,0}$. In \cite{Zin09b} and \cite{Pop12}, the Calaib-Yau case is done by some sophisticated use of symmetric functions. The approach in this paper is simple and direct: we evaluate the equivariant series at $\alpha_i=\zeta^i$ in (\ref{eq-reducedGenus1Inv-inTermsOfEquivInv}), where $\zeta$ is a primitive $n$-th root of unity.
\end{remark}

\subsection{Contributions in the torus localization}
In this section, we recall some equivariant generating functions of genus 0 Gromov-Witten invariants introduced in \cite{PoZ14} and then express the localization contributions in terms of these equivariant generating functions. Proposition \ref{prop-contribution-typeA}, \ref{prop-contribution-typeB}, and Lemma \ref{prop-regularizability-Equivariant} in the following correspond to \cite[Prop. 1.1, 1.2 Lemma 2.3]{Zin09b} respectively, and the proofs are identical and thus omitted. 

Let $\pi:\mathfrak{U}\rightarrow \Mbar_{0,l}(\mathbb{P}^{n-1},\beta)$ be the universal curve, $\mathrm{ev}:\mathfrak{U}\rightarrow \mathbb{P}^{n-1}$ be the universal morphsim, $\varsigma_i:\Mbar_{0,l}(\mathbb{P}^{n-1},\beta)\rightarrow \mathfrak{U}$ the section of $\pi$ corresponding to the $i$-th marked points of the fiber curves, and $\mathrm{ev}_i=\mathrm{ev}\circ \varsigma_i$. Let $\mathcal{V}_0=\pi_* \mathrm{ev}^* \bigoplus_{k=1}^r \mathcal{O}(d_k)$, which is locally free. The restriction of sections to the $i$-th marked points gives rise to a surjective morphism of locally free sheaves $\mathcal{V}_0\rightarrow \mathrm{ev}^* \bigoplus_{k=1}^r \mathcal{O}(d_k)$. When $l=2$, let $\mathcal{V}'_0=\mathrm{ker}(\mathrm{ev}_1)$, and $\mathcal{V}''_0=\mathrm{ker}(\mathrm{ev}_2)$. Then we define 
\begin{equation}
      \mathcal{Z}_i^*(\hbar,q):=\sum_{\beta=1}^{\infty}q^{\beta}\int_{\Mbar_{0,2}(\mathbb{P}^{n-1},\beta)} 
      \frac{\mathbf{e}(\mathcal{V}'_0)}{\hbar-\psi_1}\mathrm{ev}_1^* \phi_i,
\end{equation}
\begin{equation}
      \mathcal{Z}_{ij}^*(\hbar,q):=\hbar^{-1}\sum_{\beta=1}^{\infty}q^{\beta}\int_{\Mbar_{0,2}(\mathbb{P}^{n-1},\beta)}
      \frac{\mathbf{e}(\mathcal{V}'_0)}{\hbar-\psi_1}\mathrm{ev}_1^* \phi_i \mathrm{ev}_2^{*}\phi_j,
\end{equation}
\begin{equation}
      \widetilde{\mathcal{Z}}_{ij}^*(\hbar_1,\hbar_2,q):=\frac{1}{2\hbar_1 \hbar_2}
      \sum_{\beta=1}^{\infty}q^{\beta}\int_{\Mbar_{0,2}(\mathbb{P}^{n-1},\beta)}
      \frac{\mathbf{e}(\mathcal{V}'_0)}{(\hbar_1-\psi_1)(\hbar_2-\psi_2)}\mathrm{ev}_1^* \phi_i \mathrm{ev}_2^{*}\phi_j.
\end{equation}

\begin{equation}
      \mathcal{Z}_p^*(x,\hbar,q):=x^{r+p}+\sum_{\beta=1}^{\infty}q^{\beta} \mathrm{ev}_{1*}\Big(
      \frac{\mathbf{e}(\mathcal{V}''_0)\mathrm{ev}_2^* x^{r+p}}{\hbar-\psi_1}\Big).
\end{equation}

Let $\mathcal{Z}(\hbar,u)\in u\cdot\mathbb{Q}(\alpha_1,\dots,\alpha_n,\hbar)[[u]]$. Following \cite[Definition 2.1]{Zin09b}, we say that $\mathcal{Z}(\hbar,u)$ is \emph{regularizable} (at $\hbar=0$) if there exists $\eta(u)\in u \mathbb{Q}(\alpha_1,\dots,\alpha_n)[[u]]$ and 
$\bar{\mathcal{Z}}(\hbar,u)\in u\mathbb{Q}(\alpha_1,\dots,\alpha_n,\hbar)[[u]]$, such that $\bar{\mathcal{Z}}(\hbar,u)$ is regular at $\hbar=0$ and
\begin{equation}
      1+\mathcal{Z}(\hbar,u)=e^{\eta(u)/\hbar}\big( 1+\bar{\mathcal{Z}}(\hbar,u)\big).
\end{equation}

\begin{lemma}\label{prop-regularizability-Equivariant}
$\mathcal{Z}^*_i(\hbar,q)$ is regularizable.
\end{lemma}
This enables us to define
\begin{equation}
      \eta_i(q):=\mathrm{Res}_{\hbar=0}\Big\{\ln\big(1+\mathcal{Z}^*_i(\hbar,q)\big)\Big\},
\end{equation}
and
\begin{equation}
      \Phi_0(\alpha_i,q):=\mathrm{Res}_{\hbar=0}\big\{\hbar^{-1}e^{-\eta_i(q)/\hbar}\big(1+\mathcal{Z}^*_i(\hbar,q)\big)\big\}
      =\big\{e^{-\eta_i(q)/\hbar}\big(1+\mathcal{Z}^*_i(\hbar,q)\big)\big\}|_{\hbar=0}.
\end{equation}

The fixed loci of $\Mtd_{1,1}^0(\mathbb{P}^{n-1},\beta)$ under the $\mathbb{T}$-action are described in \cite[\S 1.3, 1.4]{Zin09b}. We exactly follow the terminology of loc. cit. So a fixed locus corresponds to a \emph{decorated graph} of type $A_i$, or type $A_{ij}$, or type $B_i$, or type $B_{ij}$, where $1\leq j\leq n$.

\begin{proposition}\label{prop-contribution-typeA}
\begin{enumerate}
      \item[(i)]   The contribution of graphs of type $A_i$ to  $\mathcal{F}(\alpha,\alpha_i,q)$ is equal to
      \begin{equation}\label{eq-contribution-Ai}
      \mathcal{A}_i(u) := \frac{1}{\Phi_0(\alpha_i,u)}\mathrm{Res}_{\hbar_1=0}\mathrm{Res}_{\hbar_2=0}\big\{
      e^{-\eta_i(u)/\hbar_1}e^{-\eta_i(u)/\hbar_2}
      \widetilde{\mathcal{Z}}^*_{ii}(\hbar_1,\hbar_2,u)\big\}.
      \end{equation}
      \item[(ii)] The contribution of graphs of type $\tilde{A}_{ij}$ to $\mathcal{F}(\alpha,\alpha_i,q)$ is equal to
\begin{equation}\label{eq-contribution-Aij}
      \tilde{\mathcal{A}}_{ij}(u):=
            \frac{\mathcal{A}_j(u)}{\prod_{k\neq j}(\alpha_j- \alpha_k)}
      \mathrm{Res}_{\hbar=0}\big\{ e^{-\eta_j(u)/\hbar} \mathcal{Z}_{ji}^*(\hbar,u)\big\}.
      \end{equation}
\end{enumerate}
\end{proposition}

\begin{proposition}\label{prop-contribution-typeB}
\begin{enumerate}
      \item   The contribution of graphs of type $B_i$ to  $\mathcal{F}(\alpha,\alpha_i,q)$ is equal to   
      \begin{eqnarray*}
      \mathcal{B}_i(u)=\frac{\prod_{k=1}^r (d_k\alpha_i)}{24}\mathrm{Res}_{\hbar=0,\infty,- \mathbf{d } \alpha_i}\big\{
      \frac{\prod_{k=1}^n(\alpha_i- \alpha_k+\hbar)}{\hbar^3\prod_{k=1}^r (d_k \alpha_i+\hbar)}
      \frac{\mathcal{Z}_i^*(\hbar,u)}{1+\mathcal{Z}_i^*(\hbar,u)}\big\}.
      \end{eqnarray*}  
      \item The contribution of graphs of type $\tilde{B}_{ij}$ to  $\mathcal{F}(\alpha,\alpha_i,q)$ is equal to  
      \begin{eqnarray*}
      \tilde{\mathcal{B}}_{ij}(u)=-\frac{\prod_{k=1}^r (d_k\alpha_j)}{24\prod_{k\neq j}(\alpha_j- \alpha_k)}
      \mathrm{Res}_{\hbar=0,\infty,- \mathbf{d} \alpha_j}\big\{
      \frac{\prod_{k=1}^n(\alpha_j- \alpha_k+\hbar)}{\hbar^2\prod_{k=1}^r(d_k\alpha_j+\hbar)}
      \frac{\mathcal{Z}_{ji}^*(\hbar,u)}{1+\mathcal{Z}_j^*(\hbar,u)}\big\}.
      \end{eqnarray*}
\end{enumerate}
\end{proposition}

Combining Proposition \ref{prop-contribution-typeA} and \ref{prop-contribution-typeB}, we get
\begin{eqnarray*}
&&\mathcal{F}(\alpha,\alpha_i,q)
      =\sum_{\beta=1}^{\infty}q^\beta 
      \int_{\Mtd_{1,1}^{0}(\mathbb{P}^{n-1},\beta)}\mathbf{e}(\mathcal{V}_1)\mathrm{ev}_1^*\phi_i\\
&=& \mathcal{A}_i(q)+\sum_{j=1}^n \tilde{\mathcal{A}}_{ij}(q)+\mathcal{B}_i(q)+\sum_{j=1}^n \tilde{\mathcal{B}}_{ij}(q)\\
&=& \frac{1}{\Phi_0(\alpha_i,q)}\mathrm{Res}_{\hbar_1=0}\mathrm{Res}_{\hbar_2=0}\big(
e^{-\eta_i(q)/\hbar_1}e^{-\eta_i(q)/\hbar_2}
\widetilde{\mathcal{Z}}^*_{ii}(\hbar_1,\hbar_2,q)\big)\\
&&+\sum_{j=1}^n\frac{\mathcal{A}_j(q)}{\prod_{k\neq j}(\alpha_j- \alpha_k)}
 \mathrm{Res}_{\hbar=0}\big( e^{-\eta_j(q)/\hbar} \mathcal{Z}_{ji}^*(\hbar,q)\big)\\
&&+\frac{\prod_{k=1}^r (d_k\alpha_i)}{24}\mathrm{Res}_{\hbar=0,\infty,- \mathbf{d } \alpha_i}\big(
\frac{\prod_{k=1}^n(\alpha_i- \alpha_k+\hbar)}{\hbar^3\prod_{k=1}^r (d_k \alpha_i+\hbar)}
\frac{\mathcal{Z}_i^*(\hbar,q)}{1+\mathcal{Z}_i^*(\hbar,q)}\big)\\
&&-\sum_{j=1}^n\frac{\prod_{k=1}^r (d_k\alpha_j)}{24\prod_{k\neq j}(\alpha_j- \alpha_k)}
\mathrm{Res}_{\hbar=0,\infty,- \mathbf{d} \alpha_j}\big(
\frac{\prod_{k=1}^n(\alpha_j- \alpha_k+\hbar)}{\hbar^2\prod_{k=1}^r(d_k\alpha_j+\hbar)}
\frac{\mathcal{Z}_{ji}^*(\hbar,q)}{1+\mathcal{Z}_j^*(\hbar,q)}\big).
\end{eqnarray*}

\subsection{Hypergeometric series associated with Fano complete intersections}
In this section, we recall some hypergeometric-type series introduced in \cite{Zin09b}, \cite{Pop12}, and \cite{PoZ14} that are related to genus 0 Gromov-Witten invariants of Fano complete intersections in $\mathbb{P}^n$. Let
\begin{equation*}
      \mathcal{P}:=\{1+\sum_{i=1}^{\infty}a_i(w)q^i\in \mathbb{Q}(w)[[q]]\big| a_i(w)\ \mbox{are regular at $w=0$ for $i\geq 1$}
      \},
\end{equation*}
\begin{equation}\label{eq-def-mathbfD}
      \mathbf{D}: \mathbb{Q}(w)[[q]]\rightarrow \mathbb{Q}(w)[[q]],\
      \mathbf{D}H(w,q):=(1+\frac{q}{w}\frac{d}{dq})H(w,q),
\end{equation}
\begin{equation*}
      \tilde{\hypergeometricF}(w,q):=\sum_{\beta=0}^{\infty}q^{\beta}w^{\nu_{\mathbf{d}}\beta}\frac{\prod_{k=1}^r\prod_{i=1}^{d_k \beta}(d_k w+i)}{\prod_{j=1}^{\beta}\big((w+j)^n-w^n\big)},
\end{equation*}
\begin{equation*}
      \hypergeometricF(w,q):=\sum_{\beta=0}^{\infty}q^{\beta}w^{\nu_{\mathbf{d}}\beta}\frac{\prod_{k=1}^r\prod_{i=1}^{d_k \beta}(d_k w+i)}{\prod_{j=1}^{\beta}(w+j)^n}.
\end{equation*}
Define $\PZc_{p,l}^{(\beta)}$, $\tilde{\PZc}_{p,l}^{(\beta)}\in \mathbb{Q}$ with $p,\beta,l\geq 0$ by
\begin{equation}\label{eq-def-numbers-c}
      \sum_{\beta=0}^{\infty}\sum_{l=0}^\infty \PZc_{p,l}^{(\beta)} w^l q^{\beta}
      =w^p \mathbf{D}^p F(w,\frac{q}{w^{\nu_{\mathbf{d}}}})
      =\sum_{\beta=0}^{\infty}q^{\beta}\frac{(w+\beta)^p \prod_{k=1}^{r}\prod_{i=1}^{d_k \beta}(d_k w+i)}
      {\prod_{j=1}^{\beta}(w+j)^n},
\end{equation}
\begin{equation}\label{eq-def-numbers-tildec}
      \sum_{\begin{subarray}{c}\beta_1+\beta_2=\beta\\ \beta_1,\beta_2\geq 0\end{subarray}}
      \sum_{k=0}^{p- \nu_{\mathbf{d}}\beta_1}\tilde{\PZc}_{p,k}^{(\beta_1)}\PZc_{k,l}^{(\beta_2)}=\delta_{\beta,0}
      \delta_{p,l},\ \mbox{for}\ \beta,l\in \mathbb{Z}_{\geq 0},\ l\leq p- \nu_{\mathbf{d}}\beta.
\end{equation}
As a convention, we set
\begin{equation}
      \tilde{\PZc}_{p,l}^{(\beta)}=0,\ \mbox{for}\ l<0\ \mbox{or}\ p<0.
\end{equation}
For $p\geq 0$ we have (\cite[(2.7)]{PoZ14})
\begin{equation}\label{eq-tildePZCoefficient-KroneckerProperty}
      c^{(0)}_{p,k}=\tilde{\PZc}^{(0)}_{p,k}=\delta_{p,k}.
\end{equation}
One can show 
\begin{eqnarray}\label{eq-PZCoefficients-manifestRecursion}
      \tilde{\PZc}_{p,l}^{(\beta)}
      =\delta_{\beta,0}\delta_{p,l}
      -\sum_{\beta_1=0}^{\beta-1}
      \sum_{k=0}^{p- \nu_{\mathbf{d}}\beta_1}\tilde{\PZc}_{p,k}^{(\beta_1)}\PZc_{k,l}^{(\beta- \beta_1)},\
       \mbox{for}\ \beta,l\in \mathbb{Z}_{\geq 0},\ l\leq p- \nu_{\mathbf{d}}\beta.
\end{eqnarray}
For $p\geq 0$, define
\begin{eqnarray}\label{eq-def-tildeFp(w,q)}
\tilde{\hypergeometricF}_{p}(w,q)&:=&\sum_{\beta=0}^{\infty}\sum_{l=0}^{p- \nu_{\mathbf{d}}\beta}\frac{\tilde{\PZc}_{p,l}^{(\beta)}q^{\beta}}{w^{p- \nu_{\mathbf{d}}\beta-l}}\mathbf{D}^l \tilde{\hypergeometricF}(w,q),
\end{eqnarray}
\begin{equation}\label{eq-def-Fp(w,q)}
      \hypergeometricF_{p}(w,q):=\sum_{\beta=0}^{\infty}\sum_{l=0}^{p- \nu_{\mathbf{d}}\beta}\frac{\tilde{\PZc}_{p,l}^{(\beta)}q^{\beta}}{w^{p- \nu_{\mathbf{d}}\beta-l}}\mathbf{D}^l \hypergeometricF(w,q).
\end{equation}
In particular,$\tilde{\hypergeometricF}_0(w,q)=\tilde{\hypergeometricF}(w,q)$ and $\hypergeometricF_0(w,q)=\hypergeometricF(w,q)$.
Then we define
\begin{equation}\label{eq-def-mathbbF}
      \mathbb{F}(w_1,w_2,q):=\sum_{p=0}^{n-1-r}\tilde{\hypergeometricF}_p(w_1,q)\tilde{\hypergeometricF}_{n-1-r-p}(w_2,q)
      +\sum_{p=1}^r \tilde{\hypergeometricF}_{n-p}(w_1,q)\tilde{\hypergeometricF}_{n-1-r+p}(w_2,q).
\end{equation}

Following \cite[page 465]{PoZ14}, we define $\mathcal{C}_{p,s}^{(b)}$, $\tilde{\mathcal{C}}_{p,s}^{(b)}\in \mathbb{Q}[\alpha_1,\dots,\alpha_n][[q]]$ with $p,b,s\geq 0$ by
\begin{eqnarray}\label{eq-def-numbers-C}
&&\sum_{s=0}^{\infty}\sum_{b=0}^{s} \mathcal{C}_{p,s}^{(b)} z^{s+r-b}w^s
      =w^p z^r\sum_{\beta=0}^{\infty}q^{\beta}\frac{(z+\frac{\beta}{w})^p \prod_{k=1}^{r}\prod_{i=1}^{d_k \beta}(d_k z+\frac{i}{w})}
      {\prod_{j=1}^{\beta}\prod_{k=1}^n(z- \alpha_k+\frac{j}{w})}\nn\\
&=& z^r\sum_{\beta=0}^{\infty}q^{\beta}w^{\nu_{\mathbf{d}}\beta}\frac{(wz+\beta)^p \prod_{k=1}^{r}\prod_{i=1}^{d_k \beta}(d_k w z+i)}
      {\prod_{j=1}^{\beta}\prod_{k=1}^n\big(w(z- \alpha_k)+j\big)},
\end{eqnarray}
and
\begin{equation}\label{eq-def-numbers-tildeC}
      \sum_{\begin{subarray}{c}b_1+b_2=b\\ 
      b_1, b_2\geq 0\end{subarray}}
      \sum_{t=0}^{p- b_1}\tilde{\mathcal{C}}_{p,t}^{(b_1)}\mathcal{C}_{t,s-b_1}^{(b_2)}=\delta_{b,0}
      \delta_{p,s},\ \mbox{for}\ b,s\in \mathbb{Z}_{\geq 0},\ b\leq s\leq p.
\end{equation}

\begin{definition}\label{def-SpecializedAtRootsOfUnity}
For a function $g$ of $\alpha_1,\dots,\alpha_n$, we denote  by $g(\zeta)$ the value of $g$  at $\alpha_k=\zeta_k=\exp(\frac{2k\pi \sqrt{-1}}{n})$ for $1\leq k\leq n$.
\end{definition}

\begin{lemma}\label{lem-tildeC-SpecializedAtRootsOfUnity}
Suppose $0\leq s\leq p$. Then
 $\tilde{\mathcal{C}}_{p,s}^{(b)}(\zeta)=0$ when $\nu_{\mathbf{d}}\nmid b$, and 
$\tilde{\mathcal{C}}_{p,s}^{(\nu_{\mathbf{d}}\gamma)}(\zeta)= \tilde{\PZc}_{p,s}^{(\gamma)}q^{\gamma}$ for $\gamma\in \mathbb{Z}_{\geq 0}$.
\end{lemma}
\begin{proof}
Manifestly $\tilde{\mathcal{C}}_{p,s}^{(b)}$ are recursively given by 
\begin{equation*}
      \mathcal{C}^{(0)}_{p,p}=0,\ \mathcal{C}^{(0)}_{p,s}=0\ \mbox{when}\ p>s,
\end{equation*}
and 
\begin{equation}\label{eq-recursion-PZseries}
\tilde{\mathcal{C}}_{p,s}^{(b)}
=\delta_{b,0}
      \delta_{p,s+b}-\sum_{t=0}^{s-1}\tilde{\mathcal{C}}_{p,t}^{(b)}\mathcal{C}_{t,s}^{(0)}
-\sum_{\begin{subarray}{c}b_1+b_2=b\\ 
      0\leq b_1\leq b-1\end{subarray}}
      \sum_{t=0}^{p- b_1}\tilde{\mathcal{C}}_{p,t}^{(b_1)}\mathcal{C}_{t,s+b-b_1}^{(b_2)},\ 
      \mbox{for}\ b,s\in \mathbb{Z}_{\geq 0},\  s\leq p.
\end{equation}
On the other hand we have
\begin{eqnarray*}
&&\sum_{s=0}^{\infty}\sum_{b=0}^{s} \mathcal{C}_{p,s}^{(b)}(\zeta) y^{s-b}w^b
=\sum_{\beta=0}^{\infty}q^{\beta}w^{\nu_{\mathbf{d}}\beta}\frac{(y+\beta)^p \prod_{k=1}^{r}\prod_{i=1}^{d_k \beta}(d_k y+i)}
      {\prod_{j=1}^{\beta}\big((y+j)^n-w^n\big)}\\
&\equiv & \sum_{\beta=0}^{\infty}q^{\beta}w^{\nu_{\mathbf{d}}\beta}\frac{(y+\beta)^p \prod_{k=1}^{r}\prod_{i=1}^{d_k \beta}(d_k y+i)}
      {\prod_{j=1}^{\beta}(y+j)^n}\ \mod w^n.
\end{eqnarray*}
Suppose $b<n$. Then $\mathcal{C}_{p,s}^{(b)}(\zeta)=0$ when $\nu_{\mathbf{d}}\nmid b$. If $b=\nu_{\mathbf{d}}\gamma$ then
\begin{equation*}
      \mathcal{C}_{p,s}^{(\nu_{\mathbf{d}}\gamma)}(\zeta)= \PZc_{p,s-\nu_{\mathbf{d}}\gamma}^{(\gamma)}q^{\gamma}.
\end{equation*}
Comparing (\ref{eq-PZCoefficients-manifestRecursion}) and (\ref{eq-recursion-PZseries}), we obtain the conclusion by induction on $b$.
\end{proof}

\subsection{Equivariant hypergeometric series}
In this section, we recall the formulae in \cite[\S 3]{PoZ14} expressing the equivariant generating functions in terms of equivariant hypergeometric series, and then compute their values at roots of unity. Define
\begin{equation}
      \mathcal{Y}(x,\hbar,q):=\sum_{\beta=0}^{\infty}q^{\beta}
      \frac{\prod_{k=1}^r\prod_{i=1}^{d_k \beta}(d_k x+i \hbar)}{\prod_{i=1}^{\beta}\big(\prod_{j=1}^n(x- \alpha_j+i\hbar)-\prod_{j=1}^n(x- \alpha_j)\big)},
\end{equation}
and for $p\geq 0$ define inductively
\begin{equation}
      \mathfrak{D}^0 \mathcal{Y}_0(x,\hbar,q):=x^r\mathcal{Y}(x,\hbar,q),
\end{equation}
\begin{equation}
      \mathfrak{D}^p \mathcal{Y}_0(x,\hbar,q):=(x+\hbar q\frac{\mathrm{d}}{\mathrm{d}q})\mathfrak{D}^{p-1}\mathcal{Y}_0(x,\hbar,q).
\end{equation}
Then we define
\begin{equation}
      \mathcal{Y}_p(x,\hbar,q):= \mathfrak{D}^p \mathcal{Y}_0(x,\hbar,q)
      +\sum_{i=1}^{p}\sum_{j=0}^{p-i}\widetilde{\mathcal{C}}^{(i)}_{p,j}\hbar^{p-i-j}\mathfrak{D}^j \mathcal{Y}_0(x,\hbar,q).
\end{equation}

\begin{proposition}\label{prop-equivariantGeneratingSeries-formulae}
      \begin{eqnarray}\label{eq-Zi-formula}
      1+\mathcal{Z}_i^*(\hbar,q)
      = \begin{cases}
      \mathcal{Y}(\alpha_i,\hbar,q),& \mbox{if}\ \nu_{\mathbf{d}}>1;\\
      e^{-\frac{q}{\hbar} \mathbf{d}!}
      \mathcal{Y}(\alpha_i,\hbar,q),& \mbox{if}\ \nu_{\mathbf{d}}=1.
      \end{cases}
      \end{eqnarray}
      Denote by $\sigma_i(\alpha)$ the $i$-th elementary symmetric polynomial in $\alpha_1,\dots,\alpha_n$. Then
      \begin{multline}\label{eq-tildeZij-formula}
      \widetilde{\mathcal{Z}}_{ij}^*(\hbar_1,\hbar_2,q)
      =\frac{1}{2\hbar_1\hbar_2(\hbar_1+\hbar_2)}\\
      \times\sum_{\begin{subarray}{c}p_1+p_2+l=n-1\\ p_1,p_2,l\geq 0\end{subarray}}(-1)^l \sigma_l(\alpha)\big(- \alpha_i^{p_1}\alpha_j^{p_2}+\frac{1}{\alpha_i^r} \mathcal{Z}_{p_1}(\alpha_1,\hbar_1,q) \mathcal{Z}_{p_2-r}(\alpha_2,\hbar_2,q)\big),
      \end{multline}
      and
      \begin{eqnarray}\label{eq-Zji-formula}
      \mathcal{Z}_{ji}^*(\hbar,q)
      =\frac{1}{\hbar} \begin{cases}
      \sum_{\begin{subarray}{c}p_1+p_2+v=n-1\\
      p_1,p_2,v\geq 0\end{subarray}}(-1)^v \sigma_v(\alpha)\big(- \alpha_j^{p_1}\alpha_i^{p_2}+\frac{\alpha_i^{p_2}}{\alpha_j^r} \mathcal{Y}_{p_1}(\alpha_j,\hbar,q)
      \big),\ \mbox{if}\ \nu_{\mathbf{d}}>1;\\
      \sum_{\begin{subarray}{c}p_1+p_2+v=n-1\\
      p_1,p_2,v\geq 0\end{subarray}}(-1)^v \sigma_v(\alpha)\big(- \alpha_j^{p_1}\alpha_i^{p_2}+\frac{\alpha_i^{p_2}}{\alpha_j^r} e^{-\frac{\mathbf{d}!}{\hbar}q}\mathcal{Y}_{p_1}(\alpha_j,\hbar,q)
      \big),\ \mbox{if}\ \nu_{\mathbf{d}}=1.
      \end{cases}
      \end{eqnarray}
\end{proposition}
\begin{proof}
From the definition we have
\begin{eqnarray*}
&&\mathcal{Z}_p^*(x,\hbar,q)|_{P_i}=
\mathcal{Z}_p^*(\alpha_i,\hbar,q)\\
&=&\alpha_i^{r+p}+\sum_{\beta=1}^{\infty}q^{\beta}\int_{\Mbar_{0,2}(\mathbb{P}^{n-1},\beta)}
      \frac{\mathbf{e}(\mathcal{V}''_0)\mathrm{ev}_2^* x^{r+p}}{\hbar-\psi_1}\mathrm{ev}_1^* \phi_i\\
&=&\alpha_i^{r+p}+\sum_{\beta=1}^{\infty}q^{\beta}\int_{\Mbar_{0,2}(\mathbb{P}^{n-1},\beta)}
      \frac{\mathbf{e}(\mathcal{V}'_0)\mathrm{ev}_2^* x^{p}}{\hbar-\psi_1}\mathrm{ev}_1^* (x^r\phi_i)\\
&=&\alpha_i^{r+p}+\alpha_i^r\sum_{\beta=1}^{\infty}q^{\beta}\int_{\Mbar_{0,2}(\mathbb{P}^{n-1},\beta)}
      \frac{\mathbf{e}(\mathcal{V}'_0)\mathrm{ev}_2^* x^{p}}{\hbar-\psi_1}\mathrm{ev}_1^* \phi_i.         
\end{eqnarray*}
In particular,
\begin{eqnarray*}
&&\mathcal{Z}_0^*(x,\hbar,q)|_{P_i}=
\mathcal{Z}_0^*(\alpha_i,\hbar,q)\\
&=&\alpha_i^{r}+\alpha_i^r\sum_{\beta=1}^{\infty}q^{\beta}\int_{\Mbar_{0,2}(\mathbb{P}^{n-1},\beta)}
      \frac{\mathbf{e}(\mathcal{V}'_0)}{\hbar-\psi_1}\mathrm{ev}_1^* \phi_i\\  
&=&\alpha_i^{r}+\alpha_i^r\sum_{\beta=1}^{\infty}q^{\beta}\int_{\Mbar_{0,1}(\mathbb{P}^{n-1},\beta)}
      \frac{\mathbf{e}(\mathcal{V}'_0)}{\hbar(\hbar-\psi_1)}\mathrm{ev}_1^* \phi_i\\ 
&=& \alpha_i^r\big(1+\mathcal{Z}_i^*(\hbar,q)\big).      
\end{eqnarray*}
So
\begin{equation}
    1+\mathcal{Z}_i^*(\hbar,q)=\frac{\mathcal{Z}_0^*(\alpha_i,\hbar,q)}{\alpha_i^r}.
\end{equation}
Then (\ref{eq-Zi-formula}) follows from \cite[Theorem 6]{PoZ14}.
Let 
\begin{equation}
      \mathcal{Z}^*(x_1,x_2,\hbar_1,\hbar_2,q):=\sum_{\beta=1}^{\infty}q^{\beta} (\mathrm{ev}_{1}\times \mathrm{ev}_2)_*\Big(
      \frac{\mathbf{e}(\mathcal{V}_0)}{(\hbar_1-\psi_1)(\hbar_2-\psi_2)}\Big).
\end{equation}
Thus
\begin{eqnarray*}
&&\mathcal{Z}^*(\alpha_i, \alpha_j,\hbar_1,\hbar_2,q)\\
&=&\sum_{\beta=1}^{\infty}q^{\beta}\int_{\Mbar_{0,2}(\mathbb{P}^{n-1},\beta)} 
      \frac{\mathbf{e}(\mathcal{V}_0)\mathrm{ev}_1^* \phi_i \mathrm{ev}_2^* \phi_j}
      {(\hbar_1-\psi_1)(\hbar_2-\psi_2)}\\
&=& \prod_{k=1}^r (d_k \alpha_i) \sum_{\beta=1}^{\infty}q^{\beta}\int_{\Mbar_{0,2}(\mathbb{P}^{n-1},\beta)} 
      \frac{\mathbf{e}(\mathcal{V}'_0)\mathrm{ev}_1^* \phi_i \mathrm{ev}_2^* \phi_j}
      {(\hbar_1-\psi_1)(\hbar_2-\psi_2)}\\
&=& 2 \hbar_1\hbar_2 \alpha_i^r (\prod_{k=1}^r d_k)\widetilde{\mathcal{Z}}_{ij}^*(\hbar_1,\hbar_2,q).
\end{eqnarray*}
So
\begin{equation}
   \widetilde{\mathcal{Z}}_{ij}^*(\hbar_1,\hbar_2,q)=\frac{\mathcal{Z}^*(\alpha_i, \alpha_j,\hbar_1,\hbar_2,q)}{2 \hbar_1\hbar_2 \alpha_i^r (\prod_{k=1}^r d_k)}.  
\end{equation}
Then (\ref{eq-tildeZij-formula}) follows from \cite[Theorem 6]{PoZ14}.

\begin{eqnarray*}
      && \mathcal{Z}_{ij}^*(\hbar,q)=2 \mathrm{Res}_{\hbar_2=0}\big\{\hbar_2 \widetilde{\mathcal{Z}}_{ij}^*(\hbar,\hbar_2,q)\big\}\\
      &=& \mathrm{Res}_{\hbar_2=0}\Big\{\frac{1}{\hbar(\hbar+\hbar_2)}\sum_{\begin{subarray}{c}p_1+p_2+l=n-1\\ p_1,p_2,l\geq 0\end{subarray}}(-1)^l \sigma_l\big(- \alpha_i^{p_1}\alpha_j^{p_2}+\frac{1}{\alpha_i^r} \mathcal{Z}_{p_1}(\alpha_1,\hbar,q) \mathcal{Z}_{p_2-r}(\alpha_2,\hbar_2,q)\big)\Big\}.
      \end{eqnarray*}

\begin{eqnarray*}
      &&\hbar_1\mathcal{Z}_{ji}^*(\hbar_1,u)=\frac{1}{\prod_{k=1}^r (d_k \alpha_j)}\mathrm{Res}_{\hbar_2=0}\mathcal{Z}^*(\alpha_j, \alpha_i,\hbar_1,\hbar_2,u)\\
&=& \mathrm{Res}_{\hbar_2=0}\Big\{ \frac{1}{\hbar_1+\hbar_2}\sum_{\begin{subarray}{c}p_1+p_2+v=n-1\\
p_1,p_2,v\geq 0\end{subarray}}(-1)^v \sigma_v\big(- \alpha_j^{p_1} \alpha_i^{p_2}+\frac{1}{\alpha_j^r} \mathcal{Z}_{p_1}(\alpha_j,\hbar_1,q) \mathcal{Z}_{p_2-r}(\alpha_i,\hbar_2,q)
\big)\Big\}.
\end{eqnarray*}
By definition,
\begin{equation}
      \mathcal{Z}_{ji}^*(\hbar,q)=-2 \mathrm{Res}_{\hbar_2=\infty}\big\{
      \hbar_2\widetilde{\mathcal{Z}}_{ji}^*(\hbar,\hbar_2,q)      \big\}.
\end{equation}
By definition, $\mathcal{Z}_p(\alpha_i,\hbar,q)$ is holomorphic at $\hbar=\infty$ and 
\begin{equation}
      \mathcal{Z}_p(\alpha_i,\hbar,q)|_{\hbar=\infty}=\alpha_i^{r+p}.
\end{equation}
So we get (\ref{eq-Zji-formula}).
\end{proof}

Now we compute the values of the equivariant generating series in Proposition \ref{prop-equivariantGeneratingSeries-formulae} at the roots of unity. 
\begin{lemma}\label{lem-calY-SpecializedAtRootsOfUnity}
\begin{equation}\label{eq-calY-SpecializedAtRootsOfUnity}
      \mathcal{Y}(x,\hbar,q)(\zeta)=\tilde{\hypergeometricF}(\frac{x}{\hbar},\frac{q}{x^{\nu_{\mathbf{d}}}}).
\end{equation}
\begin{equation}\label{eq-calYp-SpecializedAtRootsOfUnity}
\mathcal{Y}_p(x,\hbar,q)(\zeta)
      = x^{p+r}\tilde{\hypergeometricF}_{p}(\frac{x}{\hbar},\frac{q}{x^{\nu_{\mathbf{d}}}}).
\end{equation}
\end{lemma}
\begin{proof}
By definition we have
\begin{equation}
      \mathcal{Y}(x,\hbar,q)(\zeta)=\sum_{\beta=0}^{\infty}q^{\beta}
      \frac{\prod_{k=1}^r\prod_{i=1}^{d_k \beta}(d_k x+i \hbar)}{\prod_{i=1}^{\beta}\big((x+i\hbar)^n-x^n\big)}=\tilde{\hypergeometricF}(\frac{x}{\hbar},\frac{q}{x^{\nu_{\mathbf{d}}}}).
\end{equation}
By (\ref{eq-def-mathbfD}) and (\ref{eq-calY-SpecializedAtRootsOfUnity}) we have 
\begin{equation}\label{eq-DpCalY-SpecializedAtRootsOfUnity}
      \mathfrak{D}^p \mathcal{Y}_0(x,\hbar,q)(\zeta)=x^{p+r}(1+\frac{\hbar q}{x}\frac{\mathrm{d}}{\mathrm{d}q})^p\mathcal{Y}(x,\hbar,q)(\zeta)
      =x^{p+r}(\mathbf{D}^p\tilde{\hypergeometricF})(\frac{x}{\hbar},\frac{q}{x^{\nu_{\mathbf{d}}}}).
\end{equation}
Taking $p=0$ we get the $p=0$ case of (\ref{eq-calYp-SpecializedAtRootsOfUnity}). 
By Lemma \ref{lem-tildeC-SpecializedAtRootsOfUnity},
\begin{eqnarray}\label{eq-calYp-SpecializedAtRootsOfUnity-recursion}
      \mathcal{Y}_p(x,\hbar,q)(\zeta)
      =\mathfrak{D}^p \mathcal{Y}_0(x,\hbar,q)(\zeta)
      +\sum_{\gamma=1}^{\infty}\sum_{j=0}^{p-\nu_{\mathbf{d}}\gamma}\tilde{\PZc}_{p,j}^{(\gamma)}q^{\gamma}\hbar^{p-\nu_{\mathbf{d}}\gamma-j}\mathfrak{D}^j \mathcal{Y}_0(x,\hbar,q)(\zeta).
\end{eqnarray}
By (\ref{eq-tildePZCoefficient-KroneckerProperty}) and (\ref{eq-def-tildeFp(w,q)}), 
\begin{eqnarray}
\tilde{\hypergeometricF}_{p}(w,q)= \mathbf{D}^p \tilde{\hypergeometricF}(w,q)
+\sum_{\beta=1}^{\infty}\sum_{l=0}^{p- \nu_{\mathbf{d}}\beta}\frac{\tilde{\PZc}_{p,l}^{(\beta)}q^{\beta}}{w^{p- \nu_{\mathbf{d}}\beta-l}}\mathbf{D}^l \tilde{\hypergeometricF}(w,q).
\end{eqnarray}
Comparing this with (\ref{eq-calYp-SpecializedAtRootsOfUnity-recursion}), by (\ref{eq-DpCalY-SpecializedAtRootsOfUnity}) and  induction on $p$ we get (\ref{eq-calYp-SpecializedAtRootsOfUnity}).
\end{proof}

\begin{corollary}\label{cor-equivariant-Z-SpecializedAtRootsOfUnity}
\begin{equation}\label{eq-equivariant-Zi-SpecializedAtRootsOfUnity}
      1+ \mathcal{Z}_i^*(\hbar,q)(\zeta)
      =\begin{cases}
       \tilde{\hypergeometricF}(\frac{\zeta_i}{\hbar},\frac{q}{\zeta_i^{\nu_{\mathbf{d}}}}),\ \mbox{if}\ \nu_{\mathbf{d}}>1;\\
      e^{-\frac{q}{\hbar}\mathbf{d}!}\tilde{\hypergeometricF}(\frac{\zeta_i}{\hbar},\frac{q}{\zeta_i^{\nu_{\mathbf{d}}}}),\ \mbox{if}\ \nu_{\mathbf{d}}=1.
      \end{cases}
\end{equation}
\begin{eqnarray}\label{eq-equivariant-Zji-SpecializedAtRootsOfUnity}
\hbar\mathcal{Z}_{ji}^*(\hbar,q)(\zeta)
=\begin{cases}
\sum_{p=0}^{n-1}\zeta_j^{p}\zeta_i^{n-1-p}\big(-1+\tilde{\hypergeometricF}_{p}(\frac{\zeta_j}{\hbar},\frac{q}{\zeta_j^{\nu_{\mathbf{d}}}})\big),\ \mbox{if}\ \nu_{\mathbf{d}}>1;\\
\sum_{p=0}^{n-1}\zeta_j^{p}\zeta_i^{n-1-p}\big(-1 +e^{-\frac{\mathbf{d}!}{\hbar}q}
\tilde{\hypergeometricF}_{p}(\frac{\zeta_j}{\hbar},\frac{q}{\zeta_j^{\nu_{\mathbf{d}}}})
\big),\ \mbox{if}\ \nu_{\mathbf{d}}=1.
\end{cases}
\end{eqnarray}
\begin{eqnarray}\label{eq-equivariant-Zii-SpecializedAtRootsOfUnity}
&&2\hbar_1\hbar_2(\hbar_1+\hbar_2) \widetilde{\mathcal{Z}}_{ii}^*(\hbar_1,\hbar_2,q \alpha_i^{\nu_{\mathbf{d}}})(\zeta)
+n \zeta_i^{n-1}\nn\\
&&=\begin{cases}
\zeta_i^{n-1} \mathbb{F}(\frac{\zeta_i}{\hbar_1},\frac{\zeta_i}{\hbar_2},q),\
\mbox{if}\ \nu_{\mathbf{d}}>1;\\
\zeta_i^{n-1}
\exp\big(-\mathbf{d}!q \zeta_i(\frac{1}{\hbar_1}+\frac{1}{\hbar_2})\big) 
\mathbb{F}(\frac{\zeta_i}{\hbar_1},\frac{\zeta_i}{\hbar_2},q),\
\mbox{if}\ \nu_{\mathbf{d}}=1.
\end{cases}
\end{eqnarray}
\end{corollary}
\begin{proof}
(\ref{eq-equivariant-Zi-SpecializedAtRootsOfUnity}) follows from (\ref{eq-Zi-formula}) and (\ref{eq-calY-SpecializedAtRootsOfUnity}). 
Note that $\sigma_i(\zeta)=0$ for $1\leq i\leq n-1$, and $\sigma_0(\zeta)=1$. Then (\ref{eq-equivariant-Zji-SpecializedAtRootsOfUnity}) follows from (\ref{eq-tildeZij-formula}) and  (\ref{eq-calYp-SpecializedAtRootsOfUnity}). Finally (\ref{eq-equivariant-Zii-SpecializedAtRootsOfUnity}) follows from (\ref{eq-tildeZij-formula}), (\ref{eq-equivariant-Zi-SpecializedAtRootsOfUnity}), and (\ref{eq-def-mathbbF}).
\end{proof}

\begin{corollary}\label{cor-regularizability-nonEquivariant}
$\tilde{\hypergeometricF}(\frac{1}{\hbar},q)-1$ is regularizable.
\end{corollary}
\begin{proof}
By Proposition \ref{prop-regularizability-Equivariant}, the equivalent condition for regularizability \cite[(2.3)]{Zin09b} holds for $\mathcal{Z}_i^*(\hbar,q)$. Then by  (\ref{eq-equivariant-Zi-SpecializedAtRootsOfUnity}),
\cite[(2.3)]{Zin09b} also holds for $\tilde{\hypergeometricF}(\frac{\zeta_i}{\hbar},q)-1$. So $\tilde{\hypergeometricF}(\frac{1}{\hbar},q)-1$ is regularizable.
\end{proof}

Define
\begin{equation}
      \mu(q):=\mathrm{Res}_{\hbar=0}\big(\ln \tilde{\hypergeometricF}(\frac{1}{\hbar},q)\big),
\end{equation}
\begin{equation*}
      L(q)=1+q\frac{\mathrm{d}}{\mathrm{d}q} \mu(q),
\end{equation*}
and
\begin{equation}
      Q(\hbar,q):=e^{-\frac{\mu(q)}{\hbar}}  \tilde{\hypergeometricF}(\frac{1}{\hbar},q),
\end{equation}
By Corollary \ref{cor-regularizability-nonEquivariant}, $Q(\hbar,q)$ is holomorphic at $\hbar=0$. This enables us to define
\begin{equation}
      \Phi_0(q)=Q(0,q),\ \Phi_1(q):=\frac{d}{d\hbar}Q(\hbar,q)|_{\hbar=0}.
\end{equation}
The closed formulae of $\mu(q)$, $L(q)$, $\Phi_0(q)$ and $\Phi_1(q)$ will be given in Section \ref{sec:formulae-hypergeometricSeries-Fano}.

\subsection{Localization Contributions  valued at roots of unity}
In this section, we evaluate the localization contributions obtained in Proposition \ref{prop-contribution-typeA} and \ref{prop-contribution-typeB}, at the roots of unity. The results are Proposition \ref{prop-typeA-contribution-SpecializedAtRootsOfUnity-final} and 
\ref{prop-typeB-contribution-SpecializedAtRootsOfUnity-final}. Then Lemma \ref{lem-reducedGenus1GW-generatingFunc} and Remark \ref{rem:comparision-Zinger'sMethod-ours} one gets the reduced genus 1 Gromov-Witten invariants $\langle \sfh_{1+\nu_{\mathbf{d}}\beta}\rangle_{1,\beta}^0$.

\subsubsection{Contributions of type A valued at roots of unity}

\begin{proposition}\label{prop-typeA-contribution-SpecializedAtRootsOfUnity}
\begin{eqnarray}\label{eq-typeA-contribution-SpecializedAtRootsOfUnity}
&&\mathcal{A}_i(q \alpha_j^{\nu_{\mathbf{d}}})(\zeta)
+\sum_{j=1}^n\tilde{\mathcal{A}}_{ji}(q \alpha_j^{\nu_{\mathbf{d}}})(\zeta)\nn\\
&=&  \frac{\zeta_{i}^{n-2}}{2}\sum_{p=0}^{n-1}q^{\frac{p-1}{\nu_{\mathbf{d}}}}
\mathrm{Coeff}_{q^{\frac{p-1}{\nu_{\mathbf{d}}}}}\Big\{\mathrm{Res}_{\hbar_1=0}\mathrm{Res}_{\hbar_2=0}
\big(\frac{1}{\hbar_1\hbar_2(\hbar_1+\hbar_2)}
e^{-\mu(q)(\frac{1}{\hbar_1}+\frac{1}{\hbar_2})}
 \mathbb{F}(\frac{1}{\hbar_1},\frac{1}{\hbar_2},q)\big)\nn\\
 &&\times\sum_{\beta=0}^{\infty}\tilde{\PZc}_{p,p- \nu_{\mathbf{d}}\beta}^{(\beta)}q^{\beta}
L(q)^{p- \nu_{\mathbf{d}}\beta}\Big\}.
\end{eqnarray}
\end{proposition}
\begin{proof}
By (\ref{eq-equivariant-Zi-SpecializedAtRootsOfUnity}), we have
\begin{eqnarray*}
      \eta_i(q \alpha_i^{\nu_{\mathbf{d}}})(\zeta)=\zeta_i\cdot
      \begin{cases}
      \mu(q),& \mbox{if}\ \nu_{\mathbf{d}}>1,\\
      - \mathbf{d}! q+\mu(q),& \mbox{if}\ \nu_{\mathbf{d}}=1.
      \end{cases}
\end{eqnarray*}
and
\begin{equation*}
      \Phi_0(\alpha_i,q \alpha_i^{\nu_{\mathbf{d}}})(\zeta)=\Phi_0(q).
\end{equation*}
Then by (\ref{eq-contribution-Ai}) and (\ref{eq-equivariant-Zii-SpecializedAtRootsOfUnity}),
\begin{eqnarray*}
\mathcal{A}_i(q \alpha_i^{\nu_{\mathbf{d}}})(\zeta)
&=& \frac{\zeta_i^{n-2}}{2\Phi_0(q)}\mathrm{Res}_{\hbar_1=0}\mathrm{Res}_{\hbar_2=0}\big\{\frac{1}{\hbar_1\hbar_2(\hbar_1+\hbar_2)}
e^{-\mu(q)(\frac{1}{\hbar_1}+\frac{1}{\hbar_2})} \mathbb{F}(\frac{1}{\hbar_1},\frac{1}{\hbar_2},q)\big\}.
\end{eqnarray*}
By (\ref{eq-contribution-Aij}) and (\ref{eq-equivariant-Zji-SpecializedAtRootsOfUnity}), when $\nu_{\mathbf{d}}>1$,
\begin{eqnarray*}
&&\sum_{j=1}^{n}\tilde{\mathcal{A}}_{ij}(q \alpha_i^{\nu_{\mathbf{d}}})(\zeta)\\
&=&\sum_{j=1}^{n}\sum_{p=0}^{n-1}\frac{\alpha_i^{n-1-p}\alpha_j^{p+n-2}\mathrm{Res}_{\hbar_1=0}\mathrm{Res}_{\hbar_2=0}\big\{\frac{1}{\hbar_1\hbar_2(\hbar_1+\hbar_2)}
e^{-\mu(\frac{q \alpha_i^{\nu_{\mathbf{d}}}}{\alpha_j^{\nu_{\mathbf{d}}}})(\frac{1}{\hbar_1}+\frac{1}{\hbar_2})}
 \mathbb{F}(\frac{1}{\hbar_1},\frac{1}{\hbar_2},\frac{q \alpha_i^{\nu_{\mathbf{d}}}}{\alpha_j^{\nu_{\mathbf{d}}}})\big\}}{2\Phi_0(\frac{q \alpha_i^{\nu_{\mathbf{d}}}}{\alpha_j^{\nu_{\mathbf{d}}}})\prod_{k\neq j}(\alpha_j- \alpha_k)}\\
&&\times 
\mathrm{Res}_{\hbar=0}\Big\{\frac{1}{\hbar}e^{-\mu(\frac{q \alpha_i^{\nu_{\mathbf{d}}}}{\alpha_j^{\nu_\mathbf{{d}}}})\frac{1}{\hbar}}
\big(-1+\tilde{\hypergeometricF}_{p}(\frac{1}{\hbar},\frac{q \alpha_i^{\nu_{\mathbf{d}}}}{\alpha_j^{\nu_{\mathbf{d}}}})\big)\Big\}\\
&=&  \alpha_{i}^{n-2}\sum_{p=0}^{n-1}q^{\frac{p-1}{\nu_{\mathbf{d}}}}
\mathrm{Coeff}_{q^{\frac{p-1}{\nu_{\mathbf{d}}}}}\bigg\{\frac{\mathrm{Res}_{\hbar_1=0}\mathrm{Res}_{\hbar_2=0}\big\{\frac{1}{\hbar_1\hbar_2(\hbar_1+\hbar_2)}
e^{-\mu(q)(\frac{1}{\hbar_1}+\frac{1}{\hbar_2})}
 \mathbb{F}(\frac{1}{\hbar_1},\frac{1}{\hbar_2},q)\big\}}{2\Phi_0(q)}\\
&&\times 
\mathrm{Res}_{\hbar=0}\Big\{\frac{1}{\hbar}e^{-\mu(q)\frac{1}{\hbar}}
\big(-1+\tilde{\hypergeometricF}_{p}(\frac{1}{\hbar},q)\big)\Big\}\bigg\}.
\end{eqnarray*}
When $\nu_{\mathbf{d}}=1$, the computation needs to be slightly modified: we need to replace the last residue by
\begin{equation*}
      \mathrm{Res}_{\hbar=0}\Big\{\frac{1}{\hbar}e^{-\big(- \mathbf{d}! q+ \mu(q)\big)\frac{1}{\hbar}}
\big(-1+e^{-\frac{\mathbf{d}!q}{\hbar}}\tilde{\hypergeometricF}_{p}(\frac{1}{\hbar},q)\big)\Big\},
\end{equation*}
which is  still  equal to the unmodified residue. 
We have
\begin{eqnarray*}
&& \mathrm{Res}_{\hbar=0}\Big\{\frac{1}{\hbar} e^{-\mu(q)/\hbar} 
\tilde{\hypergeometricF}_p(\frac{1}{\hbar},q)\Big\}\\
&=& \mathrm{Res}_{\hbar=0}\Big\{\frac{1}{\hbar} e^{-\mu(q)/\hbar} 
\sum_{\beta=0}^{\infty}\sum_{l=0}^{p- \nu_{\mathbf{d}}\beta}\tilde{\PZc}_{p,l}^{(\beta)}q^{\beta}\hbar^{p- \nu_{\mathbf{d}}\beta-l}\mathbf{D}^l \tilde{\hypergeometricF}(\frac{1}{\hbar},q)\Big\}\\
&=& \mathrm{Res}_{\hbar=0}\Big\{\frac{1}{\hbar} 
\sum_{\beta=0}^{\infty}\sum_{l=0}^{p- \nu_{\mathbf{d}}\beta}\tilde{\PZc}_{p,l}^{(\beta)}q^{\beta}\hbar^{p- \nu_{\mathbf{d}}\beta-l}e^{-\mu(q)/\hbar} \mathbf{D}^l \tilde{\hypergeometricF}(\frac{1}{\hbar},q)\Big\}\\
&=& \mathrm{Res}_{\hbar=0}\Big\{\frac{1}{\hbar} 
\sum_{\beta=0}^{\infty}\sum_{l=0}^{p- \nu_{\mathbf{d}}\beta}\tilde{\PZc}_{p,l}^{(\beta)}q^{\beta}\hbar^{p- \nu_{\mathbf{d}}\beta-l}(1+q\frac{d \mu(q)}{dq}+\hbar q\frac{d}{d q})^l\big(e^{-\mu(q)/\hbar}\tilde{\hypergeometricF}(\frac{1}{\hbar},q)\big)\Big\}\\
&=& \sum_{\beta=0}^{\infty}\tilde{\PZc}_{p,p- \nu_{\mathbf{d}}\beta}^{(\beta)}q^{\beta}
(1+q\frac{d \mu(q)}{dq})^{p- \nu_{\mathbf{d}}\beta}\Phi_0(q),
\end{eqnarray*}
where in the last equality we use the regularizability of $\tilde{\hypergeometricF}(\frac{1}{\hbar},q)-1$. 
Hence we get (\ref{eq-typeA-contribution-SpecializedAtRootsOfUnity}).
\end{proof}

\begin{definition}
Define
\begin{eqnarray}\label{eq-def-Theta(0)}
\Theta_{p}^{(0)}(q):= \mathrm{Res}_{\hbar=0}\big\{\frac{1}{\hbar}e^{-\frac{\mu(q)}{\hbar}}\tilde{\hypergeometricF}_p(\frac{1}{\hbar},q)\big\},
\end{eqnarray}
and
\begin{eqnarray}\label{eq-def-Theta(1)}
\Theta_{p}^{(1)}(q)
:= \mathrm{Res}_{\hbar=0}\big\{\frac{1}{\hbar^2}e^{-\frac{\mu(q)}{\hbar}}\tilde{\hypergeometricF}_p(\frac{1}{\hbar},q)\big\}.
\end{eqnarray}
\end{definition}

\begin{lemma}\label{lem-Theta(0)-formula-0}
\begin{eqnarray}\label{eq-Theta(0)-formula-0}
\Theta_{p}^{(0)}(q)
= \Phi_0(q)\sum_{\beta=0}^{\infty}\tilde{\PZc}_{p,p- \nu_{\mathbf{d}}\beta}^{(\beta)}q^{\beta}
L(q)^{p- \nu_{\mathbf{d}}\beta}.
\end{eqnarray}
\end{lemma}
\begin{proof}
By computing the commutators relating to the operator $\mathbf{D}$, we have
\begin{eqnarray*}
&& e^{-\frac{\mu(q)}{\hbar}} \mathbf{D}^s \tilde{\hypergeometricF}(\frac{1}{\hbar},q)
= e^{-\frac{\mu(q)}{\hbar}}  (1+\hbar q\frac{d}{dq})^s \tilde{\hypergeometricF}(\frac{1}{\hbar},q)\\
&=& (1+q\frac{d \mu(q)}{dq}+\hbar q\frac{d}{d q})^s\big( e^{-\frac{\mu(q)}{\hbar}}  \tilde{\hypergeometricF}(\frac{1}{\hbar},q)\big).
\end{eqnarray*}
Then
\begin{eqnarray*}
\Theta_{p}^{(0)}(q)
&=&\mathrm{Res}_{\hbar=0}\Big\{\frac{1}{\hbar} 
\sum_{\beta=0}^{\infty}\sum_{s=0}^{p- \nu_{\mathbf{d}}\beta}
\tilde{\PZc}_{p,s}^{(\beta)}q^{\beta}\hbar^{p- \nu_{\mathbf{d}}\beta-s}(1+q\frac{d \mu(q)}{dq}+\hbar q\frac{d}{d q})^s\big( e^{-\frac{\mu(q)}{\hbar}}  \tilde{\hypergeometricF}(\frac{1}{\hbar},q)\big)\Big\}\nn\\
&=& \sum_{\beta=0}^{\infty}\tilde{\PZc}_{p,p- \nu_{\mathbf{d}}\beta}^{(\beta)}q^{\beta}
(1+q\frac{d \mu(q)}{dq})^{p- \nu_{\mathbf{d}}\beta}\Phi_0(q)
= \Phi_0(q)\sum_{\beta=0}^{\infty}\tilde{\PZc}_{p,p- \nu_{\mathbf{d}}\beta}^{(\beta)}q^{\beta}
L(q)^{p- \nu_{\mathbf{d}}\beta}.
\end{eqnarray*}
\end{proof}

\begin{lemma}\label{lem-Theta(1)-formula-0}
\begin{eqnarray}\label{eq-Theta(1)-formula-0}
\Theta_{p}^{(1)}(q)
&=& \Phi_0(q)\sum_{\beta=0}^{\infty}\tilde{\PZc}_{p,p- \nu_{\mathbf{d}}\beta-1}^{(\beta)}q^{\beta}
L(q)^{p- \nu_{\mathbf{d}}\beta-1}\nn\\
&&+ \Phi_1(q)\sum_{\beta=0}^{\infty}\tilde{\PZc}_{p,p- \nu_{\mathbf{d}}\beta}^{(\beta)}q^{\beta}
L(q)^{p- \nu_{\mathbf{d}}\beta}\nn\\
&&+\Phi'_0(q)\sum_{\beta=0}^{\infty}\tilde{\PZc}_{p,p- \nu_{\mathbf{d}}\beta}^{(\beta)}q^{\beta+1}
(p- \nu_{\mathbf{d}}\beta)
L(q)^{p- \nu_{\mathbf{d}}\beta-1}\nn\\
&&+L(q)'\Phi_0(q)\sum_{\beta=0}^{\infty}\tilde{\PZc}_{p,p- \nu_{\mathbf{d}}\beta}^{(\beta)}q^{\beta+1}
\binom{p- \nu_{\mathbf{d}}\beta}{2}
L(q)^{p- \nu_{\mathbf{d}}\beta-2}.
\end{eqnarray}
\end{lemma}
\begin{proof}
For any function $f(q)$ and $k\geq 0$, 
\begin{eqnarray*}
\big(f(q)+\hbar q\frac{d}{dq}\big)^k
= f(q)^k+\binom{k}{2}\hbar qf(q)^{k-2}\frac{df(q)}{dq}+k\hbar f(q)^{k-1}q\frac{d}{dq}+o(\hbar).
\end{eqnarray*}
Then
\begin{eqnarray*}
\Theta_{p}^{(1)}(q)
&=&\mathrm{Res}_{\hbar=0}\Big\{\frac{1}{\hbar^2} 
\sum_{\beta=0}^{\infty}\sum_{s=0}^{p- \nu_{\mathbf{d}}\beta}
\tilde{\PZc}_{p,s}^{(\beta)}q^{\beta}\hbar^{p- \nu_{\mathbf{d}}\beta-s}(1+q\frac{d \mu(q)}{dq}+\hbar q\frac{d}{d q})^s\big( e^{-\frac{\mu(q)}{\hbar}}  \tilde{\hypergeometricF}(\frac{1}{\hbar},q)\big)\Big\}\\
&=& \sum_{\beta=0}^{\infty}\tilde{\PZc}_{p,p- \nu_{\mathbf{d}}\beta-1}^{(\beta)}q^{\beta}
(1+q\frac{d \mu(q)}{dq})^{p- \nu_{\mathbf{d}}\beta-1}\Phi_0(q)\\
&&+ \sum_{\beta=0}^{\infty}\tilde{\PZc}_{p,p- \nu_{\mathbf{d}}\beta}^{(\beta)}q^{\beta}
(1+q\frac{d \mu(q)}{dq})^{p- \nu_{\mathbf{d}}\beta}\Phi_1(q)\\
&&+\sum_{\beta=0}^{\infty}\tilde{\PZc}_{p,p- \nu_{\mathbf{d}}\beta}^{(\beta)}q^{\beta+1}
(p- \nu_{\mathbf{d}}\beta)
(1+q\frac{d \mu(q)}{dq})^{p- \nu_{\mathbf{d}}\beta-1}\frac{d\Phi_0(q) }{dq}\\
&&+\sum_{\beta=0}^{\infty}\tilde{\PZc}_{p,p- \nu_{\mathbf{d}}\beta}^{(\beta)}q^{\beta+1}
\binom{p- \nu_{\mathbf{d}}\beta}{2}
(1+q\frac{d \mu(q)}{dq})^{p- \nu_{\mathbf{d}}\beta-2}\big(\frac{d \mu(q)}{dq}+q\frac{d^2 \mu(q)}{dq^2}\big)\Phi_0(q),
\end{eqnarray*}
So we get (\ref{eq-Theta(1)-formula-0}).

\end{proof}

\begin{proposition}\label{prop-typeA-contribution-SpecializedAtRootsOfUnity-final}
The type A Contribution to 
\begin{equation*}
      \sum_{\beta=0}^{\infty}q^{\beta}\langle \sfh_{1+(n-d)\beta}\rangle_{1,\beta}^0
=\sum_{\beta=0}^{\infty}q^{\beta} \mathcal{F}_{\beta,0}
\end{equation*}
is equal to
\begin{eqnarray*}
&& \frac{1}{2}\sum_{p=0}^{n-1}q^{\frac{p-1}{\nu_{\mathbf{d}}}}
\mathrm{Coeff}_{q^{\frac{p-1}{\nu_{\mathbf{d}}}}}\Big\{
\frac{\sum_{p=0}^{n-1-r}    \Theta_{p}^{(1)}(q)\Theta_{n-1-r-p}^{(0)}(q)
+\sum_{p=1}^{r}   \Theta_{n-p}^{(1)}(q)\Theta_{n-1-r+p}^{(0)}(q)}{\Phi_0(q)}\Theta_p^{(0)}(q)\Big\}.
\end{eqnarray*}

\end{proposition}
\begin{proof}
By Proposition \ref{prop-contribution-typeA} and \ref{prop-typeA-contribution-SpecializedAtRootsOfUnity}, and Lemma \ref{lem-Theta(0)-formula-0}, the type A contribution is equal to 
\begin{eqnarray*}
&& \frac{1}{2}\sum_{p=0}^{n-1}q^{\frac{p-1}{\nu_{\mathbf{d}}}}
\mathrm{Coeff}_{q^{\frac{p-1}{\nu_{\mathbf{d}}}}}\Big\{
\frac{A(q)}{\Phi_0(q)}\Theta_p^{(0)}(q)\Big\}
\end{eqnarray*}
where
\begin{eqnarray}\label{eq-A(q)-inTermsOfTheta}
A(q)&=&  \mathrm{Res}_{\hbar_1=0}\mathrm{Res}_{\hbar_2=0}\big\{\frac{1}{\hbar_1\hbar_2(\hbar_1+\hbar_2)}
e^{-\mu(q)(\frac{1}{\hbar_1}+\frac{1}{\hbar_2})}
 \mathbb{F}(\frac{1}{\hbar_1},\frac{1}{\hbar_2},q)\big\}\nn\\
&=&\mathrm{Res}_{\hbar_1=0}\mathrm{Res}_{\hbar_2=0}\Big\{\frac{1}{\hbar_1\hbar_2(\hbar_1+\hbar_2)}
e^{-\frac{\mu(q)}{\hbar_1}}e^{-\frac{\mu(q)}{\hbar_2}}\nn\\
&&\big(\sum_{p=0}^{n-1-r}\tilde{\hypergeometricF}_p(\frac{1}{\hbar_1},q)\tilde{\hypergeometricF}_{n-1-r-p}(\frac{1}{\hbar_2},q)
      +\sum_{p=1}^r \tilde{\hypergeometricF}_{n-p}(\frac{1}{\hbar_1},q)\tilde{\hypergeometricF}_{n-1-r+p}(\frac{1}{\hbar_2},q)\big)\Big\}\nn\\
&=& \sum_{p=0}^{n-1-r}  \Theta_{p}^{(1)}(q)\Theta_{n-1-r-p}^{(0)}(q)
+\sum_{p=1}^{r}   \Theta_{n-p}^{(1)}(q)\Theta_{n-1-r+p}^{(0)}(q).
\end{eqnarray}
\end{proof}

\subsubsection{Contributions of type B valued at roots of unity}
\begin{proposition}\label{prop-typeII-contribution-SpecializedAtRootsOfUnity}
\begin{multline}\label{eq-typeII-contribution-SpecializedAtRootsOfUnity}
\mathcal{B}_iq \alpha_i^{\nu_{\mathbf{d}}})(\zeta)
+\sum_{j=1}^n\tilde{\mathcal{B}}_{ij}(q \alpha_i^{\nu_{\mathbf{d}}})(\zeta)=\frac{\zeta_i^{n-2}\prod_{k=1}^r d_k}{24}\\
\times\sum_{p=0}^{n-1}q^{\frac{p-1}{\nu_{\mathbf{d}}}}
\mathrm{Coeff}_{q^{\frac{p-1}{\nu_{\mathbf{d}}}}}\bigg\{
\mathrm{Res}_{\hbar=0,\infty,- \mathbf{d} }\Big\{
\frac{(1+\hbar)^n-1}{\hbar^3\prod_{k=1}^r(d_k+\hbar)}
\frac{\tilde{\hypergeometricF}(\frac{1}{\hbar},q)-\tilde{\hypergeometricF}_p(\frac{1}{\hbar},q)}{\tilde{\hypergeometricF}(\frac{1}{\hbar},q)}\Big\}\bigg\}.
\end{multline}
\end{proposition}
\begin{proof}
First we suppose $\nu_{\mathbf{d}}>1$. 
By (\ref{eq-equivariant-Zi-SpecializedAtRootsOfUnity}),
\begin{eqnarray}\label{eq-Bi-SpecializedAtRootsOfUnity}
&&\mathcal{B}_i(q \alpha_i^{\nu_{\mathbf{d}}})(\zeta)=\frac{\prod_{k=1}^r (d_k\zeta_i)}{24}\mathrm{Res}_{\hbar=0,\infty,- \mathbf{d } \alpha_i}\big(
\frac{\prod_{k=1}^n(\zeta_i- \zeta_k+\hbar)}{\hbar^3\prod_{k=1}^r (d_k \zeta_i+\hbar)}
\frac{\mathcal{Z}_i^*(\hbar,q \zeta_i^{\nu_{\mathbf{d}}})}{1+\mathcal{Z}_i^*(\hbar,q \zeta_i^{\nu_{\mathbf{d}}})}\big)\nn\\
&=& \frac{\zeta_i^{n-2}\prod_{k=1}^r d_k}{24} \mathrm{Res}_{\hbar=0,\infty,- \mathbf{d}}\big(
\frac{(1+\hbar)^n-1}{\hbar^3\prod_{k=1}^r (d_k +\hbar)}
\frac{\tilde{\hypergeometricF}(\frac{1}{\hbar},q)-1}{\tilde{\hypergeometricF}(\frac{1}{\hbar},q)}\big).
\end{eqnarray}
By (\ref{eq-equivariant-Zi-SpecializedAtRootsOfUnity})  and (\ref{eq-equivariant-Zji-SpecializedAtRootsOfUnity}),
\begin{eqnarray*}
&&\tilde{\mathcal{B}}_{ij}(q \alpha_i^{\nu_{\mathbf{d}}})(\zeta)=-\frac{\prod_{k=1}^r (d_k\zeta_j)}{24\prod_{k\neq j}(\zeta_j- \zeta_k)}
\mathrm{Res}_{\hbar=0,\infty,- \mathbf{d} \zeta_j}\big(
\frac{\prod_{k=1}^n(\zeta_j- \zeta_k+\hbar)}{\hbar^2\prod_{k=1}^r(d_k\zeta_j+\hbar)}
\frac{\mathcal{Z}_{ji}^*(\hbar,q \zeta_i^{\nu_{\mathbf{d}}})}{1+\mathcal{Z}_j^*(\hbar,q \zeta_i^{\nu_{\mathbf{d}}})}\big)\\
&=& -\frac{\prod_{k=1}^r d_k}{24\prod_{k\neq j}(\zeta_j- \zeta_k)}
\sum_{p=0}^{n-1}\bigg(\zeta_j^{p+n-2}\zeta_i^{n-1-p}
\mathrm{Res}_{\hbar=0,\infty,- \mathbf{d} }\Big\{
\frac{(1+\hbar)^n-1}{\hbar^3\prod_{k=1}^r(d_k+\hbar)}
\frac{\tilde{\hypergeometricF}_p(\frac{1}{\hbar},\frac{q \zeta_i^{\nu_{\mathbf{d}}}}{\zeta_j^{\nu_{\mathbf{d}}}})-1}{\tilde{\hypergeometricF}(\frac{1}{\hbar},\frac{q \zeta_i^{\nu_{\mathbf{d}}}}{\zeta_j^{\nu_{\mathbf{d}}}})}\Big\}\bigg).
\end{eqnarray*}
So
\begin{multline}\label{eq-Bij-SpecializedAtRootsOfUnity}
\sum_{j=1}^n\tilde{\mathcal{B}}_{ij}(q \alpha_i^{\nu_{\mathbf{d}}})(\zeta)=-\frac{\zeta_i^{n-2}\prod_{k=1}^r d_k}{24}\\
\times\sum_{p=0}^{n-1}q^{\frac{p-1}{\nu_{\mathbf{d}}}}
\mathrm{Coeff}_{q^{\frac{p-1}{\nu_{\mathbf{d}}}}}\bigg\{
\mathrm{Res}_{\hbar=0,\infty,- \mathbf{d} }\Big\{
\frac{(1+\hbar)^n-1}{\hbar^3\prod_{k=1}^r(d_k+\hbar)}
\frac{\tilde{\hypergeometricF}_p(\frac{1}{\hbar},q)-1}{\tilde{\hypergeometricF}(\frac{1}{\hbar},q)}\Big\}\bigg\}.
\end{multline}
Then (\ref{eq-typeII-contribution-SpecializedAtRootsOfUnity}) follows. In the case $\nu_{\mathbf{d}}=1$, the fraction $\frac{\tilde{\hypergeometricF}(\frac{1}{\hbar},q)-1}{\tilde{\hypergeometricF}(\frac{1}{\hbar},q)}$ in the right hand side of (\ref{eq-Bi-SpecializedAtRootsOfUnity}) is replaced by 
$\frac{\tilde{\hypergeometricF}(\frac{1}{\hbar},q)-e^{\mathbf{d}!q/\hbar}}{\tilde{\hypergeometricF}(\frac{1}{\hbar},q)}$, and the fraction $\frac{\tilde{\hypergeometricF}_p(\frac{1}{\hbar},q)-1}{\tilde{\hypergeometricF}_p(\frac{1}{\hbar},q)}$ in the right hand side of (\ref{eq-Bij-SpecializedAtRootsOfUnity}) is replaced by 
$\frac{\tilde{\hypergeometricF}_p(\frac{1}{\hbar},q)-e^{\mathbf{d}!q/\hbar}}{\tilde{\hypergeometricF}_p(\frac{1}{\hbar},q)}$. So the final result remains the same.
\end{proof}

We need to compute the residues in (\ref{eq-typeII-contribution-SpecializedAtRootsOfUnity}).

\begin{lemma}\label{lem-typeII-contribution-SpecializedAtRootsOfUnity-residueAt-hbar=-d}
\begin{eqnarray}\label{eq-typeII-contribution-SpecializedAtRootsOfUnity-residueAt-hbar=-d}
&& \mathrm{Res}_{\hbar=-\mathbf{d} }\Big\{
\frac{(1+\hbar)^n-1}{\hbar^3\prod_{k=1}^r(d_k+\hbar)}
\frac{\tilde{\hypergeometricF}(\frac{1}{\hbar},q)-\tilde{\hypergeometricF}_p(\frac{1}{\hbar},q)}{\tilde{\hypergeometricF}(\frac{1}{\hbar},q)}\Big\}\nn\\
&=& - \mathrm{Res}_{\hbar=0,\infty }\Big\{
\frac{(1+\hbar)^n-1}{\hbar^3\prod_{k=1}^r(d_k+\hbar)}
(1-\sum_{\beta=0}^{\infty}\sum_{l=0}^{p- \nu_{\mathbf{d}}\beta}
      \tilde{\PZc}_{p,l}^{(\beta)}q^{\beta}\hbar^{p- \nu_{\mathbf{d}}\beta-l})\Big\}.
\end{eqnarray}
\end{lemma}
\begin{proof}
Recall
\begin{equation}
      \tilde{\hypergeometricF}(\frac{1}{\hbar},q)=\sum_{\beta=0}^{\infty}q^{\beta}
      \frac{\prod_{k=1}^r\prod_{i=1}^{d_k \beta}(d_k +i\hbar)}{\prod_{j=1}^{\beta}\big((1+j\hbar)^n-1\big)},
\end{equation}
\begin{equation}
      \mathbf{D}^p\tilde{\hypergeometricF}(\frac{1}{\hbar},q)=\sum_{\beta=0}^{\infty}q^{\beta}
      \frac{(1+\beta \hbar)^p\prod_{k=1}^r\prod_{i=1}^{d_k \beta}(d_k +i\hbar)}{\prod_{j=1}^{\beta}\big((1+j\hbar)^n-1\big)},
\end{equation}
\begin{equation}
      \tilde{\hypergeometricF}_{p}(\frac{1}{\hbar},q)=\sum_{\beta=0}^{\infty}\sum_{l=0}^{p- \nu_{\mathbf{d}}\beta}
      \tilde{\PZc}_{p,l}^{(\beta)}q^{\beta}\hbar^{p- \nu_{\mathbf{d}}\beta-l}\mathbf{D}^l \tilde{\hypergeometricF}(\frac{1}{\hbar},q).
\end{equation}
Since $d_l\geq 2$ for $1\leq l\leq r$,  $\mathbf{D}^p\tilde{\hypergeometricF}(\frac{1}{\hbar},q)$ is holomorphic at $\hbar=-d_k$, and the coefficient of each positive power of $q$ in $$\big((1+\hbar)^n-1\big)\big(\mathbf{D}^p\tilde{\hypergeometricF}(\frac{1}{\hbar},q)-1\big)$$ has a zero of the same  order as $\prod_{l=1}^r (d_l+\hbar)$. So
\begin{eqnarray*}
&& \mathrm{Res}_{\hbar=-d_k }\Big\{
\frac{(1+\hbar)^n-1}{\hbar^3\prod_{k=1}^r(d_k+\hbar)}
\frac{\tilde{\hypergeometricF}(\frac{1}{\hbar},q)-\tilde{\hypergeometricF}_p(\frac{1}{\hbar},q)}{\tilde{\hypergeometricF}(\frac{1}{\hbar},q)}\Big\}\\
&=& \mathrm{Res}_{\hbar=-d_k }\Big\{
\frac{(1+\hbar)^n-1}{\hbar^3\prod_{k=1}^r(d_k+\hbar)}
(1-\sum_{\beta=0}^{\infty}\sum_{l=0}^{p- \nu_{\mathbf{d}}\beta}
      \tilde{\PZc}_{p,l}^{(\beta)}q^{\beta}\hbar^{p- \nu_{\mathbf{d}}\beta-l})\Big\}.    
\end{eqnarray*}
Then the residue theorem yields
\begin{eqnarray*}
&& \mathrm{Res}_{\hbar=-\mathbf{d} }\Big\{
\frac{(1+\hbar)^n-1}{\hbar^3\prod_{k=1}^r(d_k+\hbar)}
\frac{\tilde{\hypergeometricF}(\frac{1}{\hbar},q)-\tilde{\hypergeometricF}_p(\frac{1}{\hbar},q)}{\tilde{\hypergeometricF}(\frac{1}{\hbar},q)}\Big\}\\
&=& - \mathrm{Res}_{\hbar=0,\infty }\Big\{
\frac{(1+\hbar)^n-1}{\hbar^3\prod_{k=1}^r(d_k+\hbar)}
(1-\sum_{\beta=0}^{\infty}\sum_{l=0}^{p- \nu_{\mathbf{d}}\beta}
      \tilde{\PZc}_{p,l}^{(\beta)}q^{\beta}\hbar^{p- \nu_{\mathbf{d}}\beta-l})\Big\}.
\end{eqnarray*}
\end{proof}

\begin{lemma}\label{lem-typeII-contribution-SpecializedAtRootsOfUnity-residueAt-hbar=-d-continued}
\begin{eqnarray}\label{eq-typeII-contribution-SpecializedAtRootsOfUnity-residueAt-hbar=-d-hbar=0}
&& \mathrm{Res}_{\hbar=0 }\Big\{
\frac{(1+\hbar)^n-1}{\hbar^3\prod_{k=1}^r(d_k+\hbar)}
(1-\sum_{\beta=0}^{\infty}\sum_{l=0}^{p- \nu_{\mathbf{d}}\beta}
      \tilde{\PZc}_{p,l}^{(\beta)}q^{\beta}\hbar^{p- \nu_{\mathbf{d}}\beta-l})\Big\}\nn\\
&=& \frac{n}{\prod_{k=1}^r d_k}\big((\frac{n-1}{2}-\sum_{k=1}^{r}\frac{1}{d_k})(1-\sum_{\beta=0}^{\infty}\tilde{\PZc}_{p,p- \nu_{\mathbf{d}}\beta}^{(\beta)}q^{\beta})
-\sum_{\beta=0}^{\infty}\tilde{\PZc}_{p,p- \nu_{\mathbf{d}}\beta-1}^{(\beta)}q^{\beta}\big).
\end{eqnarray}
\begin{eqnarray}\label{eq-typeII-contribution-SpecializedAtRootsOfUnity-residueAt-hbar=-d-hbar=infty}
&&\mathrm{Coeff}_{q^{\frac{p-1}{\nu_{\mathbf{d}}}}}\bigg\{
\mathrm{Res}_{\hbar=\infty }\Big\{
\frac{(1+\hbar)^n-1}{\hbar^3\prod_{k=1}^r(d_k+\hbar)}
(1-\sum_{\beta=0}^{\infty}\sum_{l=0}^{p- \nu_{\mathbf{d}}\beta}
      \tilde{\PZc}_{p,l}^{(\beta)}q^{\beta}\hbar^{p- \nu_{\mathbf{d}}\beta-l})\Big\}\bigg\}\nn\\
&=&\mathrm{Res}_{w=0 }\Big\{
\frac{(1+w)^n (\tilde{\PZc}_{1+\nu_{\mathbf{d}} \gamma,0}^{(\gamma)}+\tilde{\PZc}_{1+\nu_{\mathbf{d}} \gamma,1}^{(\gamma)}w)}{w^{n-r}\prod_{k=1}^r(d_k w+1)}
\Big\}.
\end{eqnarray}
\end{lemma}
\begin{proof}
A direct computation yields (\ref{eq-typeII-contribution-SpecializedAtRootsOfUnity-residueAt-hbar=-d-hbar=0}). For (\ref{eq-typeII-contribution-SpecializedAtRootsOfUnity-residueAt-hbar=-d-hbar=infty}), first we have
\begin{eqnarray*}
&& \mathrm{Res}_{\hbar=\infty }\Big\{
\frac{(1+\hbar)^n-1}{\hbar^3\prod_{k=1}^r(d_k+\hbar)}
(1-\sum_{\beta=0}^{\infty}\sum_{l=0}^{p- \nu_{\mathbf{d}}\beta}
      \tilde{\PZc}_{p,l}^{(\beta)}q^{\beta}\hbar^{p- \nu_{\mathbf{d}}\beta-l})\Big\}\\
&=& -\mathrm{Res}_{w=0 }\Big\{
\frac{(1+w)^n-w^n}{w^{n-r-1}\prod_{k=1}^r(d_k w+1)}
(1-\sum_{\beta=0}^{\infty}\sum_{l=0}^{p- \nu_{\mathbf{d}}\beta}
      \tilde{\PZc}_{p,l}^{(\beta)}q^{\beta}w^{-p+\nu_{\mathbf{d}}\beta+l})\Big\}.
\end{eqnarray*}
Then for $p=1+\nu_{\mathbf{d}} \gamma$, by counting the index in the powers of $w$ we have
\begin{eqnarray*}
&&\mathrm{Coeff}_{q^{\frac{p-1}{\nu_{\mathbf{d}}}}}\bigg\{
-\mathrm{Res}_{w=0 }\Big\{
\frac{(1+w)^n-w^n}{w^{n-r-1}\prod_{k=1}^r(d_k w+1)}
(1-\sum_{\beta=0}^{\infty}\sum_{l=0}^{p- \nu_{\mathbf{d}}\beta}
      \tilde{\PZc}_{p,l}^{(\beta)}q^{\beta}w^{-p+\nu_{\mathbf{d}}\beta+l})\Big\}\bigg\}\\
&=&\mathrm{Res}_{w=0 }\Big\{
\frac{(1+w)^n-w^n}{w^{n-r-1}\prod_{k=1}^r(d_k w+1)}
(\tilde{\PZc}_{p,0}^{(\gamma)}w^{-1}+\tilde{\PZc}_{p,1}^{(\gamma)})\Big\}\\
&=&\mathrm{Res}_{w=0 }\Big\{
\frac{(1+w)^n (\tilde{\PZc}_{1+\nu_{\mathbf{d}} \gamma,0}^{(\gamma)}+\tilde{\PZc}_{1+\nu_{\mathbf{d}} \gamma,1}^{(\gamma)}w)}{w^{n-r}\prod_{k=1}^r(d_k w+1)}
\Big\}.
\end{eqnarray*}
\end{proof}

\begin{lemma}\label{lem-typeII-contribution-SpecializedAtRootsOfUnity-residueAt-hbar=0}
\begin{eqnarray}\label{eq-typeII-contribution-SpecializedAtRootsOfUnity-residueAt-hbar=0}
&& \mathrm{Res}_{\hbar=0 }\Big\{
\frac{(1+\hbar)^n-1}{\hbar^3\prod_{k=1}^r(d_k+\hbar)}
\frac{\tilde{\hypergeometricF}(\frac{1}{\hbar},q)-\tilde{\hypergeometricF}_p(\frac{1}{\hbar},q)}{\tilde{\hypergeometricF}(\frac{1}{\hbar},q)}\Big\}\nn\\
&=&\frac{n}{\prod_{k=1}^r d_k}\Big\{ \big(\frac{n-1}{2}-\sum_{k=1}^r \frac{1}{d_k}\big)\big(1-\sum_{\beta=0}^{\infty}\tilde{\PZc}_{p,p- \nu_{\mathbf{d}}\beta}^{(\beta)}q^{\beta}L(q)^{p- \nu_{\mathbf{d}}\beta}\big)\nn\\
&&-\sum_{\beta=0}^{\infty}\tilde{\PZc}_{p,p- \nu_{\mathbf{d}}\beta}^{(\beta)}q^{\beta+1}
\binom{p- \nu_{\mathbf{d}}\beta}{2} L(q)^{p- \nu_{\mathbf{d}}\beta-2}L(q)'\nn\\
&&-\sum_{\beta=0}^{\infty}\tilde{\PZc}_{p,p- \nu_{\mathbf{d}}\beta}^{(\beta)}q^{\beta+1}
(p- \nu_{\mathbf{d}}\beta)  L(q)^{p- \nu_{\mathbf{d}}\beta-1} \frac{\Phi'_0(q)}{\Phi_0(q)}
-\sum_{\beta=0}^{\infty}\tilde{\PZc}_{p,p- \nu_{\mathbf{d}}\beta-1}^{(\beta)}q^{\beta}L(q)^{p- \nu_{\mathbf{d}}\beta-1}\Big\}.
\end{eqnarray}
\end{lemma}
\begin{proof}
By Corollary \ref{cor-regularizability-nonEquivariant}, we have
\begin{eqnarray*}
&& \frac{\tilde{\hypergeometricF}(\frac{1}{\hbar},q)-\tilde{\hypergeometricF}_p(\frac{1}{\hbar},q)}{\tilde{\hypergeometricF}(\frac{1}{\hbar},q)}
=1-\frac{e^{-\frac{\mu(q)}{\hbar}}\tilde{\hypergeometricF}_p(\frac{1}{\hbar},q)}{e^{-\frac{\mu(q)}{\hbar}}\tilde{\hypergeometricF}(\frac{1}{\hbar},q)}\\
&=&1-\frac{e^{-\frac{\mu(q)}{\hbar}}\sum_{\beta=0}^{\infty}\sum_{l=0}^{p- \nu_{\mathbf{d}}\beta}
      \tilde{\PZc}_{p,l}^{(\beta)}q^{\beta}\hbar^{p- \nu_{\mathbf{d}}\beta-l}\mathbf{D}^l \tilde{\hypergeometricF}(\frac{1}{\hbar},q)}{e^{-\frac{\mu(q)}{\hbar}}\tilde{\hypergeometricF}(\frac{1}{\hbar},q)}\\
&=&1-\frac{\sum_{\beta=0}^{\infty}\sum_{l=0}^{p- \nu_{\mathbf{d}}\beta}\tilde{\PZc}_{p,l}^{(\beta)}q^{\beta}\hbar^{p- \nu_{\mathbf{d}}\beta-l}(1+q\frac{d \mu(q)}{dq}+\hbar q\frac{d}{d q})^l\big(e^{-\mu(q)/\hbar}\tilde{\hypergeometricF}(\frac{1}{\hbar},q)\big)}
     {e^{-\frac{\mu(q)}{\hbar}}\tilde{\hypergeometricF}(\frac{1}{\hbar},q)} \\
&=&1-\frac{\sum_{\beta=0}^{\infty}\sum_{l=0}^{p- \nu_{\mathbf{d}}\beta}\tilde{\PZc}_{p,l}^{(\beta)}q^{\beta}\hbar^{p- \nu_{\mathbf{d}}\beta-l}(1+q\frac{d \mu(q)}{dq}+\hbar q\frac{d}{d q})^l
\big(\Phi_0(q)+\hbar \Phi_1(q)+o(\hbar)\big)}
     {\Phi_0(q)+\hbar \Phi_1(q)+o(\hbar)}.
\end{eqnarray*}
Then a direct computation of residues gives (\ref{eq-typeII-contribution-SpecializedAtRootsOfUnity-residueAt-hbar=0}).
\end{proof}

\begin{lemma}\label{lem-typeII-contribution-SpecializedAtRootsOfUnity-residueAt-hbar=infty}
\begin{eqnarray}\label{eq-typeII-contribution-SpecializedAtRootsOfUnity-residueAt-hbar=infty}
&&\mathrm{Coeff}_{q^{\frac{p-1}{\nu_{\mathbf{d}}}}}\bigg\{
\mathrm{Res}_{\hbar=\infty }\Big\{
\frac{(1+\hbar)^n-1}{\hbar^3\prod_{k=1}^r(d_k+\hbar)}
\frac{\tilde{\hypergeometricF}(\frac{1}{\hbar},q)-\tilde{\hypergeometricF}_p(\frac{1}{\hbar},q)}{\tilde{\hypergeometricF}(\frac{1}{\hbar},q)}\Big\}\bigg\}\nn\\
&=&-\mathrm{Coeff}_{q^{\frac{p-1}{\nu_{\mathbf{d}}}}}\bigg\{ \mathrm{Res}_{w=0 }\Big\{
\frac{(1+w)^n}{w^{n-r-1}\prod_{k=1}^r(d_k w+1)}
\frac{\hypergeometricF(w,q)-\hypergeometricF_p(w,q)}{\hypergeometricF(w,q)}\Big\}\bigg\}.
\end{eqnarray}
\end{lemma}
\begin{proof}
We write the residue to compute as
\begin{eqnarray*}
&&\mathrm{Res}_{\hbar=\infty }\Big\{
\frac{(1+\hbar)^n-1}{\hbar^3\prod_{k=1}^r(d_k+\hbar)}
\frac{\tilde{\hypergeometricF}(\frac{1}{\hbar},q)-\tilde{\hypergeometricF}_p(\frac{1}{\hbar},q)}{\tilde{\hypergeometricF}(\frac{1}{\hbar},q)}\Big\}\\
&=&-\mathrm{Res}_{w=0 }\Big\{
\frac{(1+w)^n-w^n}{w^{n-r-1}\prod_{k=1}^r(d_k w+1)}
\frac{\tilde{\hypergeometricF}(w,q)-\tilde{\hypergeometricF}_p(w,q)}{\tilde{\hypergeometricF}(w,q)}\Big\}\\
&=&-\mathrm{Res}_{w=0 }\Big\{
\frac{(1+w)^n}{w^{n-r-1}\prod_{k=1}^r(d_k w+1)}
\frac{\tilde{\hypergeometricF}(w,q)-\tilde{\hypergeometricF}_p(w,q)}{\tilde{\hypergeometricF}(w,q)}\Big\}\\
&&-\mathrm{Res}_{w=0 }\Big\{
\frac{w^{r+1}}{\prod_{k=1}^r(d_k w+1)}
\frac{\tilde{\hypergeometricF}_p(w,q)}{\tilde{\hypergeometricF}(w,q)}\Big\}.
\end{eqnarray*}
Then (\ref{eq-typeII-contribution-SpecializedAtRootsOfUnity-residueAt-hbar=infty}) follows from the following Lemma \ref{lem-typeII-contribution-SpecializedAtRootsOfUnity-residueAt-hbar=infty-coeffOfQ}.
\end{proof}

\begin{lemma}\label{lem-typeII-contribution-SpecializedAtRootsOfUnity-residueAt-hbar=infty-coeffOfQ}
We have
\begin{eqnarray*}
&&\mathrm{Coeff}_{q^{\frac{p-1}{\nu_{\mathbf{d}}}}}\bigg\{
\mathrm{Res}_{w=0 }\Big\{
\frac{(1+w)^n}{w^{n-r-1}\prod_{k=1}^r(d_k w+1)}
\frac{\tilde{\hypergeometricF}(w,q)-\tilde{\hypergeometricF}_p(w,q)}{\tilde{\hypergeometricF}(w,q)}\Big\}\bigg\}\\
&=&\mathrm{Coeff}_{q^{\frac{p-1}{\nu_{\mathbf{d}}}}}\bigg\{ \mathrm{Res}_{w=0 }\Big\{
\frac{(1+w)^n}{w^{n-r-1}\prod_{k=1}^r(d_k w+1)}
\frac{\hypergeometricF(w,q)-\hypergeometricF_p(w,q)}{\hypergeometricF(w,q)}\Big\}\bigg\}.
\end{eqnarray*}
and
\begin{eqnarray*}
&&\mathrm{Coeff}_{q^{\frac{p-1}{\nu_{\mathbf{d}}}}}\bigg\{
\mathrm{Res}_{w=0 }\Big\{
\frac{w^{r+1}}{\prod_{k=1}^r(d_k w+1)}
\frac{\tilde{\hypergeometricF}_p(w,q)}{\tilde{\hypergeometricF}(w,q)}\Big\}\bigg\}=0.
\end{eqnarray*}
\end{lemma}
\begin{proof}
Unraveling the definition of $\tilde{\hypergeometricF}_{p}(w,q)$, we have
\begin{eqnarray*}
\label{eq-def-tildeFp(w,q)-expand}
\tilde{\hypergeometricF}_{p}(w,q)=\frac{1}{w^p}\sum_{\beta=0}^{\infty}\sum_{l=0}^{p- \nu_{\mathbf{d}}\beta}\Big(\tilde{\PZc}_{p,l}^{(\beta)}q^{\beta}w^{\nu_{\mathbf{d}}\beta}
\sum_{\beta'=0}^{\infty}q^{\beta'}w^{\nu_{\mathbf{d}}\beta'}
\frac{(w+ \beta')^p\prod_{k=1}^r\prod_{i=1}^{d_k \beta'}(d_k w+i)}{\prod_{j=1}^{\beta'}\big((w+j)^n-w^n\big)}\Big).
\end{eqnarray*}
It follows that 
\begin{equation*}
      \mathrm{ord}_{w=0}\mathrm{Coeff}_{q^\beta}\tilde{\hypergeometricF}_p(w,q)\geq \nu_{\mathbf{d}}\beta-p.
\end{equation*}
Hence the second equality.
Moreover, the difference between $\tilde{\hypergeometricF}_p(w,q)$ and $\hypergeometricF_p(w,q)$ has order $\geq n-p$ at $w=0$. So the first equality follows.
\end{proof}
Applying Lemma \ref{lem-typeII-contribution-SpecializedAtRootsOfUnity-residueAt-hbar=-d}-\ref{lem-typeII-contribution-SpecializedAtRootsOfUnity-residueAt-hbar=infty}  to RHS of (\ref{eq-typeII-contribution-SpecializedAtRootsOfUnity}), some elementary manipulations yield:
\begin{proposition}\label{prop-typeB-contribution-SpecializedAtRootsOfUnity-final}
The type B Contribution to 
\begin{equation*}
      \sum_{\beta=0}^{\infty}q^{\beta}\langle \sfh_{1+(n-d)\beta}\rangle_{1,\beta}^0
=\sum_{\beta=0}^{\infty}q^{\beta} \mathcal{F}_{\beta,0}
\end{equation*}
is equal to
\begin{eqnarray*}
&& \frac{n}{24}
\sum_{p=0}^{n-1}q^{\frac{p-1}{\nu_{\mathbf{d}}}}
\mathrm{Coeff}_{q^{\frac{p-1}{\nu_{\mathbf{d}}}}}\Big\{ \big(\frac{n-1}{2}-\sum_{k=1}^r \frac{1}{d_k}\big)\Big(1-\sum_{\beta=0}^{\infty}\tilde{\PZc}_{p,p- \nu_{\mathbf{d}}\beta}^{(\beta)}q^{\beta}\big(L(q)^{p- \nu_{\mathbf{d}}\beta}-1\big)\Big)\\
&&-L(q)'\sum_{\beta=0}^{\infty}\tilde{\PZc}_{p,p- \nu_{\mathbf{d}}\beta}^{(\beta)}q^{\beta+1}
\binom{p- \nu_{\mathbf{d}}\beta}{2} L(q)^{p- \nu_{\mathbf{d}}\beta-2}\\
&&-\frac{\Phi'_0(q)}{\Phi_0(q)}\sum_{\beta=0}^{\infty}\tilde{\PZc}_{p,p- \nu_{\mathbf{d}}\beta}^{(\beta)}q^{\beta+1}
(p- \nu_{\mathbf{d}}\beta)  L(q)^{p- \nu_{\mathbf{d}}\beta-1} 
-\sum_{\beta=0}^{\infty}\tilde{\PZc}_{p,p- \nu_{\mathbf{d}}\beta-1}^{(\beta)}q^{\beta}\big(L(q)^{p- \nu_{\mathbf{d}}\beta-1}-1\big)\Big\}\\
&&-\frac{\prod_{k=1}^r d_k}{24}\sum_{p=0}^{n-1}q^{\frac{p-1}{\nu_{\mathbf{d}}}}
\mathrm{Res}_{w=0 }\Big\{
\frac{(1+w)^n (\tilde{\PZc}_{p,0}^{(\frac{p-1}{\nu_{\mathbf{d}}})}+\tilde{\PZc}_{p,1}^{(\frac{p-1}{\nu_{\mathbf{d}}})}w)}{w^{n-r}\prod_{k=1}^r(d_k w+1)}
\Big\}\\
&&-\frac{\prod_{k=1}^r d_k}{24}
\sum_{p=0}^{n-1}q^{\frac{p-1}{\nu_{\mathbf{d}}}}
\mathrm{Coeff}_{q^{\frac{p-1}{\nu_{\mathbf{d}}}}}\bigg\{
\mathrm{Res}_{w=0 }\Big\{
\frac{(1+w)^n}{w^{n-r-1}\prod_{k=1}^r(d_k w+1)}
\frac{\hypergeometricF(w,q)-\hypergeometricF_p(w,q)}{\hypergeometricF(w,q)}\Big\}\bigg\}.
\end{eqnarray*}
\end{proposition}

\section{The difference between reduced and standard genus 1 GW invariants}\label{sec:SvR}
Zinger in \cite{Zin08} proved a formula that expresses the difference between the standard genus 1 GW invariants and the reduced ones:
\begin{theorem}[{\cite[Theorem 1A]{Zin08}}]\label{thm-SvR}
For $\mu=(\mu_1,\dots,\mu_k)\in  H^*(X;\mathbb{Q})^k$, and $\beta\in \mathbb{Z}_{>0}$,
\begin{multline}\label{eq-SvR}
 \langle \mu_1,\dots, \mu_k\rangle_{1,\beta}
-\langle \mu_1,\dots,\mu_k\rangle_{1,\beta}^0\\
= \sum_{m=1}^{\infty}\sum_{J\subset [k]}\Big(
(-1)^{m+|J|}\big(m+|J|;(0_j)_{j\in J}\big)_{[m],J}\\
 \times \sum_{p=0}^{\dim X-2m-|J|}
\mathrm{GW}_{(m,J)}^{\beta} \big(\eta_p,c_{\dim X-2m-|J|-p}(T_X);\mu_1,\dots, \mu_k\big)
\Big).
\end{multline}
\end{theorem}
We refer the reader to \cite[\S 2.1, 2.2]{Zin08} for the definition of the notations in this formula. Here we only recall that $\mathrm{GW}_{(m,J)}^{\beta} \big(\eta_p,c_{i}(T_X);\mu_1,\dots, \mu_k\big)$ is a variant of the usual genus 0 Gromov-Witten invariants. The related moduli space $\mathcal{Z}_{(m,J)}(X,\beta)$ parametrizes morphisms from curves like
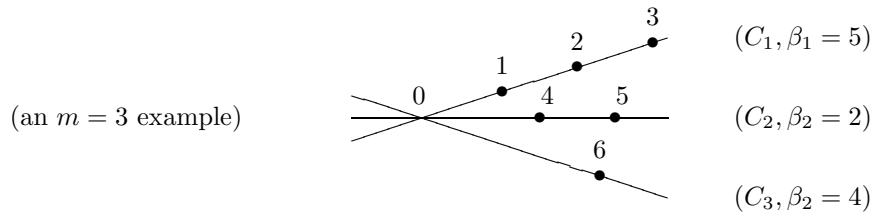
\[
\xy <1cm,0cm>:
(-4,0)*+{(\mbox{an}\ m=3\ \mbox{example})};
 (-1,0); (3.2,0) **@{-}, (5,0)*+{(C_2,\beta_2=2)};
 (-1,-0.3); (3.2,1.065) **@{-},(5,1.065)*+{(C_1,\beta_1=5)};
 (-1,0.3); (3.2,-1.065) **@{-},(5,-1.065)*+{(C_3,\beta_2=4)};
 (-0.1,0.3)*+{0}; 
 (1,0.35)*=0{\bullet};(2,0.675)*=0{\bullet};(3,1)*+{\bullet};
 (1,0.7)*+{1};(2,1.025)*+{2};(3,1.35)*+{3};
 (1.5,0)*+{\bullet};  (2.5,0)*+{\bullet};
 (1.6,0.3)*+{4};(2.6,0.3)*+{5};
  (2.3,-0.7725)*+{\bullet};
  (2.3,-0.42255)*+{6};
\endxy
\]
to $X$, where each line stands for a connected semistable curve $C_i$ of arithmetic genus 0, such that the restriction of the morphism to $C_i$ is a stable map with a strictly positive degree $\beta_i$. A formal definition is given by the cartesian diagram
\begin{equation}\label{eq-def-Mbar(m,J)}
      \xymatrix{
      \Mbar_{(m,J)}(X,\beta) \ar[r] \ar[d] & \bigsqcup_{\begin{subarray}{c}
      \beta_1+\dots+\beta_m=\beta\\
      J_1\sqcup\dots\sqcup J_m=J\end{subarray}}\Big(\prod_{i=1}^m \Mbar_{0,0\sqcup J_i}(X,\beta_i)\Big) \ar[d]^{\prod_i \mathrm{ev}_0}\\
      X \ar[r]^{\Delta} & X^m
            }
\end{equation}
where $\Mbar_{(m,J)}(X,\beta)$ is equipped with a natural $\Sigma_m$-action, and $\mathcal{Z}_{(m,J)}(X,\beta)=\Mbar_{(m,J)}(X,\beta)/\Sigma_m$. Thus the above graph represent an point in $\Mbar_{(3,J)}(X,11)$ where $J=\{1,2,3,4,5,6\}$.  The virtual cycle class $[\Mbar_{(m,J)}(X,\beta)]^{\mathrm{vir}}$ is defined by the Gysin map via $\Delta$, and $[\mathcal{Z}_{(m,J)}(X,\beta)]^{\mathrm{vir}}$ is the induced cycle class via the quotient.

By \cite[Corollary 1.2]{Zin07b},
\begin{equation}\label{eq-intersectionNum-M(1,m,J)-1}
(m+|J|;0,\dots,0)_{m,J}=\frac{m^{|J|}\cdot (m-1)!}{24}.
\end{equation}
Let $\beta\in \mathbb{Z}_{> 0}$ and $a=1+\nu_{\mathbf{d}} \beta$. Then
\begin{subequations}\label{eq-SvR-applied}
\begin{eqnarray}
&&\langle \sfh_{a}\rangle_{1}-      \langle \sfh_{a}\rangle_{1}^{0}\nn\\
&=& \sum_{m=1}^{\infty}\sum_{J\subset [1]}\Big(
(-1)^{m+|J|}\big(m+|J|;0\big)_{[m],J}
  \sum_{p=0}^{\dim X-2m-|J|}
\mathrm{GW}_{(m,J)}^{\beta} \big(\eta_p,c_{\dim X-2m-|J|-p}(T_X);\sfh_a\big)
\Big)\nn\\
&=& \sum_{m=1}^{\infty}\Big(
\frac{(-1)^{m}(m-1)!}{24}
  \sum_{p=0}^{n-1-r-2m}
\mathrm{GW}_{(m,\emptyset)}^{\beta} \big(\eta_p,c_{n-1-r-2m-p}(T_X);\sfh_a\big)
\Big)\nn\\
&&+\sum_{m=1}^{\infty}\Big(
\frac{(-1)^{m+1}m!}{24}
  \sum_{p=0}^{n-2-r-2m}
\mathrm{GW}_{(m,[1])}^{\beta} \big(\eta_p,c_{n-2-r-2m-p}(T_X);\sfh_a\big)
\Big)\nn\\
&=& \frac{1}{24}\sum_{p=2}^{n-1-r}\sum_{m=1}^{2m\leq p}(-1)^m (m-1)!
\mathrm{GW}_{(m,\emptyset)}^{\beta} \big(\eta_{p-2m},c_{n-1-r-p}(T_X);\sfh_a\big)\label{eq-SvR-emptyset}\\
&&+ \frac{1}{24}\sum_{p=2}^{n-2-r}\sum_{m=1}^{2m\leq p}(-1)^{m+1} m!
\mathrm{GW}_{(m,[1])}^{\beta} \big(\eta_{p-2m},c_{n-2-r-p}(T_X);\sfh_a\big).\label{eq-SvR-[1]}
\end{eqnarray}
\end{subequations}

In the following subsections we compute (\ref{eq-SvR-emptyset}) and (\ref{eq-SvR-[1]}) respectively.

\subsection{Computation of the first sum}
In this section we compute (\ref{eq-SvR-emptyset}).
First we compute the Chern class of $T_X$:
\begin{eqnarray*}
&& \mathrm{GW}_{(m,\emptyset)}^{\beta} \big(\eta_{p-2m},c_{n-1-r-p}(T_X);\sfh_a\big)
=\langle \eta_{p-2m} \mathrm{ev}_0^*c_{n-1-r-p}(T_X) \mathrm{ev}_1^* \sfh_a,[\mathcal{Z}_{(m,[1])}(X,\beta)]^{\mathrm{vir}}\rangle\\
&=& \mathrm{Coeff}_{x^{n-1-r-p}}\big\{ \frac{(1+x)^n}{\prod_{k=1}^{r}(1+d_k x)}\big\}
\times\langle \eta_{p-2m} \mathrm{ev}_0^* \sfh_{n-1-r-p}
\mathrm{ev}_1^* \sfh_a,[\mathcal{Z}_{(m,[1])}(X,\beta)]^{\mathrm{vir}}\rangle.
\end{eqnarray*}

\begin{lemma}\label{lem-intermediateSum-0}
\begin{eqnarray}\label{eq-intermediateSum-0}
&& \sum_{\beta=1}^{\infty}q^{\beta}\sum_{m=1}^{2m\leq p}\big((-1)^m (m-1)!
\langle \mathrm{ev}_0^* \sfh_{n-1-r-p} \prod_{i=1}^m \pi_i^* \frac{1}{1-\psi_0}\mathrm{ev}_1^* \sfh_{1+\nu_{\mathbf{d}}\beta}, [\mathcal{Z}_{(m,[1])}(X,\beta)]^{\mathrm{vir}}\rangle\big)\nn\\
&=&-\sum_{\beta>0}q^{\beta}\bigg(\mathrm{Coeff}_{q^{\beta}}\mathrm{Coeff}_{u^{\nu_{\mathbf{d}}\beta}} \mathrm{Coeff}_{(\frac{w}{\hbar})^p}\nn\\
&&\big\{
\frac{\sum_{e=0}^{n-1-r}\sum_{f=0}^{n-2-r}\big(w^{e-f}u^f\sum_{\beta=1}^{\infty}q^{\beta}w^{\nu_{\mathbf{d}}\beta}\int_{[\Mbar_{0,2}(X,\beta)]^{\mathrm{vir}}}
      \frac{1}{\hbar(\hbar-\psi_1)}\mathrm{ev}^*_1 \sfh_{n-1-r-e}\mathrm{ev}^*_2 \sfh_{1+f}\big)}
      {1+(\prod_{k=1}^r d_k)^{-1}\sum_{e=0}^{n-1-r}\big(w^e\sum_{\beta=1}^{\infty}q^{\beta}w^{\nu_{\mathbf{d}}\beta}\int_{[\Mbar_{0,1}(X,\beta)]^{\mathrm{vir}}}
      \frac{1}{\hbar(\hbar-\psi_1)}\mathrm{ev}^*_1 \sfh_{n-1-r-e}\big)}\big\}\bigg).\nn\\
\end{eqnarray}
\end{lemma}
\begin{proof}
By splitting of the virtual cycle class via (\ref{eq-def-Mbar(m,J)}),
\begin{eqnarray*}
&& \sum_{m=1}^{2m\leq p}\big((-1)^m (m-1)!
\langle \mathrm{ev}_0^* \sfh_{n-1-r-p} \prod_{i=1}^m \pi_i^* \frac{1}{1-\psi_0}\mathrm{ev}_1^* \sfh_a, [\mathcal{Z}_{(m,[1])}(X,\beta)]^{\mathrm{vir}}\rangle\big)\\
&=& \sum_{m=1}^{2m\leq p}\frac{(-1)^m}{\big(\prod_{k=1}^r d_k\big)^{m-1}} \sum_{\begin{subarray}{c}\sum_{i=1}^{m}\beta_i=\beta\\\beta_i> 0\end{subarray}}
\sum_{\begin{subarray}{c}\sum_{i=1}^{m}p_i=p\\p_i\geq 0\end{subarray}}
\langle \mathrm{ev}_1^{*}\sfh_{n-1-r-p_1} \frac{1}{1-\psi_1} \mathrm{ev}_2^* \sfh_a,[\Mbar_{0,2}(X,\beta_1)]^{\mathrm{vir}}\rangle\\
&&\times \prod_{i=2}^{m}\langle \mathrm{ev}_1^{*}\sfh_{n-1-r-p_i} \frac{1}{1-\psi_1},[\Mbar_{0,1}(X,\beta_i)]^{\mathrm{vir}}\rangle.
\end{eqnarray*}
So
\begin{eqnarray*}
&& \sum_{\beta=1}^{\infty}q^{\beta}\sum_{m=1}^{2m\leq p}\big((-1)^m (m-1)!
\langle \mathrm{ev}_0^* \sfh_{n-1-r-p} \prod_{i=1}^m \pi_i^* \frac{1}{1-\psi_0}\mathrm{ev}_1^* \sfh_{1+\nu_{\mathbf{d}}\beta}, [\mathcal{Z}_{(m,[1])}(X,\beta)]^{\mathrm{vir}}\rangle\big)\\
&=&\sum_{\beta=1}^{\infty}q^{\beta} \sum_{m=1}^{2m\leq p}\frac{(-1)^m}{\big(\prod_{k=1}^r d_k\big)^{m-1}} \sum_{\begin{subarray}{c}\sum_{i=1}^{m}\beta_i=\beta\\\beta_i> 0\end{subarray}}
\sum_{\begin{subarray}{c}\sum_{i=1}^{m}p_i=p\\p_i\geq 0\end{subarray}}\\
&&\Big(\langle \frac{1}{1-\psi_1}\mathrm{ev}_1^{*}\sfh_{n-1-r-p_1}  \mathrm{ev}_2^* \sfh_{1+\nu_{\mathbf{d}}\beta},[\Mbar_{0,2}(X,\beta_1)]^{\mathrm{vir}}\rangle
\times \prod_{i=2}^{m}\langle\frac{1}{1-\psi_1} \mathrm{ev}_1^{*}\sfh_{n-1-r-p_i} ,[\Mbar_{0,1}(X,\beta_i)]^{\mathrm{vir}}\rangle\Big)\\
&=&\sum_{m=1}^{2m\leq p}\frac{(-1)^m}{\big(\prod_{k=1}^r d_k\big)^{m-1}}\sum_{\begin{subarray}{c}\sum_{i=1}^{m}p_i=p\\p_i\geq 0\end{subarray}}
 \sum_{\beta_1,\dots,\beta_m>0}
 \Big( q^{\beta_1}
\langle \frac{1}{1-\psi_1}\mathrm{ev}_1^{*}\sfh_{n-1-r-p_1}  \mathrm{ev}_2^* \sfh_{1+\nu_{\mathbf{d}}\sum_{i=1}^m \beta_i},[\Mbar_{0,2}(X,\beta_1)]^{\mathrm{vir}}\rangle\\
&&\times \prod_{i=2}^{m}q^{\beta_i}\langle \frac{1}{1-\psi_1}\mathrm{ev}_1^{*}\sfh_{n-1-r-p_i} ,[\Mbar_{0,1}(X,\beta_i)]^{\mathrm{vir}}\rangle\Big).
\end{eqnarray*}
Since $\mathrm{dim} [\Mbar_{0,1}(X,\beta_i)]^{\mathrm{vir}}=n-r-3+\nu_{\mathbf{d}}\beta_i$, the last sum is equal to
\begin{eqnarray*}
&&\sum_{m=1}^{2m\leq p}\frac{(-1)^m}{\big(\prod_{k=1}^r d_k\big)^{m-1}} \sum_{\begin{subarray}{c}\sum_{i=1}^{m}p_i=p\\p_i\geq 0\end{subarray}}
 \sum_{\beta_1,\dots,\beta_m>0}\\
 &&\Big( q^{\beta_1}
\langle \psi_1^{p_1- \nu_{\mathbf{d}}\sum_{i=2}^{m}\beta_i-2} \mathrm{ev}_1^{*}\sfh_{n-1-r-p_1}\mathrm{ev}_2^* \sfh_{1+\nu_{\mathbf{d}}\sum_{i=1}^m \beta_i},[\Mbar_{0,2}(X,\beta_1)]^{\mathrm{vir}}\rangle\\
&&\times \prod_{i=2}^{m}q^{\beta_i}\langle  \mathrm{ev}_1^{*}\sfh_{n-1-r-p_i} \psi_1^{p_i+\nu_{\mathbf{d}}\beta_i-2},[\Mbar_{0,1}(X,\beta_i)]^{\mathrm{vir}}\rangle\Big)\\
&=&-\sum_{\beta>0}q^{\beta}\bigg(\mathrm{Coeff}_{q^{\beta}}\mathrm{Coeff}_{u^{\nu_{\mathbf{d}}\beta}} \mathrm{Coeff}_{(\frac{w}{\hbar})^p}\\
&&\big\{
\frac{\sum_{e=0}^{n-1-r}\sum_{f=0}^{n-2-r}\big(w^{e-f}u^f\sum_{\beta=1}^{\infty}q^{\beta}w^{\nu_{\mathbf{d}}\beta}\int_{[\Mbar_{0,2}(X,\beta)]^{\mathrm{vir}}}
      \frac{1}{\hbar(\hbar-\psi_1)}\mathrm{ev}^*_1 \sfh_{n-1-r-e}\mathrm{ev}^*_2 \sfh_{1+f}\big)}
      {1+(\prod_{k=1}^r d_k)^{-1}\sum_{e=0}^{n-1-r}\big(w^e\sum_{\beta=1}^{\infty}q^{\beta}w^{\nu_{\mathbf{d}}\beta}\int_{[\Mbar_{0,1}(X,\beta)]^{\mathrm{vir}}}
      \frac{1}{\hbar(\hbar-\psi_1)}\mathrm{ev}^*_1 \sfh_{n-1-r-e}\big)}\big\}\bigg).
\end{eqnarray*}
\end{proof}

\begin{lemma}\label{lem-intermediateSum-0-numerator}
\begin{eqnarray}\label{eq-intermediateSum-0-numerator}
&&\sum_{e=0}^{n-1-r}\sum_{f=0}^{n-2-r}\big(w^{e-f}u^f\sum_{\beta=1}^{\infty}q^{\beta}w^{\nu_{\mathbf{d}}\beta}\int_{[\Mbar_{0,2}(X,\beta)]^{\mathrm{vir}}}
      \frac{1}{\hbar(\hbar-\psi_1)}\mathrm{ev}^*_1 \sfh_{n-1-r-e}\mathrm{ev}^*_2 \sfh_{1+f}\big)\nn\\
&=&\prod_{k=1}^r d_k\cdot\frac{w}{\hbar}\sum_{j=1}^{n-1-r}u^{j-1}\times
\begin{cases}
e^{- \mathbf{d}! q\frac{w}{\hbar}}\hypergeometricF_j(\frac{w}{\hbar},q)-1,& \mbox{if}\  \nu_{\mathbf{d}}=1\\
\hypergeometricF_j(\frac{w}{\hbar},q)-1,& \mbox{if}\  \nu_{\mathbf{d}}\geq 2
\end{cases}  \mod w^{n-r}.
\end{eqnarray}
\end{lemma}
\begin{proof}
\begin{eqnarray*}
 &&\sum_{e=0}^{n-1-r}\sum_{f=0}^{n-2-r}\big(w^{e-f}u^{f}\sum_{\beta=1}^{\infty}q^{\beta}w^{\nu_{\mathbf{d}}\beta}\int_{[\Mbar_{0,2}(X,\beta)]^{\mathrm{vir}}}
      \frac{1}{\hbar(\hbar-\psi_1)}\mathrm{ev}^*_1 \sfh_{n-1-r-e}\mathrm{ev}^*_2 \sfh_{1+f}\big)\\
&=& \frac{1}{\hbar} \sum_{j=1}^{n-1-r}u^{j-1}\bigg(\sum_{i=0}^{n-1-r}\sum_{\beta=1}^{\infty}w^{1+i-j}q^{\beta}w^{\nu_{\mathbf{d}}\beta}\int_{[\Mbar_{0,2}(\mathbb{P}^{n-1},\beta)]^{\mathrm{vir}}}
      \frac{\mathbf{e}(\mathcal{V}_0)}{\hbar-\psi_1}\mathrm{ev}^*_1 \sfh_{n-1-r-i}\mathrm{ev}^*_2 \sfh_{j}\bigg)\\  
&=& \frac{\prod_{k=1}^r d_k}{w^r\hbar}\sum_{j=1}^{n-1-r}u^{j-1}w^{1-j}\big(Z_j^*(\hbar,qw^{\nu_{\mathbf{d}}})-\sfh_{j+r}\big)|_{\sfh_{i}\mapsto w^i, 0\leq i\leq n-r-1}.      
\end{eqnarray*}
By \cite[Theorem 5]{PoZ14},
\[
\big(Z_j^*(\hbar,q w^{\nu_{\mathbf{d}}})-\sfh_{j+r}\big)|_{\sfh_{i}\mapsto w^i, 0\leq i\leq n-r-1}
=w^{j+r}\cdot
\begin{cases}
\big(e^{- \mathbf{d}! q\frac{w}{\hbar}}\hypergeometricF_j(\frac{w}{\hbar},q)-1\big),& \mbox{if}\  \nu_{\mathbf{d}}=1\\
\big(\hypergeometricF_j(\frac{w}{\hbar},q)-1\big),& \mbox{if}\  \nu_{\mathbf{d}}\geq 2
\end{cases} \mod w^{n}.
\]
So we get (\ref{eq-intermediateSum-0-numerator}).
\end{proof}

\begin{lemma}\label{lem-intermediateSum-0-denominator}
\begin{eqnarray}\label{eq-intermediateSum-0-denominator}
&& 1+(\prod_{k=1}^r d_k)^{-1}\sum_{e=0}^{n-1-r}\big(w^{e}\sum_{\beta=1}^{\infty}q^{\beta}w^{\nu_{\mathbf{d}}\beta}\int_{[\Mbar_{0,1}(X,\beta)]^{\mathrm{vir}}}
      \frac{1}{\hbar(\hbar-\psi_1)}\mathrm{ev}^*_1 \sfh_{n-1-r-e}\big)\nn\\
&=& \begin{cases}
e^{- \mathbf{d}! q\frac{w}{\hbar}}\hypergeometricF_0(\frac{w}{\hbar},q),& \mbox{if}\  \nu_{\mathbf{d}}=1\\
\hypergeometricF_0(\frac{w}{\hbar},q),& \mbox{if}\  \nu_{\mathbf{d}}\geq 2
\end{cases}  \mod w^{n-r}.
\end{eqnarray}
\end{lemma}
\begin{proof}
\begin{eqnarray*}
&&\sum_{e=0}^{n-1-r}\big(w^{e}\sum_{\beta=1}^{\infty}q^{\beta}w^{\nu_{\mathbf{d}}\beta}\int_{[\Mbar_{0,1}(X,\beta)]^{\mathrm{vir}}}
      \frac{1}{\hbar(\hbar-\psi_1)}\mathrm{ev}^*_1 \sfh_{n-1-r-e}\big)\\
&=&\prod_{k=1}^r d_k\cdot \sum_{i=0}^{n-1-r}\sum_{\beta=1}^{\infty}w^{i}q^{\beta}w^{\nu_{\mathbf{d}}\beta}\int_{[\Mbar_{0,1}(\mathbb{P}^{n-1},\beta)]^{\mathrm{vir}}}
      \frac{\mathbf{e}(\mathcal{V}'_0)}{\hbar(\hbar-\psi)} \mathrm{ev}^*_1 \sfh_{n-1-i}\\
&=& \prod_{k=1}^r d_k\cdot\frac{1}{w^r}\big(Z_0^*(\hbar,q)-\sfh_{r}\big)|_{\sfh_{i}\mapsto w^i, 0\leq i\leq n-r-1}.   
\end{eqnarray*}
Then again by \cite[Theorem 5]{PoZ14}, we get (\ref{eq-intermediateSum-0-denominator}).
\end{proof}

\begin{proposition}\label{prop-SvR-emptyset-generatingFunc}
\begin{eqnarray}\label{eq-SvR-emptyset-generatingFunc}
&& \sum_{\beta>0}q^{\beta}\sum_{p=2}^{n-1-r}\sum_{m=1}^{2m\leq p}(-1)^m (m-1)!
\mathrm{GW}_{(m,\emptyset)}^{\beta} \big(\eta_{p-2m},c_{n-1-r-p}(T_X);\sfh_{1+\nu_{\mathbf{d}}\beta}\big)\nn\\
&=&- \prod_{k=1}^r d_k\cdot\sum_{1\leq \beta\leq \frac{n-2-r}{\nu_{\mathbf{d}}}}q^{\beta}\bigg(\mathrm{Coeff}_{q^{\beta}}\mathrm{Coeff}_{w^{n-r-2}}\Big\{
\frac{(1+w)^n}{\hypergeometricF_0(w,q)\prod_{k=1}^{r}(1+d_k w)}\nn\\
&&\times\begin{cases}
\hypergeometricF_{1+\nu_{\mathbf{d}}\beta}(w,q)-e^{\mathbf{d}! q w},& \mbox{if}\  \nu_{\mathbf{d}}=1\\
\hypergeometricF_{1+\nu_{\mathbf{d}}\beta}(w,q)-1,& \mbox{if}\  \nu_{\mathbf{d}}\geq 2
\end{cases}\Big\}\bigg).
\end{eqnarray}
\end{proposition}
\begin{proof} By Lemma \ref{lem-intermediateSum-0}, \ref{lem-intermediateSum-0-numerator} and \ref{lem-intermediateSum-0-denominator}, we have
\begin{eqnarray}\label{prop-SvR-emptyset-generatingFunc-proof-1}
&& \sum_{\beta>0}q^{\beta}\sum_{p=2}^{n-1-r}\sum_{m=1}^{2m\leq p}(-1)^m (m-1)!
\mathrm{GW}_{(m,\emptyset)}^{\beta} \big(\eta_{p-2m},c_{n-1-r-p}(T_X);\sfh_{1+\nu_{\mathbf{d}}\beta}\big)\nn\\
&=&\sum_{\beta>0}q^{\beta}\sum_{p=2}^{n-1-r}\sum_{m=1}^{2m\leq p}(-1)^m (m-1)!
\Big(\mathrm{Coeff}_{x^{n-1-r-p}}\big\{ \frac{(1+x)^n}{\prod_{k=1}^{r}(1+d_k x)}\big\}\nn\\
&&\times\langle \eta_{p-2m} \mathrm{ev}_0^* \sfh_{n-1-r-p}
\mathrm{ev}_1^* \sfh_{1+\nu_{\mathbf{d}}\beta},[\mathcal{Z}_{(m,[1])}(X,\beta)]^{\mathrm{vir}}\rangle\Big)\nn\\
&=&\sum_{p=2}^{n-1-r}\Big(\mathrm{Coeff}_{x^{n-1-r-p}}\big\{ \frac{(1+x)^n}{\prod_{k=1}^{r}(1+d_k x)}\big\}\nn\\
&&\times\sum_{\beta>0}q^{\beta}\sum_{m=1}^{2m\leq p}(-1)^m (m-1)!\langle \eta_{p-2m} \mathrm{ev}_0^* \sfh_{n-1-r-p}
\mathrm{ev}_1^* \sfh_{1+\nu_{\mathbf{d}}\beta},[\mathcal{Z}_{(m,[1])}(X,\beta)]^{\mathrm{vir}}\rangle\Big)\nn\\
&=&- \prod_{k=1}^r d_k\cdot\sum_{p=2}^{n-1-r}\bigg(\mathrm{Coeff}_{x^{n-1-r-p}}\big\{ \frac{(1+x)^n}{\prod_{k=1}^{r}(1+d_k x)}\big\}\nn\\
&& \times\sum_{\beta>0}q^{\beta}\Big(\mathrm{Coeff}_{q^{\beta}u^{\nu_{\mathbf{d}}\beta}} \mathrm{Coeff}_{(\frac{w}{\hbar})^p} \Big\{\frac{w}{\hbar}\Big(\begin{cases}
e^{- \mathbf{d}! q\frac{w}{\hbar}}\hypergeometricF_0(\frac{w}{\hbar},q),& \mbox{if}\  \nu_{\mathbf{d}}=1\\
\hypergeometricF_0(\frac{w}{\hbar},q),& \mbox{if}\  \nu_{\mathbf{d}}\geq 2
\end{cases} \Big)^{-1}\nn\\
&&\times\sum_{j=1}^{n-1-r}u^{j-1}\begin{cases}
\big(e^{- \mathbf{d}! q\frac{w}{\hbar}}\hypergeometricF_j(\frac{w}{\hbar},q)-1\big),& \mbox{if}\  \nu_{\mathbf{d}}=1\\
\big(\hypergeometricF_j(\frac{w}{\hbar},q)-1\big),& \mbox{if}\  \nu_{\mathbf{d}}\geq 2
\end{cases}   \Big\}\Big)\bigg)\nn\\
&=&- \prod_{k=1}^r d_k\cdot\sum_{p=2}^{n-1-r}\bigg(\mathrm{Coeff}_{x^{n-1-r-p}}\big\{ \frac{(1+x)^n}{\prod_{k=1}^{r}(1+d_k x)}\big\}\nn\\
&& \times\sum_{\beta>0}q^{\beta}\Big(\mathrm{Coeff}_{q^{\beta}u^{\nu_{\mathbf{d}}\beta}}\mathrm{Coeff}_{w^{p-1}} \Big\{
\hypergeometricF_0(w,q)^{-1}
\sum_{j=1}^{n-1-r}u^{j-1}
\begin{cases}
\big(\hypergeometricF_j(w,q)-e^{\mathbf{d}! q w}\big),& \mbox{if}\  \nu_{\mathbf{d}}=1\\
\big(\hypergeometricF_j(w,q)-1\big),& \mbox{if}\  \nu_{\mathbf{d}}\geq 2
\end{cases}   \Big\} \Big)\bigg).\nn\\
\end{eqnarray}
By the following Lemma \ref{lem-regularity-Fp(w,q)AtW=0}, for $p\leq 1$,
\[
\mathrm{Coeff}_{w^{p-1}} \Big\{
F_0(w,q)^{-1}
\sum_{j=1}^{n-1-r}u^{j-1}
\begin{cases}
\big(\hypergeometricF_j(w,q)-e^{\mathbf{d}! q w}\big),& \mbox{if}\  \nu_{\mathbf{d}}=1\\
\big(\hypergeometricF_j(w,q)-1\big),& \mbox{if}\  \nu_{\mathbf{d}}\geq 2
\end{cases}   \Big\}=0.
\]
So the right hand side of (\ref{prop-SvR-emptyset-generatingFunc-proof-1}) is equal to
\begin{eqnarray*}
&&- \prod_{k=1}^r d_k\cdot\sum_{\beta>0}q^{\beta}\Big(\mathrm{Coeff}_{q^{\beta}u^{\nu_{\mathbf{d}}\beta}}\mathrm{Coeff}_{w^{n-r-2}}\Big\{
\frac{(1+w)^n}{\hypergeometricF_0(w,q)\prod_{k=1}^{r}(1+d_k w)}\\
&&\times\sum_{j=1}^{n-1-r}u^{j-1}
\begin{cases}
\hypergeometricF_j(w,q)-e^{\mathbf{d}! q w},& \mbox{if}\  \nu_{\mathbf{d}}=1\\
\hypergeometricF_j(w,q)-1,& \mbox{if}\  \nu_{\mathbf{d}}\geq 2
\end{cases}\Big\}\Big)\\
&=& \mbox{RHS of (\ref{eq-SvR-emptyset-generatingFunc})}.
\end{eqnarray*}
\end{proof}

\begin{lemma}\label{lem-regularity-Fp(w,q)AtW=0}
For $p\geq 0$, $\hypergeometricF_p(w,q)$ is regular at $w=0$, i.e. $\mathrm{Coeff}_{w^i}\big\{\hypergeometricF_p(w,q)\big\}=0$ for $i<0$.
\end{lemma}
\begin{proof}
By the definition (\ref{eq-def-numbers-c}) of the coefficients $c^{(b)}_{i,j}$, 
\begin{eqnarray*}
\hypergeometricF_{p}(w,q)&=&\sum_{\beta=0}^{\infty}\sum_{l=0}^{p- \nu_{\mathbf{d}}\beta}\frac{\tilde{\PZc}_{p,l}^{(\beta)}q^{\beta}}{w^{p- \nu_{\mathbf{d}}\beta-l}}\mathbf{D}^l \hypergeometricF(w,q)\\
&=&\sum_{\beta=0}^{\infty}\sum_{l=0}^{p- \nu_{\mathbf{d}}\beta}\big(\frac{\tilde{\PZc}_{p,l}^{(\beta)}q^{\beta}}{w^{p- \nu_{\mathbf{d}}\beta}}
\sum_{\beta'=0}^{\infty}\sum_{s=0}^{\infty}\PZc_{l,s}^{(\beta')}w^{s}q^{\beta'}w^{\nu_{\mathbf{d}}\beta'}\big)\\
&=&\sum_{\beta=0}^{\infty}\sum_{\beta'=0}^{\infty}\sum_{s=0}^{\infty}\sum_{l=0}^{p- \nu_{\mathbf{d}}\beta}\tilde{\PZc}_{p,l}^{(\beta)}
\PZc_{l,s}^{(\beta')}q^{\beta+\beta'}w^{\nu_{\mathbf{d}}(\beta'+\beta)+s-p}\\
&=&\sum_{\beta=0}^{\infty}\sum_{s=0}^{\infty}q^{\beta}w^{\nu_{\mathbf{d}}\beta+s-p}
(\sum_{\begin{subarray}{c} \beta_1+\beta_2=\beta\\ \beta_1,\beta_2\geq 2\end{subarray}}\sum_{l=0}^{p- \nu_{\mathbf{d}}\beta_1}\tilde{\PZc}_{p,l}^{(\beta_1)}\PZc_{l,s}^{(\beta_2)}).
\end{eqnarray*}
For $s\leq p-\nu_{\mathbf{d}}\beta$, by (\ref{eq-def-numbers-tildec}), $\sum_{\begin{subarray}{c} \beta_1+\beta_2=\beta\\ \beta_1,\beta_2\geq 2\end{subarray}}\sum_{l=0}^{p- \nu_{\mathbf{d}}\beta_1}\tilde{\PZc}_{p,l}^{(\beta_1)}\PZc_{l,s}^{(\beta_2)}=0$. So we are done.
\end{proof}

\subsection{Computation of the second sum}
In this section we compute (\ref{eq-SvR-[1]}) on $X$. First we have
\begin{eqnarray*}
&& \mathrm{GW}_{(m,[1])}^{\beta} \big(\eta_{p-2m},c_{n-2-r-p}(T_X);\sfh_a\big)
=\langle \eta_{p-2m} \mathrm{ev}_0^*\big(c_{n-2-r-p}(T_X)  \sfh_a\big),[\mathcal{Z}_{(m,\emptyset)}(X,\beta)]^{\mathrm{vir}}\rangle\\
&=& \mathrm{Coeff}_{x^{n-2-r-p}}\big\{ \frac{(1+x)^n}{\prod_{k=1}^{r}(1+d_k x)}\big\}
\langle \eta_{p-2m} \mathrm{ev}_0^* \sfh_{n-2-r-p+a},[\mathcal{Z}_{(m,\emptyset)}(X,\beta)]^{\mathrm{vir}}\rangle.
\end{eqnarray*}

\begin{lemma}\label{lem-intermediateSum-1}
\begin{eqnarray}\label{eq-intermediateSum-1}
&& \sum_{\beta=1}^{\infty}q^{\beta}\sum_{m=1}^{2m\leq p}(-1)^{m+1} m!
\langle \mathrm{ev}_0^* \sfh_{n-1-p+\nu_{\mathbf{d}}\beta} \prod_{i=1}^m \pi_i^* \frac{1}{1-\psi_0}, [\mathcal{Z}_{(m,\emptyset)}(X,\beta)]^{\mathrm{vir}}\rangle\nn\\
&=&\mathrm{Coeff}_{\hbar^{-p}} \left\{\frac{\sum_{e} \sum_{\beta=1}^{\infty}q^{\beta}
\int_{[\Mbar_{0,1}(X,\beta)]^{\mathrm{vir}}}
      \frac{1}{\hbar(\hbar-\psi)} \mathrm{ev}^*_1 \sfh_{n-3-e}}
{1+\big(\prod_{k=1}^r d_k\big)^{-1}\sum_{e} \sum_{\beta=1}^{\infty}q^{\beta}
\int_{[\Mbar_{0,1}(X,\beta)]^{\mathrm{vir}}}
      \frac{1}{\hbar(\hbar-\psi)} \mathrm{ev}^*_1 \sfh_{n-3-e}}
\right\}.
\end{eqnarray}
\end{lemma}
\begin{proof}
By splitting of the virtual cycle class via (\ref{eq-def-Mbar(m,J)}),
\begin{eqnarray*}
&&\sum_{m=1}^{2m\leq p}\big((-1)^{m+1}m!
\langle \mathrm{ev}_0^* \sfh_{n-2-p+a} \prod_{i=1}^m \pi_i^* \frac{1}{1-\psi_0}, [\mathcal{Z}_{(m,\emptyset)}(X,\beta)]^{\mathrm{vir}}\rangle\big)\\
&=& \sum_{m=1}^{2m\leq p}\frac{(-1)^{m+1}}{\big(\prod_{k=1}^r d_k\big)^{m-1}} \sum_{\begin{subarray}{c}\sum_{i=1}^{m}\beta_i=\beta\\\beta_i> 0\end{subarray}}
\sum_{\begin{subarray}{c}\sum_{i=1}^{m}p_i=p-a+1\\p_i\geq 0\end{subarray}}
\Big( \prod_{i=1}^{m}\langle \mathrm{ev}_1^{*}\sfh_{n-1-p_i} \frac{1}{1-\psi_1},[\Mbar_{0,1}(X,\beta_i)]^{\mathrm{vir}}\rangle\Big).
\end{eqnarray*}
So
\begin{eqnarray*}
&& \sum_{\beta=1}^{\infty}q^{\beta}\sum_{m=1}^{2m\leq p}(-1)^{m+1} m!
\langle \mathrm{ev}_0^* \sfh_{n-1-p+\nu_{\mathbf{d}}\beta} \prod_{i=1}^m \pi_i^* \frac{1}{1-\psi_0}, [\mathcal{Z}_{(m,\emptyset)}(X,\beta)]^{\mathrm{vir}}\rangle\\
&=&\sum_{\beta=0}^{\infty}\sum_{m=1}^{2m\leq p}\frac{(-1)^{m+1}}{\big(\prod_{k=1}^r d_k\big)^{m-1}}\sum_{\begin{subarray}{c}\sum_{i=1}^{m}\beta_i=\beta\\\beta_i> 0\end{subarray}}
\sum_{\begin{subarray}{c}\sum_{i=1}^{m}p_i=p- \nu_{\mathbf{d}}\beta\\p_i\geq 0\end{subarray}}q^{\beta_i}
\Big( \prod_{i=1}^{m}\langle \mathrm{ev}_1^{*}\sfh_{n-1-p_i} \frac{1}{1-\psi_1},[\Mbar_{0,1}(X,\beta_i)]^{\mathrm{vir}}\rangle\Big)\\
&=&\sum_{\beta=1}^{\infty}\sum_{m=1}^{2m\leq p}\frac{(-1)^{m+1}}{\big(\prod_{k=1}^r d_k\big)^{m-1}} \sum_{\begin{subarray}{c}\sum_{i=1}^{m}\beta_i=\beta\\\beta_i > 0\end{subarray}}
\sum_{\begin{subarray}{c}\sum_{i=1}^{m}p_i=p\\p_i\geq 0\end{subarray}}q^{\beta_i}
\Big( \prod_{i=1}^{m}\langle \mathrm{ev}_1^{*}\sfh_{n-1-p_i+\nu_{\mathbf{d}} \beta_i} \frac{1}{1-\psi_1},[\Mbar_{0,1}(X,\beta_i)]^{\mathrm{vir}}\rangle\Big)\\
&=& \sum_{m=1}^{2m\leq p}\frac{(-1)^{m+1}}{\big(\prod_{k=1}^r d_k\big)^{m-1}}\sum_{\begin{subarray}{c}\sum_{i=1}^{m}p_i=p\\p_i\geq 0\end{subarray}} \prod_{i=1}^m \check{Z}_{p_i-2}(q)\\
&=&\mathrm{Coeff}_{w^p} \left\{\frac{\sum_{i=0}^{\infty}\check{Z}_i(q)w^{i+2}}{1+\big(\prod_{k=1}^r d_k\big)^{-1}\sum_{i=0}^{\infty}\check{Z}_i(q)w^{i+2}}
\right\},
\end{eqnarray*}
where
\begin{eqnarray*}
      \check{Z}_{i}(q)&:=&\sum_{\beta=1}^{\infty}q^{\beta}\int_{[\Mbar_{0,1}(X,\beta)]^{\mathrm{vir}}}
      \frac{1}{1-\psi}\mathrm{ev}^*_1 \sfh_{n-3-i+\nu_{\mathbf{d}}\beta}\\
&=&\sum_{\beta=1}^{\infty}q^{\beta}\int_{[\Mbar_{0,1}(X,\beta)]^{\mathrm{vir}}}
      \psi^{i} \mathrm{ev}^*_1 \sfh_{n-3-i+\nu_{\mathbf{d}}\beta}.
\end{eqnarray*}
So we get (\ref{eq-intermediateSum-1}).
\end{proof}

\begin{lemma}\label{lem-intermediateSum-1-fraction}
For $2\leq p\leq n-2-r$,
\begin{eqnarray}\label{eq-intermediateSum-1-fraction}
&& \mathrm{Coeff}_{\hbar^{-p}} \left\{\frac{\sum_{e} \sum_{\beta=1}^{\infty}q^{\beta}
\int_{[\Mbar_{0,1}(X,\beta)]^{\mathrm{vir}}}
      \frac{1}{\hbar(\hbar-\psi)} \mathrm{ev}^*_1 \sfh_{n-3-e}}
{1+\big(\prod_{k=1}^r d_k\big)^{-1}\sum_{e} \sum_{\beta=1}^{\infty}q^{\beta}
\int_{[\Mbar_{0,1}(X,\beta)]^{\mathrm{vir}}}
      \frac{1}{\hbar(\hbar-\psi)} \mathrm{ev}^*_1 \sfh_{n-3-e}}
\right\}\nn\\
&=& \big(\prod_{k=1}^r d_k\big)\cdot\mathrm{Coeff}_{w^{p}} \big\{1-
\begin{cases}
e^{\mathbf{d}! q w}\hypergeometricF(w,q)^{-1},& \mbox{if}\ \nu_{\mathbf{d}}=1\\
\hypergeometricF(w,q)^{-1},& \mbox{if}\ \nu_{\mathbf{d}}\geq 2
\end{cases}\big\}.
\end{eqnarray}
\end{lemma}
\begin{proof}
Define 
\begin{equation}
      Z_p^*(\hbar,q):=\sfh_{r+p}+\sum_{\beta=1}^{\infty}q^{\beta} \mathrm{ev}_{1*}\Big(
      \frac{\mathbf{e}(\mathcal{V}''_0)\mathrm{ev}_2^* \sfh_{r+p}}{\hbar-\psi_1}\Big),
\end{equation}
which is a non-equivariant version of $\mathcal{Z}_p^*(x,\hbar,q)$. Then by the string equation we have
\begin{eqnarray*}
Z_0^*(\hbar,q)
= \sfh_{r}\Big(1+\sum_{\beta=1}^{\infty}q^{\beta} 
\sum_{i=0}^{n-1-r}\sfh_{i} 
\int_{[\Mbar_{0,1}(\mathbb{P}^{n-1},\beta)]}\frac{\mathbf{e}(\mathcal{V}'_0)\mathrm{ev}_1^* \sfh_{n-1-i}}{\hbar(\hbar-\psi_1)}\Big).
\end{eqnarray*}
By \cite[Theorem 5]{PoZ14},
\begin{eqnarray*}
&& Z_0^*(\hbar,q\sfh_{\nu_{\mathbf{d}}})=\begin{cases}
e^{-\mathbf{d}!q\frac{\sfh}{\hbar}}\sfh_{r}\cdot \hypergeometricF(\frac{\sfh}{\hbar},q),& \mbox{if}\ \nu_{\mathbf{d}}=1,\\
\sfh_{r}\cdot \hypergeometricF(\frac{\sfh}{\hbar},q),& \mbox{if}\ \nu_{\mathbf{d}}\geq 2.
\end{cases}
\end{eqnarray*}
This means that, when  $\nu_{\mathbf{d}}\geq 2$,
\begin{eqnarray*}
&& 1+\sum_{\beta=1}^{\infty}q^{\beta}w^{\nu_{\mathbf{d}}\beta} 
\sum_{i=0}^{n-1-r}w^{i} 
\int_{[\Mbar_{0,1}(\mathbb{P}^{n-1},\beta)]^{\mathrm{vir}}}\frac{\mathbf{e}(\mathcal{V}'_0)\mathrm{ev}_1^* \sfh_{n-1-i}}{\hbar(\hbar-\psi_1)}\\
&=& \sum_{\beta=0}^{\infty}q^{\beta}(\frac{w}{\hbar})^{\nu_{\mathbf{d}}\beta}\frac{\prod_{k=1}^r\prod_{i=1}^{d_k \beta}(\frac{d_k w}{\hbar}+i)}{\prod_{j=1}^{\beta}(\frac{w}{\hbar}+j)^n}\ \mod w^{n-r}.
\end{eqnarray*}
Then for $2\leq p\leq n-2-r$,
\begin{eqnarray*}
&& \mathrm{Coeff}_{\hbar^{-p}} \left\{\frac{\sum_{e} \sum_{\beta=1}^{\infty}q^{\beta}
\int_{[\Mbar_{0,1}(\mathbb{P}^{n-1},\beta)]}
      \frac{\mathbf{e}(\mathcal{V}'_0)}{\hbar(\hbar-\psi)} \mathrm{ev}^*_1 \sfh_{n-3-e}}
{1+\sum_{e} \sum_{\beta=1}^{\infty}q^{\beta}
\int_{[\Mbar_{0,1}(\mathbb{P}^{n-1},\beta)]}
      \frac{\mathbf{e}(\mathcal{V}'_0)}{\hbar(\hbar-\psi)} \mathrm{ev}^*_1 \sfh_{n-3-e}}
\right\}\\
&=& \mathrm{Coeff}_{\hbar^{-p}} \big\{1-\frac{1}{\hypergeometricF(\frac{1}{\hbar},q)}\big\}
=\mathrm{Coeff}_{w^{p}} \big\{1-\frac{1}{\hypergeometricF(w,q)}\big\}.
\end{eqnarray*}
So we complete the proof in the cases $\nu_{\mathbf{d}}\geq 2$.
The case $\nu_{\mathbf{d}}=1$ is similar.
\end{proof}

\begin{proposition}\label{prop-SvR-[1]-generatingFunc}
\begin{eqnarray}\label{eq-SvR-[1]-generatingFunc}
&&\sum_{\beta>0}\sum_{p=2}^{n-2-r}\sum_{m=1}^{2m\leq p}(-1)^{m+1} m!
\mathrm{GW}_{(m,[1])}^{\beta} \big(\eta_{p-2m},c_{n-2-r-p}(T_X);\sfh_{1+\nu_{\mathbf{d}}\beta}\big)\\
&=&\prod_{k=1}^r d_k\cdot \mathrm{Coeff}_{w^{n-2-r}}\Big\{
\frac{(1+w)^n}{\hypergeometricF_0(w,q)\prod_{k=1}^{r}(1+d_k w)}
\begin{cases}
\big(\hypergeometricF_0(w,q)-e^{\mathbf{d}! q w}\big),& \mbox{if}\  \nu_{\mathbf{d}}=1\\
\big(\hypergeometricF_0(w,q)-1\big),& \mbox{if}\  \nu_{\mathbf{d}}\geq 2
\end{cases}\Big\}.
\end{eqnarray}
\end{proposition}
\begin{proof}
By Lemma \ref{lem-intermediateSum-1} and \ref{lem-intermediateSum-1-fraction} we have
\begin{eqnarray}\label{prop-SvR-[1]-generatingFunc-proof-1}
&&\sum_{\beta>0}\sum_{p=2}^{n-2-r}\sum_{m=1}^{2m\leq p}(-1)^{m+1} m!
\mathrm{GW}_{(m,[1])}^{\beta} \big(\eta_{p-2m},c_{n-2-r-p}(T_X);\sfh_{1+\nu_{\mathbf{d}}\beta}\big)\nn\\
&=&\prod_{k=1}^r d_k\cdot \sum_{p=2}^{n-2-r}\Big(\mathrm{Coeff}_{x^{n-2-r-p}}\big\{ \frac{(1+x)^n}{\prod_{k=1}^{r}(1+d_k x)}\big\}\nn\\
&&\times \mathrm{Coeff}_{w^{p}} \big\{1-
\begin{cases}
e^{\mathbf{d}! q w}\hypergeometricF_0(w,q)^{-1},& \mbox{if}\ \nu_{\mathbf{d}}=1\\
\hypergeometricF_0(w,q)^{-1},& \mbox{if}\ \nu_{\mathbf{d}}\geq 2
\end{cases}\big\}\Big).
\end{eqnarray}
By the definition of $\hypergeometricF_0(w,q)$ we have
\begin{equation*}
      \mathrm{Coeff}_{w} \big\{1-
\begin{cases}
e^{\mathbf{d}! q w}\hypergeometricF_0(w,q)^{-1},& \mbox{if}\ \nu_{\mathbf{d}}=1\\
\hypergeometricF_0(w,q)^{-1},& \mbox{if}\ \nu_{\mathbf{d}}\geq 2
\end{cases}\big\}=0.
\end{equation*}
So the right hand side of (\ref{prop-SvR-[1]-generatingFunc-proof-1}) is equal to
\begin{eqnarray*}
&&\prod_{k=1}^r d_k\cdot \bigg(\mathrm{Coeff}_{w^{n-2-r}}\Big\{
\frac{(1+w)^n}{\hypergeometricF_0(w,q)\prod_{k=1}^{r}(1+d_k w)}
\begin{cases}
\big(\hypergeometricF_0(w,q)-e^{\mathbf{d}! q w}\big),& \mbox{if}\  \nu_{\mathbf{d}}=1\\
\big(\hypergeometricF_0(w,q)-1\big),& \mbox{if}\  \nu_{\mathbf{d}}\geq 2
\end{cases}\Big\}\\
&&-\mathrm{Coeff}_{x^{n-3-r}}\big\{ \frac{(1+x)^n}{\prod_{k=1}^{r}(1+d_k x)}\big\}
\times \mathrm{Coeff}_{w} \big\{1-
\begin{cases}
e^{\mathbf{d}! q w}\hypergeometricF_0(w,q)^{-1},& \mbox{if}\ \nu_{\mathbf{d}}=1\\
\hypergeometricF_0(w,q)^{-1},& \mbox{if}\ \nu_{\mathbf{d}}\geq 2
\end{cases}\big\}\bigg)\\
&=&\prod_{k=1}^r d_k\cdot \mathrm{Coeff}_{w^{n-2-r}}\Big\{
\frac{(1+w)^n}{\hypergeometricF_0(w,q)\prod_{k=1}^{r}(1+d_k w)}
\begin{cases}
\big(\hypergeometricF_0(w,q)-e^{\mathbf{d}! q w}\big),& \mbox{if}\  \nu_{\mathbf{d}}=1\\
\big(\hypergeometricF_0(w,q)-1\big),& \mbox{if}\  \nu_{\mathbf{d}}\geq 2
\end{cases}\Big\}.
\end{eqnarray*}

\end{proof}

\subsection{The final formula}
\begin{proposition}\label{prop-SvR-final}
For a Fano complete intersection $X\in \mathbb{P}^{n-1}$ of multidegree $\mathbf{d}$ and Fano index $\nu_{\mathbf{d}}$, and for $1\leq \beta\leq \frac{n-2-r}{\nu_{\mathbf{d}}}$,
\begin{gather}\label{eq-SvR-final}
\langle \sfh_{1+\nu_{\mathbf{d}}\beta}\rangle_{1}-\langle \sfh_{1+\nu_{\mathbf{d}}\beta}\rangle_{1}^0\nn\\
=\frac{\prod_{k=1}^r d_k}{24}\cdot \mathrm{Coeff}_{q^{\beta}} \mathrm{Coeff}_{w^{n-2-r}}\Big\{
\frac{(1+w)^n\big(\hypergeometricF_0(w,q)-\hypergeometricF_{1+\nu_{\mathbf{d}}\beta}(w,q)\big)}{\hypergeometricF_0(w,q)\prod_{k=1}^{r}(1+d_k w)}\Big\}.
\end{gather}
\end{proposition}
\begin{proof}
Applying Proposition \ref{prop-SvR-emptyset-generatingFunc} and \ref{prop-SvR-[1]-generatingFunc} to (\ref{eq-SvR-applied}) we get
\begin{eqnarray*}
&& \sum_{\beta>0}q^{\beta}\langle \sfh_{1+\nu_{\mathbf{d}}\beta}\rangle_{1}
-\sum_{\beta>0}q^{\beta}\langle \sfh_{1+\nu_{\mathbf{d}}\beta}\rangle_{1}^0\\
&=&- \frac{\prod_{k=1}^r d_k}{24}\cdot\sum_{1\leq \beta\leq \frac{n-2-r}{\nu_{\mathbf{d}}}}q^{\beta}\bigg(\mathrm{Coeff}_{q^{\beta}}\mathrm{Coeff}_{w^{n-r-2}}\\
&&\Big\{
\frac{(1+w)^n}{\hypergeometricF_0(w,q)\prod_{k=1}^{r}(1+d_k w)}
\begin{cases}
\hypergeometricF_{1+\nu_{\mathbf{d}}\beta}(w,q)-e^{\mathbf{d}! q w},& \mbox{if}\  \nu_{\mathbf{d}}=1\\
\hypergeometricF_{1+\nu_{\mathbf{d}}\beta}(w,q)-1,& \mbox{if}\  \nu_{\mathbf{d}}\geq 2
\end{cases}\Big\}\bigg)\\
&&+\frac{\prod_{k=1}^r d_k}{24}\cdot  \mathrm{Coeff}_{w^{n-2-r}}\Big\{
\frac{(1+w)^n}{\hypergeometricF_0(w,q)\prod_{k=1}^{r}(1+d_k w)}
\begin{cases}
\hypergeometricF_0(w,q)-e^{\mathbf{d}! q w},& \mbox{if}\  \nu_{\mathbf{d}}=1\\
\hypergeometricF_0(w,q)-1,& \mbox{if}\  \nu_{\mathbf{d}}\geq 2
\end{cases}\Big\}.
\end{eqnarray*}
\end{proof}

\begin{remark}
  The difference (\ref{eq-SvR-final}) vanishes, a priori,  when $\beta> \lfloor\frac{n-2-r}{\nu_{\mathbf{d}}}\rfloor$.
\end{remark}

\begin{theorem}\label{thm-genus1-GWInv-Fano}
Let $X$ be a smooth complete intersection of multidegree $\mathbf{d}$ in $\mathbb{P}^{n-1}$, with Fano index $\nu_{\mathbf{d}}\geq 1$. 
 Then for $0\leq b\leq \frac{n-1}{\nu_{\mathbf{d}}}$,
\begin{eqnarray}\label{eq-genus1-GWInv}
&& \langle \sfh_{1+\nu_{\mathbf{d}}b}\rangle_{1,b}\nn\\
&=& \frac{1}{2}
\mathrm{Coeff}_{q^{b}}\Big\{
\frac{\Theta_{1+\nu_{\mathbf{d}}b}^{(0)}(q)\big(\sum_{p=0}^{n-1-r}      \Theta_{p}^{(1)}(q)\Theta_{n-1-r-p}^{(0)}(q)
+\sum_{p=1}^{r}   \Theta_{n-p}^{(1)}(q)\Theta_{n-1-r+p}^{(0)}(q)\big)}{\Phi_0(q)}\Big\}\nn\\
&&+\frac{n}{24}
\mathrm{Coeff}_{q^{b}}\bigg\{ \big(\frac{n-1}{2}-\sum_{k=1}^r \frac{1}{d_k}\big)\Big(1-\sum_{\beta=0}^{\infty}\tilde{\PZc}_{1+\nu_{\mathbf{d}}b,1+\nu_{\mathbf{d}}b- \nu_{\mathbf{d}}\beta}^{(\beta)}q^{\beta}\big(L(q)^{1+\nu_{\mathbf{d}}b- \nu_{\mathbf{d}}\beta}-1\big)\Big)\nn\\
&&-L(q)'\sum_{\beta=0}^{\infty}\tilde{\PZc}_{1+\nu_{\mathbf{d}}b,1+\nu_{\mathbf{d}}b- \nu_{\mathbf{d}}\beta}^{(\beta)}q^{\beta+1}
\binom{1+\nu_{\mathbf{d}}b- \nu_{\mathbf{d}}\beta}{2} L(q)^{1+\nu_{\mathbf{d}}b- \nu_{\mathbf{d}}\beta-2}\nn\\
&&-\frac{\Phi'_0(q)}{\Phi_0(q)}\sum_{\beta=0}^{\infty}\tilde{\PZc}_{1+\nu_{\mathbf{d}}b,1+\nu_{\mathbf{d}}b- \nu_{\mathbf{d}}\beta}^{(\beta)}q^{\beta+1}
(1+\nu_{\mathbf{d}}b- \nu_{\mathbf{d}}\beta)  L(q)^{1+\nu_{\mathbf{d}}b- \nu_{\mathbf{d}}\beta-1} \nn\\
&&-\sum_{\beta=0}^{\infty}\tilde{\PZc}_{1+\nu_{\mathbf{d}}b,1+\nu_{\mathbf{d}}b- \nu_{\mathbf{d}}\beta-1}^{(\beta)}q^{\beta}\big(L(q)^{1+\nu_{\mathbf{d}}b- \nu_{\mathbf{d}}\beta-1}-1\big)\bigg\}\nn\\
&&-\frac{\prod_{k=1}^r d_k}{24}
\mathrm{Res}_{w=0 }\Big\{
\frac{(1+w)^n (\tilde{\PZc}_{1+\nu_{\mathbf{d}}b,0}^{(b)}+\tilde{\PZc}_{1+\nu_{\mathbf{d}}b,1}^{(b)}w)}{w^{n-r}\prod_{k=1}^r(d_k w+1)}
\Big\}.
\end{eqnarray}
\end{theorem} 
\begin{proof}
For $1\leq b\leq \frac{n-1}{\nu_{\mathbf{d}}}$, (\ref{eq-genus1-GWInv}) follows from Proposition \ref{prop-SvR-final}, 
\ref{prop-typeA-contribution-SpecializedAtRootsOfUnity-final} and \ref{prop-typeB-contribution-SpecializedAtRootsOfUnity-final}. By definition and (\ref{eq-tildePZCoefficient-KroneckerProperty}) we have
\begin{equation}
      \Theta_p^{(1)}(0)=0.
\end{equation}
So when $b=0$ (\ref{eq-genus1-GWInv}) reads
\begin{eqnarray}\label{eq-genus1-GWInv-deg0}
\langle \sfh_{1}\rangle_{1,0}
= -\frac{\prod_{k=1}^r d_k}{24}
\mathrm{Res}_{w=0 }\Big\{
\frac{(1+w)^n }{w^{n-r-1}\prod_{k=1}^r(d_k w+1)}
\Big\}.
\end{eqnarray}
This is also true by a genus 1 axiom of GW invariants.
\end{proof}


\section{Some formulae related to hypergeometric series associated with Fano complete intersections}\label{sec:formulae-hypergeometricSeries-Fano}
In this section, we give closed formulae for the series $\mu(q)$, $L(q)$, $\Phi_0(q)$, and $\Phi_1(q)$ that appear in (\ref{eq-genus1-GWInv}).
\subsection{Series and their formulae}
Let
\begin{equation}
      \varphi(z,q):=\sum_{\beta=0}^{\infty}q^{\beta}
      \frac{\prod_{k=1}^r\prod_{i=1}^{d_k \beta}(d_k +iz)}{\prod_{j=1}^{\beta}\big((1+jz)^n-1\big)}.
\end{equation}
Then
\begin{equation*}
      \mu(q)=\mathrm{Res}_{z=0}\big(\ln \varphi(z,q)\big),
\end{equation*}
\begin{equation*}
      Q(z,q)=e^{-\frac{\mu(q)}{z}}  \varphi(z,q),
\end{equation*}
and
\begin{equation*}
       \Phi_i(q)=\frac{\partial^i}{\partial z^i}Q(z,q)|_{z=0}.
\end{equation*}
Recall $L(q)=1+q \mu'(q)$.
We define
\begin{equation*}
      L_1(q):=q L'(q)=q \mu'(q)+q^2 \mu''(q),
\end{equation*}
\begin{equation*}
      L_2(q):=q L'_1(q)=q \mu'(q)+3q^2 \mu''(q)+q^3 \mu'''(q).
\end{equation*}
Note that
\begin{equation*}
      Q(z,0)=\varphi(z,0)=1.
\end{equation*}
So
\begin{equation*}
      \Phi_0(q)\in 1+q \mathbb{C}[[q]],\ \Phi_1(q)\in q \mathbb{C}[[q]].
\end{equation*}
Define an operator
\begin{equation*}
      D_z:=1+zq \frac{\partial}{\partial q}.
\end{equation*}
Induction on $n$ yields
\begin{eqnarray}\label{eq-formula-Dz}
&&  e^{-\frac{\mu(q)}{z}}\circ D_z^n\circ e^{\frac{\mu(q)}{z}}Q(z,q)\nn\\
&=& L^n\Phi_0(q)
+z\Big(L^n\Phi_1(q)+nqL^{n-1}\Phi'_0(q)
+\binom{n}{2}L^{n-2}L_1\Phi_0(q)\Big)\nn\\
&&+z^2\Big(L^n\Phi_2(q)
+\binom{n}{2}L^{n-2}L_1\Phi_1(q)+nq L^{n-1}\Phi'_1(q)\nn\\
&&+\big(\frac{n(n-1)(n-2)(n-3)}{8} 
L^{n-4}L_1^2+\frac{n(n-1)(n-2)}{6} L^{n-3}L_2\big)\Phi_0(q)\nn\\
&&+\big(\binom{n}{2}q L^{n-2}
+\frac{n(n-1)(n-2)}{2}qL^{n-3}L_1\big) \Phi'_0(q)
+\binom{n}{2}q^2L^{n-2}\Phi''_0(q)\Big)+o(z^2).
\end{eqnarray}

\begin{lemma}\label{eq-algebraicEq-L0}
\begin{equation}\label{eq-diffEq-mu-completeIntersection}
 L(q)^n- q \mathbf{d}^{\mathbf{d}} L(q)^{|\mathbf{d}|}=1.
 \end{equation}
\end{lemma}
\begin{proof}
We have
\begin{equation}\label{eq-hypergeometricEq-varphi}
      \big(D_z^n-q\prod_{k=1}^r\prod_{i=1}^{d_k }(d_k D_z+iz)\big)\varphi(z,q)=\varphi(z,q).
\end{equation}
We expand
\begin{eqnarray*}
&& q\prod_{k=1}^r\prod_{i=1}^{d_k }(d_k D_z+iz)\varphi(z,q)\\
&=&q \Big( \big(\prod_{k=1}^r d_k^{d_k}\big) D_z^{|\mathbf{d}|}+(\prod_{k=1}^{r}d_k^{d_k})(\sum_{k=1}^r\sum_{i=1}^{d_k}\frac{i}{d_k}\big) zD_z^{|\mathbf{d}|-1}\\
&&+\mathbf{d}^{\mathbf{d}}\cdot \frac{1}{2}(\sum_{k_1=1}^r\sum_{k_2=1}^{r}\sum_{i_1=1}^{k_1}\sum_{i_2=1}^{k_2}\frac{i_1 i_2}{d_{k_1}d_{k_2}}
-\sum_{k=1}^r \sum_{i=1}^k \frac{i^2}{d_k^2})z^2D_z^{|\mathbf{d}|-2}
+o(z^2)\Big)\varphi(z,q)\\
&=&q \Big( \mathbf{d}^{\mathbf{d}} D_z^{|\mathbf{d}|}+\frac{\mathbf{d}^{\mathbf{d}}(|\mathbf{d}|+r)}{2}zD_z^{|\mathbf{d}|-1}+\mathbf{d}^{\mathbf{d}}\vartheta(\mathbf{d})z^2 D_z^{|\mathbf{d}|-2}
+o(z^2)\Big)\varphi(z,q),
\end{eqnarray*}
where
\begin{eqnarray*}
\vartheta(\mathbf{d})&:=&\frac{1}{2}\left(\sum_{k_1=1}^r\sum_{k_2=1}^{r}\sum_{i_1=1}^{k_1}\sum_{i_2=1}^{k_2}\frac{i_1 i_2}{d_{k_1}d_{k_2}}
-\sum_{k=1}^r \sum_{i=1}^k \frac{i^2}{d_k^2}\right)\\
&=& \frac{1}{8}|\mathbf{d}|^2+(\frac{r}{4}-\frac{1}{6})|\mathbf{d}|+(\frac{r^2}{8}-\frac{r}{4})
-\frac{1}{12}\sum_{k=1}^r \frac{1}{d_k}.
\end{eqnarray*}
Then we apply (\ref{eq-formula-Dz}) to LHS of (\ref{eq-hypergeometricEq-varphi}), and extract the coefficient of $z^0$. We obtain
\begin{equation*}
      \big(L(q)^n- q \mathbf{d}^{\mathbf{d}} L(q)^{|\mathbf{d}|}\big)\Phi_0(q)=\Phi_0(q).
\end{equation*}
But $\Phi_0(q)=1+O(q)$ is invertible in $\mathbb{C}[[q]]$, so we get (\ref{eq-diffEq-mu-completeIntersection}).
\end{proof}

\begin{proposition}\label{prop-formula-L0}
We have
\begin{equation}\label{eq-formula-L0(q)}
      L(q)=\sum_{k=0}^{\infty}\frac{\prod_{i=1}^{k-1}(k|\mathbf{d}|+1-in)}{k!}(\frac{\mathbf{d}^{\mathbf{d}}q}{n})^k,
\end{equation}
and
\begin{eqnarray}
\mu(q)
=\sum_{k=1}^{\infty}\frac{(\mathbf{d}^{\mathbf{d}})^k\prod_{i=1}^{k-1}(k|\mathbf{d}|+1-in)}{k! k n^{k}}q^k.
\end{eqnarray}
\end{proposition}
\begin{proof}
Let $\tau(q):=L(\frac{n q}{\mathbf{d}^{\mathbf{d}}})$.
Then  (\ref{eq-diffEq-mu-completeIntersection}) reads
\begin{equation*}
      \tau(q)^n-n q\tau(q)^{|\mathbf{d}|}=1,
\end{equation*}
and (\ref{eq-formula-L0(q)}) is equivalent to 
\begin{equation*}
      \tau(q)=\sum_{k=0}^{\infty}\frac{\prod_{i=1}^{k-1}(k|\mathbf{d}|+1-in)}{k!}q^k.
\end{equation*}
Wo show this by  Lagrange inversion. For brevity, we treat only the case $r=1$, namely $\mathbf{d}=(d)$. For a function $f(z)$ holomorphic in $U\backslash\{z_0\}$, where $U$ is an open neighborhood of $z_0$, we denote by 
\begin{equation*}
      \oint_{z_0}f(z)\mathrm{d}z
\end{equation*}
the integral along a circle in $U$ centered at $z_0$, in the counterclockwise direction. Since
\[
\mathrm{d}q=\frac{(n^2-nd)\tau^{n+d-1}+nd \tau^{d-1}}{n^2 \tau^{2d}}\mathrm{d} \tau,
\]
and that 1 is a simple root of $\tau^n-1=0$, we have
\begin{eqnarray*}
 \frac{1}{2\pi \sqrt{-1}}\oint_{0}\frac{\tau(q)}{q^{k+1}}\mathrm{d}q
=\frac{ n^k(n-d)}{2\pi \sqrt{-1}}\oint_{1} \frac{\tau^{dk}}{(\tau^n-1)^{k}}\mathrm{d} \tau
+\frac{n^{k+1}}{2\pi \sqrt{-1}}\oint_{1} \frac{\tau^{dk}}{(\tau^n-1)^{k+1}}\mathrm{d} \tau.
\end{eqnarray*}
Let $\sigma=\tau^n$. Then
\begin{eqnarray*}
\frac{1}{2\pi \sqrt{-1}}\oint_{1} \frac{\tau^{a}}{(\tau^n-1)^{k}}\mathrm{d} \tau
=\frac{1}{2n\pi \sqrt{-1}}\oint_{1} \frac{\sigma^{\frac{a-n+1}{n}}}{(\sigma-1)^{k}}\mathrm{d} \sigma
= \frac{1}{n}\binom{\frac{a-n+1}{n}}{k-1}.
\end{eqnarray*}
Hence
\begin{eqnarray*}
&& \frac{1}{2\pi \sqrt{-1}}\oint_{0}\frac{\tau(q)}{q^{k+1}}\mathrm{d}q
= n^{k-1}(n-d)\binom{\frac{dk-n+1}{n}}{k-1}+ n^k \binom{\frac{dk-n+1}{n}}{k}\\
&=& \frac{\prod_{i=1}^{k-1}(kd+1-in)}{k!}.
\end{eqnarray*}
\end{proof}

\begin{proposition}\label{prop-formula-Phi0-Phi1}
\begin{enumerate}
\item[\mbox{(i)}]
      \begin{equation}\label{eq-formula-Phi0}
      \Phi_0(q)=L(q)^{\frac{r+1}{2}}\cdot \big(1+ \mathbf{d}^{\mathbf{d}}(1-\frac{|\mathbf{d}|}{n})q L(q)^{|\mathbf{d}|}\big)^{-\frac{1}{2}}.
      \end{equation}
\item[(ii)]
\begin{eqnarray}\label{eq-formula-Phi1}
\Phi_1(q)&=&\frac{3 r^2-2  |\mathbf{d}|\sum_{k=1}^{r}\frac{1}{d_k} -1}{24|\mathbf{d}|}L(q)^{\frac{r-1}{2}}\big(L(q)-1\big) \big(1+ \mathbf{d}^{\mathbf{d}}(1-\frac{|\mathbf{d}|}{n})q L(q)^{|\mathbf{d}|}\big)^{-\frac{1}{2}}\nn\\
&&+\frac{L(q)^{\frac{r-1}{2}}\cdot \big(1+ \mathbf{d}^{\mathbf{d}}(1-\frac{|\mathbf{d}|}{n})q L(q)^{|\mathbf{d}|}\big)^{-\frac{7}{2}} }{24|\mathbf{d}|n^3}
\times\Big(|\mathbf{d}| ^3   \left(|\mathbf{d}|  n-|\mathbf{d}| -3 r^2+1\right)L(q)\nn\\
&&+|\mathbf{d}| ^2 n  \left(2 |\mathbf{d}| ^2-6 |\mathbf{d}|  n-6 |\mathbf{d}|  r+3 n^2+6 n r+n+3 r^2-1\right)L(q)^n\nn\\
&&+3 |\mathbf{d}| ^2  (n-|\mathbf{d}| ) \left(|\mathbf{d}|  n-|\mathbf{d}| -3 r^2+1\right) L(q)^{n+1}\nn\\
&&+|\mathbf{d}|  n (n-|\mathbf{d}| ) \left(4 |\mathbf{d}| ^2-5 |\mathbf{d}|  n-12 |\mathbf{d}|  r-2 n^2+6 n r+n+6 r^2-2\right)L(q)^{2 n} \nn\\
&&+3 |\mathbf{d}|   (n-|\mathbf{d}| )^2  \left(|\mathbf{d}|  n-|\mathbf{d}| -3 r^2+1\right)L(q)^{2 n+1}\nn\\
&&+n (n-|\mathbf{d}| )^2 \left(2 |\mathbf{d}| ^2+|\mathbf{d}|  n-6 |\mathbf{d}|  r+3 r^2-1\right)L(q)^{3 n}\nn \\
&&+ (n-|\mathbf{d}| )^3  \left(|\mathbf{d}|  n-|\mathbf{d}| -3 r^2+1\right)L(q)^{3 n+1}\Big). 
\end{eqnarray}
\end{enumerate}

\end{proposition}

The proof is given in the next section.

\subsection{ODEs from hypergeometric series}
As in the proof of Lemma \ref{eq-algebraicEq-L0}, extracting the coefficient of $z$ and $z^2$ in (\ref{eq-hypergeometricEq-varphi}), we obtain respectively
\begin{eqnarray}\label{eq-ODE-Phi0-0}
&& \Big(L(q)^n\Phi_1(q)+nqL(q)^{n-1}\Phi'_0(q)+\binom{n}{2}L(q)^{n-2}L_1(q)\Phi_0(q)\Big)\nn\\
&&-q \mathbf{d}^\mathbf{d}\Big(L(q)^{|\mathbf{d}|}\Phi_1(q)+q|\mathbf{d}|L(q)^{|\mathbf{d}|-1}\Phi'_0(q)+\binom{|\mathbf{d}|}{2}L(q)^{|\mathbf{d}|-2}L_1(q)\Phi_0(q)\Big)\nn\\
&&-\frac{q \mathbf{d}^{\mathbf{d}}(|\mathbf{d}|+1)}{2}\big(1+q \mu'(q)\big)^{|\mathbf{d}|-1}\Phi_0(q)=\Phi_1(q),
\end{eqnarray}
and
\begin{eqnarray}\label{eq-ODE-Phi1-0}
&& \Big(L^n\Phi_2(q)
+\binom{n}{2}L^{n-2}L_1\Phi_1(q)+nq L^{n-1}\Phi'_1(q)\nn\\
&&+\big(\frac{n(n-1)(n-2)(n-3)}{8} 
L^{n-4}L_1^2+\frac{n(n-1)(n-2)}{6} L^{n-3}L_2\big)\Phi_0(q)\nn\\
&&+\big(\binom{n}{2}q L^{n-2}
+\frac{n(n-1)(n-2)}{2}qL^{n-3}L_1\big) \Phi'_0(q)
+\binom{n}{2}q^2L^{n-2}\Phi''_0(q)\Big)\nn\\
&&-q \mathbf{d}^{\mathbf{d}} \Big(L^{|\mathbf{d}|}\Phi_2(q)
+\binom{|\mathbf{d}|}{2}L^{|\mathbf{d}|-2}L_1\Phi_1(q)+|\mathbf{d}|q L^{|\mathbf{d}|-1}\Phi'_1(q)\nn\\
&&+\big(\frac{|\mathbf{d}|(|\mathbf{d}|-1)(|\mathbf{d}|-2)(|\mathbf{d}|-3)}{8} 
L^{|\mathbf{d}|-4}L_1^2+\frac{|\mathbf{d}|(|\mathbf{d}|-1)(|\mathbf{d}|-2)}{6} L^{|\mathbf{d}|-3}L_2\big)\Phi_0(q)\nn\\
&&+\big(\binom{|\mathbf{d}|}{2}q L^{|\mathbf{d}|-2}
+\frac{|\mathbf{d}|(|\mathbf{d}|-1)(|\mathbf{d}|-2)}{2}qL^{|\mathbf{d}|-3}L_1\big) \Phi'_0(q)
+\binom{|\mathbf{d}|}{2}q^2L^{|\mathbf{d}|-2}\Phi''_0(q)\Big)\nn\\
&&-q \frac{\mathbf{d}^{\mathbf{d}}(|\mathbf{d}|+r)}{2}\Big(L^{|\mathbf{d}|-1}\Phi_1(q)+(|\mathbf{d}|-1)qL^{|\mathbf{d}|-2}\Phi'_0(q)
+\binom{|\mathbf{d}|-1}{2}L^{|\mathbf{d}|-3}L_1\Phi_0(q)\Big)\nn\\
&&-q\mathbf{d}^{\mathbf{d}}\vartheta(\mathbf{d})L^{|\mathbf{d}|-2}\Phi_0(q)= \Phi_2(q).
\end{eqnarray}
Applying (\ref{eq-algebraicEq-L0}) to (\ref{eq-ODE-Phi0-0}), the sum of terms involving $\Phi_1(q)$ vanishes, and we obtain an ODE for $\Phi_0(q)$. Similarly from (\ref{eq-ODE-Phi1-0}) we obtain an ODE for $\Phi_1(q)$.
Set
\begin{equation}
      A(q):=\binom{n}{2}L^{n-2}L_1-q \mathbf{d}^{\mathbf{d}}\binom{|\mathbf{d}|}{2}L^{|\mathbf{d}|-2}L_1 
-q \frac{\mathbf{d}^{\mathbf{d}}(|\mathbf{d}|+r)}{2}L^{|\mathbf{d}|-1},
\end{equation}
\begin{equation}
B(q):=nq L^{n-1}-q^2 \mathbf{d}^{\mathbf{d}}|\mathbf{d}| L^{|\mathbf{d}|-1},
\end{equation}
and
\begin{eqnarray}
C(q)&:=& \Big(\frac{n(n-1)(n-2)(n-3)}{8} L^{n-4}L_1^2+\frac{n(n-1)(n-2)}{6} L^{n-3}L_2\nn\\
&&-q \mathbf{d}^{\mathbf{d}}\big(\frac{|\mathbf{d}|(|\mathbf{d}|-1)(|\mathbf{d}|-2)(|\mathbf{d}|-3)}{8} 
L^{|\mathbf{d}|-4}L_1^2+\frac{|\mathbf{d}|(|\mathbf{d}|-1)(|\mathbf{d}|-2)}{6} L^{|\mathbf{d}|-3}L_2\big)\nn\\
&&-q \frac{\mathbf{d}^{\mathbf{d}}(|\mathbf{d}|+r)}{2}\binom{|\mathbf{d}|-1}{2}L^{|\mathbf{d}|-3}L_1
-q\mathbf{d}^{\mathbf{d}}\vartheta(\mathbf{d})L^{|\mathbf{d}|-2}\Big)\Phi_0(q)\nn\\
&&+\Big(\binom{n}{2}q L^{n-2}
+\frac{n(n-1)(n-2)}{2}qL^{n-3}L_1\nn\\
&&- \frac{\mathbf{d}^{\mathbf{d}}|\mathbf{d}|(|\mathbf{d}|-1)(|\mathbf{d}|-2)}{2}q^2 L^{|\mathbf{d}|-3}L_1
-\frac{ \mathbf{d}^\mathbf{d}(2|\mathbf{d}|+r)(|\mathbf{d}|-1)}{2}q^2 L^{|\mathbf{d}|-2}\Big) \Phi'_0(q)\nn\\
&&+\big(\binom{n}{2}L^{n-2}-q \mathbf{d}^{\mathbf{d}} \binom{|\mathbf{d}|}{2}L^{|\mathbf{d}|-2})q^2\Phi''_0(q).
\end{eqnarray}
Then
\begin{equation*}
      A(q)\Phi_0(q)+B(q)\Phi'_0(q)=0,
\end{equation*}
and
\begin{equation*}
      A(q)\Phi_1(q)+B(q)\Phi'_1(q)+C(q)=0.
\end{equation*}
Since $\Phi_0(0)=1$ and $\Phi_1(0)=0$, we obtain
\begin{equation}
      \Phi_0(q)=e^{-\int_0^q \frac{A(q)}{B(q)}\mathrm{d}q},
\end{equation}
and
\begin{equation}\label{eq-integral-Phi1}
      \Phi_1(q)=-\Phi_0(q)\int_0^q \frac{C(q)}{B(q)\Phi_0(q)}\mathrm{d}q.
\end{equation}
We are left to integrate the integrals in the following two subsections.

\begin{proof}[Proof of Proposition \ref{prop-formula-Phi0-Phi1} (i)]
Our strategy is to manipulate the integrand into a total differential. In the process we need to write the functions $L_i(q)$ and their derivatives in terms of $L(q)$. So we need repeatedly use
\begin{eqnarray*}
L'(q)=\frac{\mathbf{d}^{\mathbf{d}} L^{|\mathbf{d}|+1}}{n+q \mathbf{d}^{\mathbf{d}}(n-|\mathbf{d}|)L^{|\mathbf{d}|}},
\end{eqnarray*}
\begin{equation}
      L_1(q)=\frac{\mathbf{d}^{\mathbf{d}} qL^{|\mathbf{d}|+1}}{n+ (n-|\mathbf{d}|)\mathbf{d}^{\mathbf{d}}qL^{|\mathbf{d}|}},
\end{equation}
and 
\begin{eqnarray*}
 L_2(q)=\frac{n^2+(n^2+n) \mathbf{d}^{\mathbf{d}}q L^{|\mathbf{d}|}+(n-|\mathbf{d}|) (\mathbf{d}^{\mathbf{d}})^2 q^2  L^{2|\mathbf{d}|}}{\big(n+ \mathbf{d}^{\mathbf{d}}(n-|\mathbf{d}|)q L^{|\mathbf{d}|}\big)^2}L_1,
\end{eqnarray*}
which are all derived from  (\ref{eq-diffEq-mu-completeIntersection}). Then computations yield
\begin{eqnarray}\label{eq-A(q)-inTermsOfL0}
A(q)= \frac{q \mathbf{d}^\mathbf{d} L^{|\mathbf{d}|-1}}{2}\Big(n-r-1
-\frac{n|\mathbf{d}|}{n+q \mathbf{d}^{\mathbf{d}}(n-|\mathbf{d}|)L^{|\mathbf{d}|}}\Big),
\end{eqnarray}
\begin{eqnarray}\label{eq-B(q)-inTermsOfL0}
B(q)= \frac{q\big(n + \mathbf{d}^{\mathbf{d}}(n-|\mathbf{d}|)qL^{|\mathbf{d}|}\big)}{L},
\end{eqnarray}
\begin{eqnarray}
 \frac{\partial}{\partial q}\big(n+ \mathbf{d}^{\mathbf{d}}(n-|\mathbf{d}|)qL^{|\mathbf{d}|}\big)
= \mathbf{d}^{\mathbf{d}} n(n-|\mathbf{d}|)\cdot \frac{(1+ \mathbf{d}^{\mathbf{d}} q L^{|\mathbf{d}|})L^{|\mathbf{d}|}}{n+ \mathbf{d}^{\mathbf{d}}(n-|\mathbf{d}|)q L^{|\mathbf{d}|}},
\end{eqnarray}
\begin{eqnarray*}
 \frac{\frac{\partial}{\partial q}\big(n+ \mathbf{d}^{\mathbf{d}}(n-|\mathbf{d}|)qL^{|\mathbf{d}|}\big)}{n+ \mathbf{d}^{\mathbf{d}}(n-|\mathbf{d}|)qL^{|\mathbf{d}|}}
= \frac{n L'}{L}
-\frac{\mathbf{d}^{\mathbf{d}}|\mathbf{d}|n L^{|\mathbf{d}|}}{\big(n+ \mathbf{d}^{\mathbf{d}}(n-|\mathbf{d}|) q L^{|\mathbf{d}|}\big)^2}.
\end{eqnarray*}
Then
\begin{eqnarray*}
\frac{A(q)}{B(q)}&=& \frac{\mathbf{d}^{\mathbf{d}} L^{|\mathbf{d}|}}{2}\cdot\Big(\frac{n-r-1}{n+\mathbf{d}^{\mathbf{d}}(n-|\mathbf{d}|) q L^{|\mathbf{d}|}}
-\frac{n|\mathbf{d}|}{\big(n+ \mathbf{d}^{\mathbf{d}}(n-|\mathbf{d}|)q L^{|\mathbf{d}|}\big)^2}\Big)\\
&=& -\frac{r+1}{2}  \frac{L'}{L}+\frac{1}{2}\frac{\frac{\partial}{\partial q}\big(n+ \mathbf{d}^{\mathbf{d}}(n-|\mathbf{d}|)qL^{|\mathbf{d}|}\big)}{n+ \mathbf{d}^{\mathbf{d}}(n-|\mathbf{d}|)qL^{|\mathbf{d}|}}.
\end{eqnarray*}
Hence
\[
-\int_0^q \frac{A(q)}{B(q)}dq=\frac{r+1}{2}\ln L-\frac{1}{2}\ln \big(n+ \mathbf{d}^{\mathbf{d}}(n-|\mathbf{d}|)qL^{|\mathbf{d}|}\big),
\]
and (\ref{eq-formula-Phi0}) follows.
\end{proof}

\begin{proof}[Proof of Proposition \ref{prop-formula-Phi0-Phi1} (ii)]
 From (\ref{eq-A(q)-inTermsOfL0}) and (\ref{eq-B(q)-inTermsOfL0}) we can compute the derivatives
 \begin{eqnarray*}
&A'(q)&=  \frac{\mathbf{d}^{\mathbf{d}} L^{|\mathbf{d}|-1}+\mathbf{d}^{\mathbf{d}} (|\mathbf{d}|-1) L^{|\mathbf{d}|-2}L_1}{2}\Big(n-r-1
-\frac{n|\mathbf{d}|}{n+ \mathbf{d}^{\mathbf{d}}(n-|\mathbf{d}|)qL^{|\mathbf{d}|}}\Big)\\
&&+ \frac{\mathbf{d}^{\mathbf{d}} q L^{|\mathbf{d}|-1}}{2} \cdot \frac{n|\mathbf{d}|\big( \mathbf{d}^{\mathbf{d}}(n-|\mathbf{d}|)L^{|\mathbf{d}|}+ \mathbf{d}^{\mathbf{d}}(n-|\mathbf{d}|)|\mathbf{d}| L^{|\mathbf{d}|-1}L_1
\big)}{\big(n+\mathbf{d}^{\mathbf{d}}(n-|\mathbf{d}|) q L^{|\mathbf{d}|}\big)^2},
\end{eqnarray*}
\begin{eqnarray*}
B'(q)&=& \frac{n +2 \mathbf{d}^{\mathbf{d}}(n-|\mathbf{d}|)q L^{|\mathbf{d}|}+\mathbf{d}^{\mathbf{d}}|\mathbf{d}|(n-|\mathbf{d}|)q L^{|\mathbf{d}|-1} L_1}{L}\\
&&-\frac{\big(n + \mathbf{d}^{\mathbf{d}}(n-|\mathbf{d}|)q L^{|\mathbf{d}|}\big)L_1}{L^2},
\end{eqnarray*}
and thus also
\begin{equation*}
      \frac{\Phi'_0(q)}{\Phi_0(q)}=-\frac{A(q)}{B(q)},
\end{equation*}
and
\begin{eqnarray*}
&&\frac{\Phi''_0(q)}{\Phi_0(q)}=\frac{-A'(q)B(q)+A(q)B'(q)+A(q)^2}{B(q)^2}.
\end{eqnarray*}
We put
\begin{equation}
      w(q)=\mathbf{d}^{\mathbf{d}}q  L(q)^{|\mathbf{d}|}.
\end{equation}
Then one can see that $\frac{C(q)}{\Phi_0(q)}$ can be written as a polynomial of $w(q)$. Moreover, using (\ref{eq-algebraicEq-L0}) and
\begin{eqnarray*}
 \frac{\partial}{\partial q}\big(\mathbf{d}^{\mathbf{d}}q  L(q)^{|\mathbf{d}|}\big)
= \frac{ n \mathbf{d}^{\mathbf{d}} L^{|\mathbf{d}|}(1+ \mathbf{d}^{\mathbf{d}} q L^{|\mathbf{d}|})}{n+\mathbf{d}^{\mathbf{d}} (n-|\mathbf{d}|) qL^{|\mathbf{d}|}}
\end{eqnarray*}
we get
\begin{eqnarray*}
&&\int_0^q \frac{C(q)}{B(q)\Phi_0(q)}\mathrm{d}q\\
&=& \frac{1}{24n} \int_0^w \frac{1}{(1+w)^{\frac{1}{n}}}\times\Big( -\frac{15 n^4 |\mathbf{d}|^2}{(n w+n-|\mathbf{d}| w)^4}
+\frac{n^3 \left(14 n |\mathbf{d}|-12 r |\mathbf{d}|+14 |\mathbf{d}|^2+3 |\mathbf{d}|\right)}{(n w+n-|\mathbf{d}| w)^3}\\
&&+\frac{-2 n^4+6 n^3 r-8 n^3 |\mathbf{d}|-n^3-3 n^2 r^2+6 n^2 r |\mathbf{d}|+6 n^2 r-2 n^2 |\mathbf{d}|^2-7 n^2 |\mathbf{d}|+2 n^2}{(n w+n-|\mathbf{d}| w)^2}\\
&&+\frac{n^2 |\mathbf{d}|+3 n r^2-6 n r |\mathbf{d}|+2 n |\mathbf{d}|^2-n}{|\mathbf{d}| (n w+n-|\mathbf{d}| w)}
+\frac{-3 r^2+2  |\mathbf{d}|\sum_{k=1}^{r}\frac{1}{d_k} +1}{|\mathbf{d}| (w+1)}\Big)
\mathrm{d}w.
\end{eqnarray*}
To compute this integration, one can make a further change of variables $x=(1+w)^{\frac{1}{n}}-1$. Then  a lengthy computation yields
\begin{eqnarray*}
&& \int_0^{q} \frac{C(q)}{B(q)\Phi_0(q)} \mathrm{d}q\\
&=& \frac{1}{24n}\bigg(\frac{1}{|\mathbf{d}|  L(q) \left(|\mathbf{d}| +(n-|\mathbf{d}| ) L(q)^n\right)^3}\Big(-|\mathbf{d}| ^3 n  \left(|\mathbf{d}|  n-|\mathbf{d}| -3 r^2+1\right)L(q)\\
&&-|\mathbf{d}| ^2 n^2  \left(2 |\mathbf{d}| ^2-6 |\mathbf{d}|  n-6 |\mathbf{d}|  r+3 n^2+6 n r+n+3 r^2-1\right)L(q)^n\\
&&-3 |\mathbf{d}| ^2 n (n-|\mathbf{d}| ) \left(|\mathbf{d}|  n-|\mathbf{d}| -3 r^2+1\right) L(q)^{n+1}\\
&&-|\mathbf{d}|  n^2 (n-|\mathbf{d}| ) \left(4 |\mathbf{d}| ^2-5 |\mathbf{d}|  n-12 |\mathbf{d}|  r-2 n^2+6 n r+n+6 r^2-2\right)L(q)^{2 n} \\
&&-3 |\mathbf{d}|  n (n-|\mathbf{d}| )^2  \left(|\mathbf{d}|  n-|\mathbf{d}| -3 r^2+1\right)L(q)^{2 n+1}\\
&&-n^2 (n-|\mathbf{d}| )^2 \left(2 |\mathbf{d}| ^2+|\mathbf{d}|  n-6 |\mathbf{d}|  r+3 r^2-1\right)L(q)^{3 n} \\
&&-n (n-|\mathbf{d}| )^3  \left(|\mathbf{d}|  n-|\mathbf{d}| -3 r^2+1\right)L(q)^{3 n+1}
\Big)\\
&&+\frac{-3 r^2+2  |\mathbf{d}|\sum_{k=1}^{r}\frac{1}{d_k} +1}{|\mathbf{d}|}\cdot n\big(1-L(q)^{-1}\big)\bigg).
\end{eqnarray*}
So by (\ref{eq-formula-Phi0}) and (\ref{eq-integral-Phi1}) we obtain (\ref{eq-formula-Phi1}).

\end{proof}

\section{Further simplifications of  \texorpdfstring{$A(q)$}{A(q)} }\label{sec:simplifications-A(q)}
Recall the expression (\ref{eq-A(q)-inTermsOfTheta}) of the series $A(q)$, which appears in (\ref{eq-genus1-GWInv}). In this section, we make some simplifications of this expression and propose some conjectural closed formulae. We need to deal with the sum 
\begin{equation}\label{eq-A(q)-inTermsOfTheta-firstSum}
\sum_{p=0}^{n-1-r}  \Theta_{p}^{(1)}(q)\Theta_{n-1-r-p}^{(0)}(q)
\end{equation}
and
\begin{equation}\label{eq-A(q)-inTermsOfTheta-secondSum}
      \sum_{p=1}^{r}   \Theta_{n-p}^{(1)}(q)\Theta_{n-1-r+p}^{(0)}(q).
\end{equation}

\subsection{The first sum}
By the definition of $\Theta_{p}^{(0)}(q)$ and $\Theta_{p}^{(1)}(q)$, we have
\begin{eqnarray}\label{eq-Theta(0)-firstSum-expansion}
&& \sum_{p=0}^{n-1-r}   \Theta_{p}^{(1)}(q)\Theta_{n-1-r-p}^{(0)}(q)\nn\\
&=& \sum_{\beta=0}^{\infty}q^{\beta}(1+q\frac{d \mu(q)}{dq})^{n-1-r- \nu_{\mathbf{d}}\beta-1}\Phi_0(q)^2
\big(\sum_{p=0}^{n-1-r}\sum_{\begin{subarray}{c}\beta_1+\beta_2=\beta\\ \beta_1,\beta_2\geq 0\end{subarray}}
\tilde{\PZc}_{n-1-r-p,n-1-r-p- \nu_{\mathbf{d}}\beta_1}^{(\beta_1)}\tilde{\PZc}_{p,p- \nu_{\mathbf{d}}\beta_2-1}^{(\beta_2)}\big)\nn\\
&&+\sum_{\beta=0}^{\infty}q^{\beta}(1+q\frac{d \mu(q)}{dq})^{n-1-r- \nu_{\mathbf{d}}\beta}\Phi_0(q)\Phi_1(q)
\big(\sum_{p=0}^{n-1-r}\sum_{\begin{subarray}{c}\beta_1+\beta_2=\beta\\ \beta_1,\beta_2\geq 0\end{subarray}}
\tilde{\PZc}_{n-1-r-p,n-1-r-p- \nu_{\mathbf{d}}\beta_1}^{(\beta_1)}\tilde{\PZc}_{p,p- \nu_{\mathbf{d}}\beta_2}^{(\beta_2)}\big)\nn\\
&&+\sum_{\beta=0}^{\infty}q^{\beta+1}(1+q\frac{d \mu(q)}{dq})^{n-1-r- \nu_{\mathbf{d}}\beta-1}\Phi_0(q)\Phi'_0(q)\nn\\
&&\times \big(\sum_{p=0}^{n-1-r}\sum_{\begin{subarray}{c}\beta_1+\beta_2=\beta\\ \beta_1,\beta_2\geq 0\end{subarray}}
\tilde{\PZc}_{n-1-r-p,n-1-r-p- \nu_{\mathbf{d}}\beta_1}^{(\beta_1)}\tilde{\PZc}_{p,p- \nu_{\mathbf{d}}\beta_2}^{(\beta_2)}(p- \nu_{\mathbf{d}}\beta_2)\big)\nn\\
&&+\sum_{\beta=0}^{\infty}q^{\beta+1}(1+q\frac{d \mu(q)}{dq})^{n-1-r- \nu_{\mathbf{d}}\beta-2}\Phi_0(q)^2\big(\mu'(q)+q \mu''(q)\big)\nn\\
&&\times \big(\sum_{p=0}^{n-1-r}\sum_{\begin{subarray}{c}\beta_1+\beta_2=\beta\\ \beta_1,\beta_2\geq 0\end{subarray}}
\tilde{\PZc}_{n-1-r-p,n-1-r-p- \nu_{\mathbf{d}}\beta_1}^{(\beta_1)}\tilde{\PZc}_{p,p- \nu_{\mathbf{d}}\beta_2}^{(\beta_2)}\binom{p- \nu_{\mathbf{d}}\beta_2}{2}\big).
\end{eqnarray}
Let
\begin{equation*}
      U_1(n,\mathbf{d},\beta):=\sum_{p=0}^{n-1-r}\sum_{\begin{subarray}{c}\beta_1+\beta_2=\beta\\ \beta_1,\beta_2\geq 0\end{subarray}}
\tilde{\PZc}_{n-1-r-p,n-1-r-p- \nu_{\mathbf{d}}\beta_1}^{(\beta_1)}\tilde{\PZc}_{p,p- \nu_{\mathbf{d}}\beta_2-1}^{(\beta_2)},
\end{equation*}
\begin{equation*}
      U_2(n,\mathbf{d},\beta):=\sum_{p=0}^{n-1-r}\sum_{\begin{subarray}{c}\beta_1+\beta_2=\beta\\ \beta_1,\beta_2\geq 0\end{subarray}}
\tilde{\PZc}_{n-1-r-p,n-1-r-p- \nu_{\mathbf{d}}\beta_1}^{(\beta_1)}\tilde{\PZc}_{p,p- \nu_{\mathbf{d}}\beta_2}^{(\beta_2)},
\end{equation*}
\begin{equation*}
      U_3(n,\mathbf{d},\beta):=\sum_{p=0}^{n-1-r}\sum_{\begin{subarray}{c}\beta_1+\beta_2=\beta\\ \beta_1,\beta_2\geq 0\end{subarray}}
\tilde{\PZc}_{n-1-r-p,n-1-r-p- \nu_{\mathbf{d}}\beta_1}^{(\beta_1)}\tilde{\PZc}_{p,p- \nu_{\mathbf{d}}\beta_2}^{(\beta_2)}(p- \nu_{\mathbf{d}} \beta_2)(n-1-r-p- \nu_{\mathbf{d}}\beta_1).
\end{equation*}

\begin{lemma}
\begin{eqnarray}\label{eq-U2-formula}
      U_2(n,\mathbf{d},\beta)
= (n-r)\delta_{\beta,0}-(n-r- \nu_{\mathbf{d}})\mathbf{d}^{\mathbf{d}}\delta_{\beta,1}. 
\end{eqnarray}
\end{lemma}
\begin{proof}
By \cite[(2.9)]{PoZ14},
\begin{eqnarray*}
U_2(n,\mathbf{d},\beta)
&=& \sum_{\begin{subarray}{c}\beta_1+\beta_2=\beta\\ \beta_1,\beta_2\geq 0\end{subarray}}\sum_{p=0}^{n-1-r}
\tilde{\PZc}_{n-1-r-p,n-1-r-p- \nu_{\mathbf{d}}\beta_1}^{(\beta_1)}\tilde{\PZc}_{p,p- \nu_{\mathbf{d}}\beta_2}^{(\beta_2)}\\
&=& \sum_{\begin{subarray}{c}\beta_1+\beta_2=\beta\\ \beta_1,\beta_2\geq 0\end{subarray}}\sum_{p}
\tilde{\PZc}_{n-1-r-p,n-1-r-p- \nu_{\mathbf{d}}\beta_1}^{(\beta_1)}\tilde{\PZc}_{p,p- \nu_{\mathbf{d}}\beta_2}^{(\beta_2)}\\
&=& \sum_{\begin{subarray}{c}\beta_1+\beta_2=\beta\\ \beta_1,\beta_2\geq 0\end{subarray}}\sum_{p}
\tilde{\PZc}_{n-1-r-p+ \nu_{\mathbf{d}} \beta_1,n-1-r-p}\tilde{\PZc}_{p-\nu_{\mathbf{d}} \beta_1,p- \nu_{\mathbf{d}}\beta}^{(\beta_2)}\\
&=& \sum_{p= \nu_{\mathbf{d}} \beta}^{n-1-r}\sum_{\begin{subarray}{c}\beta_1+\beta_2=\beta\\ \beta_1,\beta_2\geq 0\end{subarray}}
\tilde{\PZc}_{n-1-r-p+ \nu_{\mathbf{d}} \beta_1,n-1-r-p}\tilde{\PZc}_{p-\nu_{\mathbf{d}} \beta_1,p- \nu_{\mathbf{d}}\beta}^{(\beta_2)}\\
&=&   \sum_{p= \nu_{\mathbf{d}} \beta}^{n-1-r} (\delta_{\beta,0}- \mathbf{d}^{\mathbf{d}}\delta_{\beta,1})
=(n-r)\delta_{\beta,0}-(n-r- \nu_{\mathbf{d}})\mathbf{d}^{\mathbf{d}}\delta_{\beta,1}. 
\end{eqnarray*}
\end{proof}
Using the symmetry of the range in summations we compute
\begin{eqnarray}\label{eq-U4-formula}
&&\sum_{p=0}^{n-1-r}\sum_{\begin{subarray}{c}\beta_1+\beta_2=\beta\\ \beta_1,\beta_2\geq 0\end{subarray}}
\tilde{\PZc}_{n-1-r-p,n-1-r-p- \nu_{\mathbf{d}}\beta_1}^{(\beta_1)}\tilde{\PZc}_{p,p- \nu_{\mathbf{d}}\beta_2}^{(\beta_2)}(p- \nu_{\mathbf{d}}\beta_2)\nn\\
&=& \sum_{p=0}^{n-1-r}\sum_{\begin{subarray}{c}\beta_1+\beta_2=\beta\\ \beta_1,\beta_2\geq 0\end{subarray}}
\tilde{\PZc}_{n-1-r-p,n-1-r-p- \nu_{\mathbf{d}}\beta_1}^{(\beta_1)}\tilde{\PZc}_{p,p- \nu_{\mathbf{d}}\beta_2}^{(\beta_2)}(n-1-r-p- \nu_{\mathbf{d}}\beta_1)\nn\\
&=& \frac{1}{2}\sum_{p=0}^{n-1-r}\sum_{\begin{subarray}{c}\beta_1+\beta_2=\beta\\ \beta_1,\beta_2\geq 0\end{subarray}}
\tilde{\PZc}_{n-1-r-p,n-1-r-p- \nu_{\mathbf{d}}\beta_1}^{(\beta_1)}\tilde{\PZc}_{p,p- \nu_{\mathbf{d}}\beta_2}^{(\beta_2)}\big((p- \nu_{\mathbf{d}}\beta_2)+(n-1-r-p- \nu_{\mathbf{d}}\beta_1)\big)\nn\\
&=& \frac{n-1-r- \nu_{\mathbf{d}}\beta}{2}\sum_{p=0}^{n-1-r}\sum_{\begin{subarray}{c}\beta_1+\beta_2=\beta\\ \beta_1,\beta_2\geq 0\end{subarray}}
\tilde{\PZc}_{n-1-r-p,n-1-r-p- \nu_{\mathbf{d}}\beta_1}^{(\beta_1)}\tilde{\PZc}_{p,p- \nu_{\mathbf{d}}\beta_2}^{(\beta_2)}\nn\\
&=& \frac{(n-r)(n-r-1)}{2}\delta_{\beta,0}-\frac{(n-r- \nu_{\mathbf{d}})(n-r- \nu_{\mathbf{d}}-1)}{2}\mathbf{d}^{\mathbf{d}}\delta_{\beta,1},
\end{eqnarray}
and
\begin{eqnarray}\label{eq-U5-formula}
&& \sum_{p=0}^{n-1-r}\sum_{\begin{subarray}{c}\beta_1+\beta_2=\beta\\ \beta_1,\beta_2\geq 0\end{subarray}}
\tilde{\PZc}_{n-1-r-p,n-1-r-p- \nu_{\mathbf{d}}\beta_1}^{(\beta_1)}\tilde{\PZc}_{p,p- \nu_{\mathbf{d}}\beta_2}^{(\beta_2)}\binom{p- \nu_{\mathbf{d}}\beta_2}{2}\nn\\
&=& \frac{1}{2}\sum_{p=0}^{n-1-r}\sum_{\begin{subarray}{c}\beta_1+\beta_2=\beta\\ \beta_1,\beta_2\geq 0\end{subarray}}
\tilde{\PZc}_{n-1-r-p,n-1-r-p- \nu_{\mathbf{d}}\beta_1}^{(\beta_1)}\tilde{\PZc}_{p,p- \nu_{\mathbf{d}}\beta_2}^{(\beta_2)}\big(\binom{p- \nu_{\mathbf{d}}\beta_2}{2}+\binom{n-1-r-p- \nu_{\mathbf{d}}\beta_1}{2}\big)\nn\\
&=&\frac{1}{4}\sum_{p=0}^{n-1-r}\sum_{\begin{subarray}{c}\beta_1+\beta_2=\beta\\ \beta_1,\beta_2\geq 0\end{subarray}}
\tilde{\PZc}_{n-1-r-p,n-1-r-p- \nu_{\mathbf{d}}\beta_1}^{(\beta_1)}\tilde{\PZc}_{p,p- \nu_{\mathbf{d}}\beta_2}^{(\beta_2)}\nn\\
&&\times\big(-2(p- \nu_{\mathbf{d}} \beta_2)(n-1-r-p- \nu_{\mathbf{d}}\beta_1)
+(n-1-r- \nu_{\mathbf{d}} \beta)^2-(n-1-r- \nu_{\mathbf{d}} \beta)\big)\nn\\
&=&\frac{1}{2}\binom{n-1-r- \nu_{\mathbf{d}} \beta}{2}U_2(n,\mathbf{d},\beta)
-\frac{1}{2}U_3(n,\mathbf{d},\beta)\nn\\
&=&\frac{3}{2}\binom{n-r}{3} \delta_{\beta,0}
-\frac{3\mathbf{d}^{\mathbf{d}}}{2}\binom{n-r- \nu_{\mathbf{d}}}{3} \delta_{\beta,1}
-\frac{1}{2}U_3(n,\mathbf{d},\beta).
\end{eqnarray}
By (\ref{eq-U2-formula}), (\ref{eq-U4-formula}), (\ref{eq-U5-formula}), and the results in Section \ref{sec:formulae-hypergeometricSeries-Fano}, the computation of (\ref{eq-A(q)-inTermsOfTheta-firstSum}) amounts to that of $U_1(n,\mathbf{d},\beta)$ and $U_3(n,\mathbf{d},\beta)$.
For the latter, we have the following conjecture.
\begin{conjecture}
\begin{eqnarray}
      U_3(n,\mathbf{d},\beta)=\binom{n-r}{3} \delta_{\beta,0}-\binom{|\mathbf{d}|-r}{3} \mathbf{d}^{\mathbf{d}} \delta_{\beta,1}.
\end{eqnarray}
\end{conjecture}

The sum $U_1(n,\mathbf{d},\beta)$ is the most difficult one. We have not even a conjectural formula. In the following, we give some partial results.
\begin{lemma}
Suppose  $n\geq |\mathbf{d}|+1$ and $(\beta-1)n\geq \beta|\mathbf{d}|-r-1$. Then $ U_1(n,\mathbf{d},\beta)=0$.
\end{lemma}
\begin{proof}
Then condition $(\beta-1)n\geq \beta|\mathbf{d}|-r-1$ is equivalent to $n-r-1\leq \beta \nu_{\mathbf{d}}$. Recall
\begin{equation}
      U_1(n,\mathbf{d},\beta)=\sum_{p=0}^{n-1-r}\sum_{\begin{subarray}{c}\beta_1+\beta_2=\beta\\ \beta_1,\beta_2\geq 0\end{subarray}}
\tilde{\PZc}_{n-1-r-p,n-1-r-p- \nu_{\mathbf{d}}\beta_1}^{(\beta_1)}\tilde{\PZc}_{p,p- \nu_{\mathbf{d}}\beta_2-1}^{(\beta_2)}.
\end{equation}
If $U_1(n,\mathbf{d},\beta)$ does not vanish,  both the subscripts $n-1-r-p- \nu_{\mathbf{d}}\beta_1$ and $p- \nu_{\mathbf{d}}\beta_2-1$ are nonnegative,  which imply  $n-2-r-\nu_{\mathbf{d}}\beta\geq 0$.
\end{proof}

\begin{conjecture}
For $n\geq |\mathbf{d}|$,
\begin{equation}
      U_1(n,\mathbf{d},\beta)=0\ \mbox{if and only if}\ (\beta-1)n\geq \beta|\mathbf{d}|-r-1.
\end{equation}
\end{conjecture}

\begin{lemma}
\begin{eqnarray}
      U_1(n,\mathbf{d},1)
=- \frac{\mathbf{d}^{\mathbf{d}}}{2}\big((\sum_{i=1}^{r}\frac{d_i-1}{2})^2
-\sum_{i=1}^{r}\frac{(d_i-1)(2d_i-1)}{6d_i}\big).
\end{eqnarray}
In particular,
\begin{eqnarray*}
U_1(n,d,1)=-\frac{(d-1)(d-2)(3d-1)d^{d-1}}{24}.
\end{eqnarray*}
\end{lemma}
\begin{proof}
\begin{eqnarray*}
      U_1(n,\mathbf{d},1)&=&\sum_{p=0}^{n-1-r}\sum_{\begin{subarray}{c}\beta_1+\beta_2=1\\ \beta_1,\beta_2\geq 0\end{subarray}}
\tilde{\PZc}_{n-1-r-p,n-1-r-p- \nu_{\mathbf{d}}\beta_1}^{(\beta_1)}\tilde{\PZc}_{p,p- \nu_{\mathbf{d}}\beta_2-1}^{(\beta_2)}\\
&=& \sum_{p=0}^{n-1-r}
(\tilde{\PZc}_{n-1-r-p,n-1-r-p- \nu_{\mathbf{d}}}^{(1)}\tilde{\PZc}_{p,p-1}^{(0)}
+\tilde{\PZc}_{n-1-r-p,n-1-r-p}^{(0)}\tilde{\PZc}_{p,p- \nu_{\mathbf{d}}-1}^{(1)})\\
&=& \sum_{p=0}^{n-1-r}\tilde{\PZc}_{p,p- \nu_{\mathbf{d}}-1}^{(1)},
\end{eqnarray*}
where we have used \cite[(2.7)]{PoZ14}. Again using \cite[(2.7)]{PoZ14} we have
\begin{eqnarray*}
\tilde{\PZc}_{p,p- \nu_{\mathbf{d}}-1}^{(1)}= -\prod_{i=1}^{r}d_i\cdot\mathrm{Res}_{w=0}\frac{w^{\nu_{\mathbf{d}}-p}\prod_{k=1}^{r}\prod_{j=1}^{d_k-1}(d_k w+j)}{(w+1)^{n-r-p}}.
\end{eqnarray*}
So
\begin{eqnarray*}
      U_1(n,\mathbf{d},1)
&=& -\prod_{i=1}^{r}d_i\cdot\mathrm{Res}_{w=0}\sum_{p=0}^{n-1-r}\frac{w^{\nu_{\mathbf{d}}-p}\prod_{k=1}^{r}\prod_{j=1}^{d_k-1}(d_k w+j)}{(w+1)^{n-r-p}}\\
&=&- \frac{\mathbf{d}^{\mathbf{d}}}{2}\big((\sum_{i=1}^{r}\frac{d_i-1}{2})^2
-\sum_{i=1}^{r}\frac{(d_i-1)(2d_i-1)}{6d_i}\big).
\end{eqnarray*}
\end{proof}

\begin{conjecture}
Suppose  $n\geq |\mathbf{d}|$. Then
\begin{multline*}
U_1(n,\mathbf{d},2)= \frac{(\mathbf{d}!)^2 (-1)^{n+r}}{2}\sum_{\sum j k_j\leq 2|\mathbf{d}|-n-r-2}
\Big\{\binom{2|\mathbf{d}|-2r-3-\sum j k_j}{2|\mathbf{d}|-n-r-2-\sum j k_j}\\
\times \prod_{j=1}^{\infty}\frac{1}{k_j!}(-\frac{2 \sum_{i=1}^r d_i^j h_j(d_i)}{j})^{k_j}\Big\},
\end{multline*}
where the summation is taken over the set of lists $K=(k_1,k_2,\dots)$ of nonzero integers $k_j$ satisfying 
\[
\sum_{j=1}^{\infty}j k_j\leq 2|\mathbf{d}|-n-r-2.
\]
\end{conjecture}
For example, when $r=2$, for $d_1+d_2\geq 6$, this conjecture reads 
\begin{eqnarray*}
&&U_1(2|\mathbf{d}|-6,\mathbf{d},2)\\
&=&(d_1!)^2(d_2!)^2\Big( (|\mathbf{d}|-\frac{7}{2})(|\mathbf{d}|-4)+(8-2|\mathbf{d}|)\big(d_1 h_1(d_1)+d_2h_1(d_2)\big)\\
&&+\big(d_1 h_1(d_1)+d_2h_1(d_2)\big)^2-\frac{1}{2}\big(d_1^2h_2(d_1)+d_2^2h_2(d_2)\big)\Big).
\end{eqnarray*}

\begin{remark}
For given $r$ and $2d-n$, one can in principle show this conjecture by computations of residues. We have a proof  in the following cases (i) $2|\mathbf{d}|-n=r+2$; (ii) $r=1$ and $2d-n=5$. 
\end{remark}

\subsection{The second sum}

\begin{eqnarray}\label{eq-Theta(0)-secondSum-expansion}
&& \sum_{p=1}^{r}   \Theta_{n-p}^{(1)}(q)\Theta_{n-1-r+p}^{(0)}(q)\nn\\
&=& \sum_{p=1}^{r}  \big(\sum_{\beta_1=0}^{\infty}\tilde{\PZc}_{n-1-r+p,n-1-r+p- \nu_{\mathbf{d}}\beta_1}^{(\beta_1)}q^{\beta_1}
(1+q\frac{d \mu(q)}{dq})^{n-1-r+p- \nu_{\mathbf{d}}\beta_1}\Phi_0(q)\big)\nn\\
&& \times\Big(\sum_{\beta_2=0}^{\infty}\tilde{\PZc}_{n-p,n-p- \nu_{\mathbf{d}}\beta_2-1}^{(\beta_2)}q^{\beta_2}
(1+q\frac{d \mu(q)}{dq})^{n-p- \nu_{\mathbf{d}}\beta_2-1}\Phi_0(q)\nn\\
&&+ \sum_{\beta_2=0}^{\infty}\tilde{\PZc}_{n-p,n-p- \nu_{\mathbf{d}}\beta_2}^{(\beta_2)}q^{\beta_2}
(1+q\frac{d \mu(q)}{dq})^{n-p- \nu_{\mathbf{d}}\beta_2}\Phi_1(q)\nn\\
&&+\sum_{\beta_2=0}^{\infty}\tilde{\PZc}_{n-p,n-p- \nu_{\mathbf{d}}\beta_2}^{(\beta_2)}q^{\beta_2+1}
(n-p- \nu_{\mathbf{d}}\beta_2)
(1+q\frac{d \mu(q)}{dq})^{n-p- \nu_{\mathbf{d}}\beta_2-1}\frac{d\Phi_0(q) }{dq}\nn\\
&&+\sum_{\beta_2=0}^{\infty}\tilde{\PZc}_{n-p,n-p- \nu_{\mathbf{d}}\beta_2}^{(\beta_2)}q^{\beta_2+1}
\binom{n-p- \nu_{\mathbf{d}}\beta_2}{2}
(1+q\frac{d \mu(q)}{dq})^{n-p- \nu_{\mathbf{d}}\beta_2-2}\big(\frac{d \mu(q)}{dq}+q\frac{d^2 \mu(q)}{dq^2}\big)\Phi_0(q)\Big)\nn\\
&=& \sum_{\beta=0}^{\infty}q^{\beta}(1+q\frac{d \mu(q)}{dq})^{2n-1-r- \nu_{\mathbf{d}}\beta-1}\Phi_0(q)^2
\big(\sum_{p=1}^{r}\sum_{\begin{subarray}{c}\beta_1+\beta_2=\beta\\ \beta_1,\beta_2\geq 0\end{subarray}}
\tilde{\PZc}_{n-1-r+p,n-1-r+p- \nu_{\mathbf{d}}\beta_1}^{(\beta_1)}\tilde{\PZc}_{n-p,n-p- \nu_{\mathbf{d}}\beta_2-1}^{(\beta_2)}\big)\nn\\
&&+\sum_{\beta=0}^{\infty}q^{\beta}(1+q\frac{d \mu(q)}{dq})^{2n-1-r- \nu_{\mathbf{d}}\beta}\Phi_0(q)\Phi_1(q)
\big(\sum_{p=1}^{r}\sum_{\begin{subarray}{c}\beta_1+\beta_2=\beta\\ \beta_1,\beta_2\geq 0\end{subarray}}
\tilde{\PZc}_{n-1-r+p,n-1-r+p- \nu_{\mathbf{d}}\beta_1}^{(\beta_1)}\tilde{\PZc}_{n-p,n-p- \nu_{\mathbf{d}}\beta_2}^{(\beta_2)}\big)\nn\\
&&+\sum_{\beta=0}^{\infty}q^{\beta+1}(1+q\frac{d \mu(q)}{dq})^{2n-1-r- \nu_{\mathbf{d}}\beta-1}\Phi_0(q)\Phi'_0(q)\nn\\
&&\times \big(\sum_{p=1}^{r}\sum_{\begin{subarray}{c}\beta_1+\beta_2=\beta\\ \beta_1,\beta_2\geq 0\end{subarray}}
\tilde{\PZc}_{n-1-r+p,n-1-r+p- \nu_{\mathbf{d}}\beta_1}^{(\beta_1)}\tilde{\PZc}_{n-p,n-p- \nu_{\mathbf{d}}\beta_2}^{(\beta_2)}(n-p- \nu_{\mathbf{d}}\beta_2)\big)\nn\\
&&+\sum_{\beta=0}^{\infty}q^{\beta+1}(1+q\frac{d \mu(q)}{dq})^{2n-1-r- \nu_{\mathbf{d}}\beta-2}\Phi_0(q)^2\big(\mu'(q)+q \mu''(q)\big)\nn\\
&&\times \big(\sum_{p=1}^{r}\sum_{\begin{subarray}{c}\beta_1+\beta_2=\beta\\ \beta_1,\beta_2\geq 0\end{subarray}}
\tilde{\PZc}_{n-1-r+p,n-1-r+p- \nu_{\mathbf{d}}\beta_1}^{(\beta_1)}\tilde{\PZc}_{n-p,n-p- \nu_{\mathbf{d}}\beta_2}^{(\beta_2)}\binom{n-p- \nu_{\mathbf{d}}\beta_2}{2}\big).
\end{eqnarray}
Let
\begin{equation*}
      V_1(n,\mathbf{d},\beta)=\sum_{p=1}^{r}\sum_{\begin{subarray}{c}\beta_1+\beta_2=\beta\\ \beta_1,\beta_2\geq 0\end{subarray}}
\tilde{\PZc}_{n-1-r+p,n-1-r+p- \nu_{\mathbf{d}}\beta_1}^{(\beta_1)}\tilde{\PZc}_{n-p,n-p- \nu_{\mathbf{d}}\beta_2-1}^{(\beta_2)},
\end{equation*}
\begin{equation*}
      V_2(n,\mathbf{d},\beta)=\sum_{p=1}^{r}\sum_{\begin{subarray}{c}\beta_1+\beta_2=\beta\\ \beta_1,\beta_2\geq 0\end{subarray}}
\tilde{\PZc}_{n-1-r+p,n-1-r+p- \nu_{\mathbf{d}}\beta_1}^{(\beta_1)}\tilde{\PZc}_{n-p,n-p- \nu_{\mathbf{d}}\beta_2}^{(\beta_2)},
\end{equation*}
\begin{equation*}
      V_3(n,\mathbf{d},\beta)=\sum_{p=1}^{r}\sum_{\begin{subarray}{c}\beta_1+\beta_2=\beta\\ \beta_1,\beta_2\geq 0\end{subarray}}
\tilde{\PZc}_{n-1-r+p,n-1-r+p- \nu_{\mathbf{d}}\beta_1}^{(\beta_1)}\tilde{\PZc}_{n-p,n-p- \nu_{\mathbf{d}}\beta_2}^{(\beta_2)}(n-1-r+p- \nu_{\mathbf{d}}\beta_1)(n-p- \nu_{\mathbf{d}}\beta_2).
\end{equation*}
Using the symmetry of the range in summations we compute
\begin{eqnarray*}
&&\sum_{p=1}^{r}\sum_{\begin{subarray}{c}\beta_1+\beta_2=\beta\\ \beta_1,\beta_2\geq 0\end{subarray}}
\tilde{\PZc}_{n-1-r+p,n-1-r+p- \nu_{\mathbf{d}}\beta_1}^{(\beta_1)}\tilde{\PZc}_{n-p,n-p- \nu_{\mathbf{d}}\beta_2}^{(\beta_2)}(n-p- \nu_{\mathbf{d}}\beta_2)\\
  &=& \frac{1}{2}\sum_{p=1}^{r}\sum_{\begin{subarray}{c}\beta_1+\beta_2=\beta\\ \beta_1,\beta_2\geq 0\end{subarray}}
\tilde{\PZc}_{n-1-r+p,n-1-r+p- \nu_{\mathbf{d}}\beta_1}^{(\beta_1)}\tilde{\PZc}_{n-p,n-p- \nu_{\mathbf{d}}\beta_2}^{(\beta_2)}\big((n-p- \nu_{\mathbf{d}}\beta_2)
+(n-1-r+p- \nu_{\mathbf{d}}\beta_1)\big)\\
&=& \frac{2n-1-r- \nu_{\mathbf{d}}  \beta}{2}\sum_{p=1}^{r}\sum_{\begin{subarray}{c}\beta_1+\beta_2=\beta\\ \beta_1,\beta_2\geq 0\end{subarray}}
\tilde{\PZc}_{n-1-r+p,n-1-r+p- \nu_{\mathbf{d}}\beta_1}^{(\beta_1)}\tilde{\PZc}_{n-p,n-p- \nu_{\mathbf{d}}\beta_2}^{(\beta_2)}\\
&=& \frac{2n-1-r- \nu_{\mathbf{d}}  \beta}{2}V_2(n,\mathbf{d},\beta),
\end{eqnarray*}
and
\begin{eqnarray*}
&&\sum_{p=1}^{r}\sum_{\begin{subarray}{c}\beta_1+\beta_2=\beta\\ \beta_1,\beta_2\geq 0\end{subarray}}
\tilde{\PZc}_{n-1-r+p,n-1-r+p- \nu_{\mathbf{d}}\beta_1}^{(\beta_1)}\tilde{\PZc}_{n-p,n-p- \nu_{\mathbf{d}}\beta_2}^{(\beta_2)}\binom{n-p- \nu_{\mathbf{d}}\beta_2}{2}\\
  &=& \frac{1}{2} \sum_{p=1}^{r}\sum_{\begin{subarray}{c}\beta_1+\beta_2=\beta\\ \beta_1,\beta_2\geq 0\end{subarray}}
\tilde{\PZc}_{n-1-r+p,n-1-r+p- \nu_{\mathbf{d}}\beta_1}^{(\beta_1)}\tilde{\PZc}_{n-p,n-p- \nu_{\mathbf{d}}\beta_2}^{(\beta_2)}\big(\binom{n-p- \nu_{\mathbf{d}}\beta_2}{2}+\binom{n-1-r+p- \nu_{\mathbf{d}}\beta_1}{2}\big)\\
&=&\frac{1}{2}\binom{2n-1-r- \nu_{\mathbf{d}}  \beta}{2}\sum_{p=1}^{r}\sum_{\begin{subarray}{c}\beta_1+\beta_2=\beta\\ \beta_1,\beta_2\geq 0\end{subarray}}
\tilde{\PZc}_{n-1-r+p,n-1-r+p- \nu_{\mathbf{d}}\beta_1}^{(\beta_1)}\tilde{\PZc}_{n-p,n-p- \nu_{\mathbf{d}}\beta_2}^{(\beta_2)}\\
&&-\frac{1}{2} \sum_{p=1}^{r}\sum_{\begin{subarray}{c}\beta_1+\beta_2=\beta\\ \beta_1,\beta_2\geq 0\end{subarray}}
\tilde{\PZc}_{n-1-r+p,n-1-r+p- \nu_{\mathbf{d}}\beta_1}^{(\beta_1)}\tilde{\PZc}_{n-p,n-p- \nu_{\mathbf{d}}\beta_2}^{(\beta_2)}(n-p- \nu_{\mathbf{d}}\beta_2)(n-1-r+p- \nu_{\mathbf{d}}\beta_1)\\
&=&\frac{1}{2}\binom{2n-1-r- \nu_{\mathbf{d}}  \beta}{2}V_2(n,\mathbf{d},\beta)-\frac{1}{2}V_3(n,\mathbf{d},\beta).
\end{eqnarray*} 
Hence the computation of (\ref{eq-A(q)-inTermsOfTheta-secondSum}) amounts to that of $V_1(n,\mathbf{d},\beta)$, $V_2(n,\mathbf{d},\beta)$ and $V_3(n,\mathbf{d},\beta)$. We have the following conjecture.
\begin{conjecture}
   \begin{equation}
      V_1(n,\mathbf{d},\beta)=-\frac{r(|\mathbf{d}|-1)\mathbf{d}^{\mathbf{d}}}{2} \delta_{\beta,1}
      + \frac{r(|\mathbf{d}|-1)\mathbf{d}^{2\mathbf{d}}}{2}\delta_{\beta,2}.
\end{equation}
\begin{equation}
      V_2(n,\mathbf{d},\beta)=r\delta_{\beta,0}-2r \mathbf{d}^{\mathbf{d}} \delta_{\beta,1}+ r \mathbf{d}^{2 \mathbf{d}}\delta_{\beta,2}.
\end{equation}
\begin{eqnarray}
V_3(n,\mathbf{d},\beta)&=&\big(n^2 r-n (r^2+r)+\frac{r(r+1)(r+2)}{6}\big) \delta_{\beta,0}\nn\\
      &&-\Big( (2r |\mathbf{d}|-r^2-r)n-\big((r^2+r)|\mathbf{d}|-\frac{r(r+1)(r+2)}{3}\big)\Big)\mathbf{d}^{\mathbf{d}} \delta_{\beta,1}\nn\\
      &&      +\big(n^2 r-n (r^2+r)+\frac{r(r+1)(r+2)}{6}\big) \mathbf{d}^{2 \mathbf{d}}\delta_{\beta,2}.
\end{eqnarray}
\end{conjecture}

\begin{appendices}

\section{Proof of Lemma \ref{lem-symmetricReduction-0<=a,b,c<=n,n+1<=d<=n+m}}\label{sec:proof-lem-symmetricReduction-0<=a,b,c<=n,n+1<=d<=n+m}
From (\ref{eq-wdvv230}) we get
\begin{equation}\label{eq-wdvv1-der-s}
      F_{abes}\tauG^{ef}F_{sf}+F_{abe}\tauG^{ef}F_{ssf}+2F_{sab}F_{ss}+2sF_{ssab}F_{ss}+2sF_{sab}F_{sss}
      =F_{ssa}F_{sb}+F_{sa}F_{ssb},
\end{equation}
\begin{equation}\label{eq-wdvv1-der-c}
      F_{abce}\tauG^{ei}F_{si}+F_{abe}\tauG^{ei}F_{sci}+2sF_{sabc}F_{ss}+2sF_{sab}F_{ssc}
      =F_{sac}F_{sb}+F_{sa}F_{sbc},
\end{equation}
\begin{eqnarray}\label{eq-wdvv1-der-cf}
&& F_{abcfe}\tauG^{ei}F_{si}+F_{abce}\tauG^{ei}F_{sif}+F_{abef}\tauG^{ei}F_{sci}+F_{abe}\tauG^{ei}F_{scif}\nn\\
&&+2sF_{sabcf}F_{ss}+2sF_{sabc}F_{ssf}+2sF_{sabf}F_{ssc}+2sF_{sab}F_{sscf}\nn\\
&=& F_{sacf}F_{sb}+F_{sac}F_{sbf}+F_{saf}F_{sbc}+F_{sa}F_{sbcf},
\end{eqnarray}
and
\begin{eqnarray}\label{eq-wdvv1-der-sc}
&& F_{sabce}\tauG^{ei}F_{si}+F_{abce}\tauG^{ei}F_{ssi}+F_{sabe}\tauG^{ei}F_{sci}+F_{abe}\tauG^{ei}F_{ssci}\nn\\
&&+2F_{sabc}F_{ss}+2sF_{ssabc}F_{ss}+2sF_{sabc}F_{sss}
+2F_{sab}F_{ssc}+2sF_{ssab}F_{ssc}+2sF_{sab}F_{sssc}\nn\\
&=& F_{ssac}F_{sb}+F_{sac}F_{ssb}+F_{ssa}F_{sbc}+F_{sa}F_{ssbc}.
\end{eqnarray}
From (\ref{eq-wdvv240}) we get
\begin{equation}\label{eq-wdvv2-der-ab}
      F_{sae}\tauG^{ef}F_{sfb}+F_{se}\tauG^{ef}F_{sfab}+2sF_{ssa}F_{ssb}+2sF_{ss}F_{ssab}=0,
\end{equation}
and
\begin{eqnarray}\label{eq-wdvv2-der-abc}
&& F_{sace}\tauG^{ef}F_{sfb}+F_{sae}\tauG^{ef}F_{sfbc}+F_{sce}\tauG^{ef}F_{sfab}+F_{se}\tauG^{ef}F_{sfabc}\nn\\
&&+2sF_{ssac}F_{ssb}+2sF_{ssa}F_{ssbc}+2sF_{ssc}F_{ssab}+2sF_{ss}F_{ssabc}=0.
\end{eqnarray}

\begin{proof}[Proof of Lemma \ref{lem-symmetricReduction-0<=a,b,c<=n,n+1<=d<=n+m}]
Suppose that the dimension $n$ is even. Then
\begin{multline*}
\sum_{\mu=0}^{n+m} \sum_{\nu=0}^{n+m}\sum_{\lambda=0}^{n+m} \sum_{\rho=0}^{n+m}F_{ab\lambda}\tauG^{\lambda\mu}G_{\mu\nu}\tauG^{\nu\rho}F_{\rho cd}
=\tau^d\big(F_{abe}\tauG^{ei}F_{scf}\tauG^{fj}G_{ij}+2sF_{sab}F_{scf}\tauG^{fj}G_{sj}\\+F_{abe}\tauG^{ei}F_{sc}G_{si}
+2sF_{abe}\tauG^{ei}F_{ssc}G_{si}
+F_{sab}(F_{sc}G_{s}+2s F_{sc}G_{ss}+2s F_{ssc}G_s+4s^2 F_{ssc}G_{ss})\big),
\end{multline*}

\begin{multline*}
\sum_{\mu=0}^{n+m} \sum_{\nu=0}^{n+m}\sum_{\lambda=0}^{n+m} \sum_{\rho=0}^{n+m}F_{ab\lambda}\tauG^{\lambda\mu}F_{\mu c \rho}\tauG^{\rho\nu}G_{\nu d}
=\tau^d\big(F_{abe}\tauG^{ei}F_{cif}\tauG^{fj}G_{sj}+2sF_{sab}F_{scf}\tauG^{fj}G_{sj}\\
+ F_{abe}\tauG^{ei} F_{sci}G_s
+2sF_{abe}\tauG^{ei} F_{sci}G_{ss}
+F_{sab}(F_{sc}G_s+2sF_{sc}G_{ss}+2sF_{ssc}G_s+4s^2 F_{ssc}G_{ss})\big),
\end{multline*}

\begin{multline*}
\sum_{\mu=0}^{n+m} \sum_{\nu=0}^{n+m}\sum_{\lambda=0}^{n+m} \sum_{\rho=0}^{n+m}F_{ab\lambda}\tauG^{\lambda\mu}F_{\mu d \rho}\tauG^{\rho\nu}G_{\nu c}
=\tau^d\big( F_{abe}\tauG^{ei}  F_{sif}\tauG^{fj}G_{cj}
+F_{sab}F_{sf}\tauG^{fj}G_{cj}\\
+2s F_{sab}F_{ssf}\tauG^{fj}G_{cj}
+F_{abe}\tauG^{ei}F_{si}G_{sc}+2sF_{abe}\tauG^{ei}F_{ssi}G_{sc}
+F_{sab}(6sF_{ss}G_{sc}+4s^2 F_{sss}G_{sc})\big),
\end{multline*}

\begin{multline*}
\sum_{\mu=0}^{n+m} \sum_{\nu=0}^{n+m}\sum_{\lambda=0}^{n+m} \sum_{\rho=0}^{n+m}F_{ad\lambda}\tauG^{\lambda\mu}F_{\mu b \rho}\tauG^{\rho\nu}G_{\nu c}
= \tau^d\big( F_{sae}\tauG^{ei} F_{bif}\tauG^{fj}G_{cj}
+F_{sa}F_{sbf}\tauG^{fj}G_{cj}\\
+2sF_{ssa}F_{sbf}\tauG^{fj}G_{cj}
+2s F_{sae}\tauG^{ei} F_{sbi}G_{sc}+ (F_{sa}F_{sb}+2sF_{sa}F_{ssb}+2sF_{ssa}F_{sb}+4s^2 F_{ssa}F_{ssb})G_{sc}\big),
\end{multline*}

\begin{multline*}
\sum_{\mu=0}^{n+m} \sum_{\nu=0}^{n+m}\sum_{\lambda=0}^{n+m} \sum_{\rho=0}^{n+m}F_{ab\lambda}\tauG^{\lambda\mu}F_{\mu cd\rho}\tauG^{\rho\nu}G_{\nu}
= \tau^d( F_{abe}\tauG^{ei} F_{scif}\tauG^{fj}G_j+F_{sab}F_{scf}\tauG^{fj}G_j\\
+2sF_{sab}F_{sscf}\tauG^{fj}G_j
+F_{abe}\tauG^{ei}F_{sci}G_s+2sF_{abe}\tauG^{ei}F_{ssci}G_s+6sF_{sab}F_{ssc}G_s+4s^2 F_{sab}F_{sssc}G_s),
\end{multline*}

\begin{multline*}
\sum_{\mu=0}^{n+m} \sum_{\nu=0}^{n+m}\sum_{\lambda=0}^{n+m} \sum_{\rho=0}^{n+m}F_{ad\lambda}\tauG^{\lambda\mu}F_{\mu bc\rho}\tauG^{\rho\nu}G_{\nu}
= \tau^d\big( F_{sae}\tauG^{ei}F_{bcif}\tauG^{fj}G_{j} + F_{sa}F_{sbcf}\tauG^{fj}G_j\\
+2sF_{ssa}F_{sbcf}\tauG^{fj}G_j
+(2s F_{sae}\tauG^{ei} F_{sbci}+ F_{sa}F_{sbc}+2s F_{sa}F_{ssbc}+2s F_{ssa}F_{sbc}+4s^2 F_{ssa}F_{ssbc})G_s\big),
\end{multline*}

\begin{multline*}   
\sum_{\mu=0}^{n+m} \sum_{\nu=0}^{n+m}\sum_{\lambda=0}^{n+m} \sum_{\rho=0}^{n+m}F_{abc\lambda}\tauG^{\lambda\mu}F_{\mu d\rho}\tauG^{\rho\nu}G_{\nu}
= \tau^d( F_{abce}\tauG^{ei} F_{sif}\tauG^{fj}G_j
+F_{sabc}F_{sf}\tauG^{fj}G_j\\
+2s F_{sabc}F_{ssf}\tauG^{fj}G_j
+F_{abce}\tauG^{ei}F_{si}G_s+2sF_{abce}\tauG^{ei}F_{ssi}G_s
+6sF_{sabc}F_{ss}G_s+4s^2 F_{sabc}F_{sss}G_s),
\end{multline*}

\begin{multline*}   
\sum_{\mu=0}^{n+m} \sum_{\nu=0}^{n+m}\sum_{\lambda=0}^{n+m} \sum_{\rho=0}^{n+m}F_{abd\lambda}\tauG^{\lambda\mu}F_{\mu c\rho}\tauG^{\rho\nu}G_{\nu}
=\tau^d\big( F_{sabe}\tauG^{ei}F_{cif}\tauG^{fj}G_j
+F_{sab}F_{scf}\tauG^{fj}G_j\\
+2s F_{ssab}F_{scf}\tauG^{fj}G_j
+2sF_{sabe}\tauG^{ei}F_{sci}G_s+(F_{sab}F_{sc}+2sF_{sab}F_{ssc}+2sF_{ssab}F_{sc}+4s^2 F_{ssab}F_{ssc})G_s\big).
\end{multline*}
Summing over the permutations, the sum of the first two rows of (\ref{Grelation-F&G}) is equal to
\begin{eqnarray}\label{eq-proof-symmetricReduction-0<=a,b,c<=n,n+1<=d<=n+m-genus1Part}
&&  \tau^d\Big(12\sum_{\mathscr{P}(a,b,c)}F_{abe}\tauG^{ei}F_{scf}\tauG^{fj}G_{ij}\nn\\
&&+\sum_{\mathscr{P}(a,b,c)}(12F_{abf}F_{sc}-4F_{abe}\tauG^{ei}F_{cif}+16sF_{sab}F_{scf}
+24sF_{abf}F_{ssc})\tauG^{fj}G_{sj}\nn\\
&&-4\sum_{\mathscr{P}(a,b,c)}(F_{bce}\tauG^{ei} F_{sif}+F_{sbc}F_{sf}+2s F_{sbc}F_{ssf}+F_{sce}\tauG^{ei} F_{bif}+F_{sc}F_{sbf}\nn\\
&&+2sF_{ssc}F_{sbf}+F_{sbe}\tauG^{ei} F_{cif}+F_{sb}F_{scf}+2sF_{ssb}F_{scf})\tauG^{fj}G_{aj}\nn\\
&&-4\sum_{\mathscr{P}(a,b,c)}(F_{bce}\tauG^{ei}F_{si}+2sF_{bce}\tauG^{ei}F_{ssi}+6sF_{sbc}F_{ss}+4s^2 F_{sbc}F_{sss}\nn\\
&&+2F_{sb}F_{sc}+4sF_{sb}F_{ssc}+4sF_{ssb}F_{sc}+8s^2 F_{ssb}F_{ssc}+4s F_{sbe}\tauG^{ef} F_{scf})G_{sa}\nn\\
&&+\big(\sum_{\mathscr{P}(a,b,c)}(-2 F_{abe}\tauG^{ei} F_{scif}+4F_{sab}F_{scf}-4sF_{sab}F_{sscf}
-2 F_{sce}\tauG^{ei}F_{abif} -2 F_{sc}F_{sabf}-4sF_{ssc}F_{sabf}\nn\\
&&+6 F_{sabe}\tauG^{ei}F_{cif}+12s F_{ssab}F_{scf})+12F_{abce}\tauG^{ei} F_{sif}+12F_{sabc}F_{sf}+24s F_{sabc}F_{ssf}\big)\tauG^{fj}G_j\nn\\
&&+\big(\sum_{\mathscr{P}(a,b,c)}(12F_{sab}F_{sc}+12sF_{sab}F_{ssc}-6F_{abe}\tauG^{ei} F_{sci}+8sF_{ssbc}F_{sa}\nn\\
&&+16s^2 F_{ssbc}F_{ssa}-4sF_{abe}\tauG^{ei}F_{ssci}-8s^2 F_{sab}F_{sssc}+8sF_{sabe}\tauG^{ei}F_{sci})\nn\\
&&+12F_{abce}\tauG^{ei}F_{si}+24sF_{abce}\tauG^{ei}F_{ssi}+72sF_{sabc}F_{ss}+48s^2 F_{sabc}F_{sss}\big)G_s\nn\\
&&+8s\sum_{\mathscr{P}(a,b,c)}(2F_{sab}F_{sc}+4sF_{sab} F_{ssc}-F_{abe}\tauG^{ei} F_{sci})G_{ss}\Big).
\end{eqnarray}

By (\ref{eq-wdvv1-der-c}) we get
\begin{eqnarray}\label{eq-proof-symmetricReduction-0<=a,b,c<=n,n+1<=d<=n+m-genus1Part-simplification-1}
&& F_{bce}\tauG^{ei} F_{sif}+F_{sbc}F_{sf}+2s F_{sbc}F_{ssf}\nn\\
&&+F_{sce}\tauG^{ei} F_{bif}+F_{sc}F_{sbf}+2sF_{ssc}F_{sbf}
+F_{sbe}\tauG^{ei} F_{cif}+F_{sb}F_{scf}+2sF_{ssb}F_{scf}\nn\\
&=&3F_{sbc}F_{sf}+3F_{sc}F_{sbf}+3F_{sb}F_{scf}
-3F_{bcfe}\tauG^{ei}F_{si}-6sF_{sbcf}F_{ss}.
\end{eqnarray}
Using (\ref{eq-wdvv230}), (\ref{eq-wdvv1-der-s}), and (\ref{eq-wdvv2-der-ab}) successively we get
\begin{eqnarray}\label{eq-proof-symmetricReduction-0<=a,b,c<=n,n+1<=d<=n+m-genus1Part-simplification-2}
&& F_{bce}\tauG^{ei}F_{si}+2sF_{bce}\tauG^{ei}F_{ssi}+6sF_{sbc}F_{ss}+4s^2 F_{sbc}F_{sss}\nn\\
&&+2F_{sb}F_{sc}+4sF_{sb}F_{ssc}+4sF_{ssb}F_{sc}+8s^2 F_{ssb}F_{ssc}+4s F_{sbe}\tauG^{ef} F_{scf}\nn\\
&=& 3F_{sb}F_{sc}+6sF_{sb}F_{ssc}+6sF_{ssb}F_{sc}+12s^2 F_{ssb}F_{ssc}+6s F_{sbe}\tauG^{ef} F_{scf}.
\end{eqnarray}

Using (\ref{eq-wdvv1-der-cf})  we get
\begin{eqnarray}\label{eq-proof-symmetricReduction-0<=a,b,c<=n,n+1<=d<=n+m-genus1Part-simplification-3}
&&\sum_{\mathscr{P}(a,b,c)}(-2 F_{abe}\tauG^{ei} F_{scif}+4F_{sab}F_{scf}-4sF_{sab}F_{sscf}
-2 F_{sce}\tauG^{ei}F_{abif} -2 F_{sc}F_{sabf}-4sF_{ssc}F_{sabf}\nn\\
&&+6 F_{sabe}\tauG^{ei}F_{cif}+12s F_{ssab}F_{scf})+12F_{abce}\tauG^{ei} F_{sif}+12F_{sabc}F_{sf}+24s F_{sabc}F_{ssf}\nn\\
&=& \sum_{\mathscr{P}(a,b,c)}(12F_{sbc}F_{saf}-12F_{bcfe}\tauG^{ei}F_{sia}-24sF_{sbcf}F_{ssa})\nn\\
&&-24 F_{abcfe}\tauG^{ei}F_{si}+24F_{abce}\tauG^{ei} F_{sif}+48F_{sabc}F_{sf}
-48s F_{sabcf}F_{ss}+48s F_{sabc}F_{ssf}.
\end{eqnarray}
Using (\ref{eq-wdvv1-der-c}), (\ref{eq-wdvv1-der-sc}) and (\ref{eq-wdvv2-der-abc}) successively we get
\begin{eqnarray}\label{eq-proof-symmetricReduction-0<=a,b,c<=n,n+1<=d<=n+m-genus1Part-simplification-4}
&&\sum_{\mathscr{P}(a,b,c)}(12F_{sab}F_{sc}+12sF_{sab}F_{ssc}-6F_{abe}\tauG^{ei} F_{sci}+8sF_{ssbc}F_{sa}\nn\\
&&+16s^2 F_{ssbc}F_{ssa}-4sF_{abe}\tauG^{ei}F_{ssci}-8s^2 F_{sab}F_{sssc}+8sF_{sabe}\tauG^{ei}F_{sci})\nn\\
&&+12F_{abce}\tauG^{ei}F_{si}+24sF_{abce}\tauG^{ei}F_{ssi}+72sF_{sabc}F_{ss}+48s^2 F_{sabc}F_{sss}\nn\\
&=& 24s\sum_{\mathscr{P}(a,b,c)}F_{sbc}F_{ssa}+ 48(F_{abce}\tauG^{ei}F_{si}
+sF_{abce}\tauG^{ei}F_{ssi}+3sF_{sabc}F_{ss}+2s^2 F_{sabc}F_{sss}).
\end{eqnarray}
Finally using (\ref{eq-wdvv1-der-c}) we get
\begin{eqnarray}\label{eq-proof-symmetricReduction-0<=a,b,c<=n,n+1<=d<=n+m-genus1Part-simplification-5}
&&\sum_{\mathscr{P}(a,b,c)}(2F_{sab}F_{sc}+4sF_{sab} F_{ssc}-F_{abe}\tauG^{ei} F_{sci})\nn\\
&=& 6s\sum_{\mathscr{P}(a,b,c)}F_{sbc}F_{ssa}+6(F_{abce}\tauG^{ei}F_{si}+2sF_{sabc}F_{ss}).
\end{eqnarray}
Applying (\ref{eq-proof-symmetricReduction-0<=a,b,c<=n,n+1<=d<=n+m-genus1Part-simplification-1})-(\ref{eq-proof-symmetricReduction-0<=a,b,c<=n,n+1<=d<=n+m-genus1Part-simplification-5}) to (\ref{eq-proof-symmetricReduction-0<=a,b,c<=n,n+1<=d<=n+m-genus1Part}) yields that  the sum of the first two rows of (\ref{Grelation-F&G}) is equal to $\tau^d\times \mathrm{LHS}$ of (\ref{eq-symmetricReduction-0<=a,b,c<=n,n+1<=d<=n+m}).

A direct computation of  the sum of the third row of (\ref{Grelation-F&G}) gives RHS of (\ref{eq-symmetricReduction-0<=a,b,c<=n,n+1<=d<=n+m}). We omit the details because  RHS is not used in the application  in this paper of Lemma \ref{lem-symmetricReduction-0<=a,b,c<=n,n+1<=d<=n+m}.

If the dimension $n$ of $X$ is odd, we take $n+1\leq d\leq n+\frac{m}{2}$. We adopt the conventions introduced in the odd dimensional cases of the proof of Lemma \ref{lem-symmetridReduction-0<=a,b,c,d<=n}.
 Then all the equations in the above proof are still true, with the factor $\tau^d$ replaced by $-\tau^{d+\frac{m}{2}}$. Since $\tau^d s^{\frac{m}{2}}=0$, (\ref{eq-symmetricReduction-0<=a,b,c<=n,n+1<=d<=n+m}) holds in the $\mod s^{\frac{m}{2}}$ sense.
\end{proof}

\section{Proof of Lemma \ref{lem-CubicHypersurfaces-(2,2)-sumZ}}\label{sec:proof-eq-CubicHypersurfaces-(2,2)-sumZ}
In this section, we prove Lemma \ref{lem-CubicHypersurfaces-(2,2)-sumZ}.
Let $X$ be a smooth complete intersection in $\mathbb{P}^{n+r}$  of dimension $n\geq 3$ and multidegree $\mathbf{d}=(d_1,\dots,d_r)$.
We recall the results in \cite[\S D.8, D.9]{Hu15}. 
For $I=(i_1,\dots,i_r)$ and $J=(j_1,\dots,j_r)\in \mathbb{Z}^{r}$, we say $J\leq I$ if $j_k\leq i_k$ for $1\leq k\leq r$. Define
\begin{equation}
      I-J=(i_1-j_1,\dots,i_r-j_r),
\end{equation}
and
\begin{equation}\label{eq-binomOfLists}
\binom{I}{J}=\prod_{k=1}^r \binom{i_k}{j_k}.
\end{equation}
We denote $(0,\dots,0)\in \mathbb{Z}^r$ by 0, when no confusion arises in the context.
For $I=(i_0,i_1,\dots,i_{n})\in \mathbb{Z}_{\geq 0}^{n+1}$, we define
\begin{equation}\label{eq-def-tauI}
\partial_{\tau^I}:=(\partial_{\tau^0})^{i_0}\circ\dots\circ (\partial_{\tau^{n}})^{i_{n}}.
\end{equation}
For $0\leq j\leq n$ let $\mathbf{e}_j$ be the $j$-th unit vector in $\mathbb{Z}^{n+1}$.
Then for $I=(p_0,\dots,p_n)\in \mathbb{Z}_{\geq 0}^{n+1}$ with $|I|\geq 1$,
\begin{eqnarray}\label{eq-recursion-F(0)-EulerField-I}
\partial_{\tau^{I+\mathbf{e}_1}}F^{(0)}(0)
&=&\frac{1}{\mathsf{a}_X}\big(-\sum_{i=0}^{n}\sum_{j=0}^n\sum_{k=0}^n
(1-i)W_j^i  M_{i}^k p_j\partial_{\tau^{I+\mathbf{e}_k- \mathbf{e}_j}}F^{(0)}(0)\nn\\
&&+(3-n)\partial_{\tau^{I}}F^{(0)}(0)\big),
\end{eqnarray}
and  for $I\in \mathbb{Z}_{\geq 0}^{n+1}$,
\begin{eqnarray}\label{eq-recursion-F(0)-WDVV-tau-I}
&&\partial_{\tau^i}\partial_{\tau^j}\partial_{\tau^k}\partial_{\tau^I}F^{(0)}(0)\nn\\
&=& -\sum_{\begin{subarray}{c}0^{n+1}\leq J\leq I\\J\neq 0^{n+1}\end{subarray}}\binom{I}{J}\sum_{a=0}^n\sum_{b=0}^n
\big(\partial_{\tau^{J}}\partial_{\tau^1} \partial_{\tau^{i-1}}\partial_{\tau^a}F^{(0)}(0)\big)\tauG^{ab}\big(\partial_{\tau^{b}}\partial_{\tau^j}\partial_{\tau^k}
\partial_{\tau^{I-J}}F^{(0)}(0)\big)\nn\\
&&+\sum_{0^{n+1}\leq  J\leq I}\binom{I}{J}\sum_{a=0}^n\sum_{b=0}^n\big(\partial_{\tau^{J}}\partial_{\tau^1} \partial_{\tau^{j}}\partial_{\tau^a}F^{(0)}(0)\big)\tauG^{ab}\big(\partial_{\tau^{b}}\partial_{\tau^{i-1}}\partial_{\tau^k}
\partial_{\tau^{I-J}}F^{(0)}(0)\big).
\end{eqnarray}
For $I=(p_0,\dots,p_n)\in \mathbb{Z}_{\geq 0}^{n+1}$ with $|I|\geq 1$,
\begin{equation}\label{eq-recursion-F(1)-EulerField}
\mathsf{a}_X\partial_{\tau^1}\partial_{\tau^I}F^{(1)}(0)=
- \sum_{i=0}^{n}\sum_{j=0}^n\sum_{k=0}^n
(1-i)W_j^i  M_{i}^k p_j \partial_{\tau^{I+\mathbf{e}_k- \mathbf{e}_j}}F^{(1)}(0)+\partial_{\tau^I}F^{(1)}(0).
\end{equation}
For $i\geq 1$ and $I\in \mathbb{Z}_{\geq 0}^{n+1}$, 
\begin{eqnarray}\label{eq-recursion-F(1)-WDVV-tau}
&& \partial_{\tau^I}\partial_{\tau^i}F^{(1)}(0)\nn\\
&=&-\sum_{\begin{subarray}{c}0^{n+1}\leq J\leq I\\J\neq 0^{n+1}\end{subarray}}\binom{I}{J}\sum_{a=0}^n\sum_{b=0}^n
\big(\partial_{\tau^{J}}\partial_{\tau^1} \partial_{\tau^{i-1}}\partial_{\tau^a}F^{(0)}(0)\big)\tauG^{ab}\big(\partial_{\tau^{b}}\partial_{\tau^{I-J}}F^{(1)}(0)\big)\nn\\
&&+ \sum_{0^{n+1}\leq  J\leq I}\binom{I}{J}\big(\partial_{\tau^{J}}\partial_{\tau^1} F^{(1)}(0)\big)\big(
\partial_{\tau^{I-J}}\partial_{\tau^{i-1}}F^{(1)}(0)\big).
\end{eqnarray}
By \cite[Example 6.11]{Hu15}, we have
\begin{lemma}\label{lemeq-M&W-FanoIndex=n-1}
In the case $\nu_{\mathbf{d}}=n-1$, 
 the only off-diagonal entries of the matrix $(M_i^j)$ are
\begin{equation}\label{eq-M-FanoIndex=n-1}
      M_{n-1}^0=- \mathbf{d}!,\
      M_{n}^1=\mathbf{d}!- \mathbf{d}^{\mathbf{d}}.
\end{equation}
Accordingly, the only off-diagonal entries of the matrix $(W_i^j)$ are
\begin{equation}\label{eq-W-FanoIndex=n-1}
      W_{n-1}^0=\mathbf{d}!,\
      W_{n}^1=-\mathbf{d}!+\mathbf{d}^{\mathbf{d}}.
\end{equation}
Moreover,
\begin{equation}\label{eq-F1-leadingTerms-largeFanoIndex}
\sfF^{(1)}(\tau)
=\tau^0-\frac{1}{2}\mathbf{d}!\sum_{\begin{subarray}{c}1\leq i,j\leq n\\
i+j=n\end{subarray}
}
\tau^i \tau^{j}
-\mathbf{d}! \mathbf{d}^{\mathbf{d}}\tau^{n-1}\tau^{n}
+O\big((\tau)^3\big).
\end{equation}
\end{lemma}

We need to compute some higher order leading terms of $F^{(0)}$ and $F^{(1)}$. For $a_1,\dots,a_k\in \mathbb{Z}$, we set
\[
p(a_1,\dots,a_k)=\sharp\{i\in [1,k]: a_i=1\}.
\]
\begin{lemma}\label{lem-proof-sumZ-F(0)}
Suppose $\mathbf{d}=(3)$ or $(2,2)$.
\begin{enumerate}
      \item[(i)] Let $1\leq a,b,c,d\leq n$ satisfying $a+b+c+d=n+1$. Then $F_{a,b,c,d}^{(0)}(0)=0$.
      \item[(ii)] Let $1\leq a,b,c\leq n$ satisfying $a+b+c=n$. Then $F_{a,b,c,n}^{(0)}(0)=\mathbf{d}!\prod_{i=1}^r d_i$.
      \item[(iii)] Let $1\leq a,b,c\leq n$ satisfying $a+b+c=n+1$.  
      Then 
      \begin{equation}
            F_{a,b,c,n-1}^{(0)}(0)=\big(\mathbf{d}^{\mathbf{d}}- p(a,b,c)\mathbf{d}!\big)\prod_{i=1}^r d_i.
      \end{equation}
      \item[(iv)] $F_{1,n-1,n-1,n}(0)=(\mathbf{d}^{\mathbf{d}}- \mathbf{d}!)\mathbf{d}^{\mathbf{d}}\prod_{i=1}^r d_i$. 
      $F_{2,n-1,n-1,n-1}(0)=(\mathbf{d}^{\mathbf{d}})^2 \prod_{i=1}^r d_i$.
      \item[(v)] $F_{1,1,1,n-1,n-1}^{(0)}(0)=(\mathbf{d}^{\mathbf{d}}-2 \mathbf{d}!)\prod_{i=1}^r d_i$,
      $F_{1,1,n-1,n-1,n}^{(0)}(0)
      =(\mathbf{d}^{\mathbf{d}}-\mathbf{d}!) (\mathbf{d}^{\mathbf{d}}+2 \mathbf{d}!)\prod_{i=1}^r d_i$.
\end{enumerate}
\end{lemma}
\begin{proof}
Suppose $I=(p_0,\dots,p_n)\in \mathbb{Z}_{\geq 0}^{n+1}$, satisfying $|I|=\sum_i p_i=3$. Then by (\ref{eq-recursion-F(0)-EulerField-I}) and Lemma \ref{lemeq-M&W-FanoIndex=n-1},
\begin{eqnarray*}
&&\partial_{\tau^{I+\mathbf{e}_1}}F^{(0)}(0)\\
&=& \frac{1}{n-1}\big(-\sum_{i=0}^{n}(1-i)p_i\partial_{\tau^{I}}F^{(0)}(0)
-W_{n-1}^0p_{n-1}\partial_{\tau^{I+\mathbf{e}_0- \mathbf{e}_{n-1}}}F^{(0)}(0)\\
&&-(2-n) M_{n-1}^0 p_{n-1}\partial_{\tau^{I+\mathbf{e}_0- \mathbf{e}_{n-1}}}F^{(0)}(0)
-(1-n)M_{n}^1 p_n\partial_{\tau^{I+\mathbf{e}_1- \mathbf{e}_n}}F^{(0)}(0)
\nn\\
&&+(3-n)\partial_{\tau^{I}}F^{(0)}(0)\big)\\
&=& \frac{\sum_i ip_i-n}{n-1}\partial_{\tau^{I}}F^{(0)}(0)+M_{n-1}^0 p_{n-1}\partial_{\tau^{I+\mathbf{e}_0- \mathbf{e}_{n-1}}}F^{(0)}(0)
+M_{n}^1 p_n\partial_{\tau^{I+\mathbf{e}_1- \mathbf{e}_n}}F^{(0)}(0).
\end{eqnarray*}

(i) Suppose $\sum_i ip_i=n$. By (\ref{eq-qp1.5}), $\partial_{\tau^{I+\mathbf{e}_0- \mathbf{e}_{n-1}}}F^{(0)}(0)=\partial_{\tau^{I+\mathbf{e}_0- \mathbf{e}_{n-1}}}F^{(0)}(0)=\partial_{\tau^{I+\mathbf{e}_1- \mathbf{e}_n}}F^{(0)}(0)=0$, thus
\begin{eqnarray*}
\partial_{\tau^{I+\mathbf{e}_1}}F^{(0)}(0)
=\frac{1}{n-1}\big(-\sum_{i=0}^{n}(1-i)p_i\partial_{\tau^{I}}F^{(0)}(0)
+(3-n)\partial_{\tau^{I}}F^{(0)}(0)\big)=0.
\end{eqnarray*}
Then by (\ref{eq-recursion-F(0)-WDVV-tau-I}) and (\ref{eq-qp1.5}) one can show that $\partial_{\tau^{I}}F^{(0)}(0)=0$  for all $I=(p_0,\dots,p_n)$ satisfying 
$|I|=4$ and $\sum_i ip_i=n+1$.

(ii) Suppose $\sum_i ip_i=2n-1$. Then by (\ref{eq-qp1.5}) and Lemma \ref{lemeq-M&W-FanoIndex=n-1},
\begin{eqnarray}\label{lem-proof-sumZ-F(0)-sumipi=2n-1}
\partial_{\tau^{I+\mathbf{e}_1}}F^{(0)}(0)
= \mathbf{d}^{\mathbf{d}}\prod_{i=1}^r d_i- p_{n-1}\mathbf{d}!\prod_{i=1}^r d_i
+p_n(\mathbf{d}!- \mathbf{d}^{\mathbf{d}})\prod_{i=1}^r d_i.
\end{eqnarray}
If $p_n=1$ and $p_{n-1}=0$, then 
$\partial_{\tau^{I+\mathbf{e}_1}}F^{(0)}(0)
= \mathbf{d}!\prod_{i=1}^r d_i$.
This shows (ii) when $a=1$. Then (\ref{eq-recursion-F(0)-WDVV-tau-I}) and (\ref{eq-qp1.5}) yield (ii).

(iii) If $a=1$, then $p(a,b,c)=p_{n-1}$ and $p_n=0$. So (\ref{lem-proof-sumZ-F(0)-sumipi=2n-1}) yields (iii) when $a=1$. Then (\ref{eq-recursion-F(0)-WDVV-tau-I}) and (\ref{eq-qp1.5}) yield (iii).

(iv) If $p_{n-1}=2$ and $p_n=1$, then
\begin{eqnarray*}
&&\partial_{\tau^{I+\mathbf{e}_1}}F^{(0)}(0)\\
&=& 2\partial_{\tau^{I}}F^{(0)}(0)+2M_{n-1}^0 \partial_{\tau^{I+\mathbf{e}_0- \mathbf{e}_{n-1}}}F^{(0)}(0)
+M_{n}^1 \partial_{\tau^{I+\mathbf{e}_1- \mathbf{e}_n}}F^{(0)}(0)\\
&=& 2(\mathbf{d}^{\mathbf{d}})^2\prod_{i=1}^r d_i-2 \mathbf{d}!\mathbf{d}^{\mathbf{d}}\prod_{i=1}^r d_i
+(\mathbf{d}!- \mathbf{d}^{\mathbf{d}})\mathbf{d}^{\mathbf{d}}\prod_{i=1}^r d_i
=(\mathbf{d}^{\mathbf{d}}- \mathbf{d}!)\mathbf{d}^{\mathbf{d}}\prod_{i=1}^r d_i.
\end{eqnarray*}
 (\ref{eq-recursion-F(0)-WDVV-tau-I}) and (\ref{eq-qp1.5}) yield
\begin{eqnarray*}
&&\partial_{\tau^2}\partial_{\tau^{n-1}}\partial_{\tau^{n-1}}\partial_{\tau^{n-1}}F^{(0)}(0)\nn\\
&=& -\partial_{\tau^{n-1}}\partial_{\tau^1} \partial_{\tau^{1}}\partial_{\tau^{n-1}}F^{(0)}(0)\tauG^{n-1,1}\partial_{\tau^{1}}\partial_{\tau^{n-1}}\partial_{\tau^{n-1}}F^{(0)}(0)\\
&&+2\big(\partial_{\tau^{n-1}}\partial_{\tau^1} \partial_{\tau^{n-1}}\partial_{\tau^1}F^{(0)}(0)\tauG^{1,0}\partial_{\tau^{0}}\partial_{\tau^{1}}\partial_{\tau^{n-1}}F^{(0)}(0)\\
&&+\partial_{\tau^{n-1}}\partial_{\tau^1} \partial_{\tau^{n-1}}\partial_{\tau^n}F^{(0)}(0)\tauG^{n,0}\partial_{\tau^{0}}\partial_{\tau^{1}}\partial_{\tau^{n-1}}F^{(0)}(0)\\
&&+\partial_{\tau^{n-1}}\partial_{\tau^1} \partial_{\tau^{n-1}}\partial_{\tau^1}F^{(0)}(0)\tauG^{1,n-1}\partial_{\tau^{n-1}}\partial_{\tau^{1}}\partial_{\tau^{n-1}}F^{(0)}(0)\big)\\
&=& (\mathbf{d}^{\mathbf{d}})^2 \prod_{i=1}^r d_i.
\end{eqnarray*}

(v) Suppose $I=(p_0,\dots,p_n)\in \mathbb{Z}_{\geq 0}^{n+1}$, satisfying $|I|=\sum_i p_i=4$. Then
\begin{eqnarray*}
&&\partial_{\tau^{I}} \partial_{\tau^1}F^{(0)}(0)\\
&=&\frac{1}{n-1}\big(-\sum_{i=0}^{n}(1-i) p_i\partial_{\tau^{I}}F^{(0)}(0)
-(1-n) M_{n}^1 p_n\partial_{\tau^{I+\mathbf{e}_1- \mathbf{e}_n}}F^{(0)}(0)
-W_{n-1}^0  p_{n-1}\partial_{\tau^{I+\mathbf{e}_0- \mathbf{e}_{n-1}}}F^{(0)}(0)\\
&&-(2-n)  M_{n-1}^0 p_{n-1}\partial_{\tau^{I+\mathbf{e}_0- \mathbf{e}_{n-1}}}F^{(0)}(0)
+(3-n)\partial_{\tau^{I}}F^{(0)}(0)\big).
\end{eqnarray*}
By the fundamental class axiom, $\partial_{\tau^{I+\mathbf{e}_0- \mathbf{e}_{n-1}}}F^{(0)}(0)=0$. Thus
\begin{eqnarray*}
\partial_{\tau^{I}} \partial_{\tau^1}F^{(0)}(0)
=\frac{\sum_i i p_i-n-1}{n-1}\partial_{\tau^{I}}F^{(0)}(0)
+ M_{n}^1 p_n\partial_{\tau^{I+\mathbf{e}_1- \mathbf{e}_n}}F^{(0)}(0).
\end{eqnarray*}
If $p_1=p_{n-1}=2$, then
\begin{eqnarray*}
\partial_{\tau^{I}} \partial_{\tau^1}F^{(0)}(0)=(\mathbf{d}^{\mathbf{d}}-2 \mathbf{d}!)\prod_{i=1}^r d_i.
\end{eqnarray*}
If $p_1=p_n=1$ and $p_{n-1}=2$, then
\begin{eqnarray*}
\partial_{\tau^{I}} \partial_{\tau^1}F^{(0)}(0)
=2(\mathbf{d}^{\mathbf{d}}- \mathbf{d}!)\mathbf{d}^{\mathbf{d}}\prod_{i=1}^r d_i
+ (\mathbf{d}!- \mathbf{d}^{\mathbf{d}}) (\mathbf{d}^{\mathbf{d}}-2 \mathbf{d}!)\prod_{i=1}^r d_i
=(\mathbf{d}^{\mathbf{d}}-\mathbf{d}!) (\mathbf{d}^{\mathbf{d}}+2 \mathbf{d}!)\prod_{i=1}^r d_i.
\end{eqnarray*}
\end{proof}

\begin{lemma}\label{lem-proof-sumZ-F(1)}
Let $X$ be a cubic hypersurface of dimension $n\geq 3$. 
\begin{enumerate}
       \item[(i)] Let $1\leq a,b,c\leq n$ satisfying $a+b+c=n+1$. Then $F^{(1)}_{a,b,c}(0)=- \mathbf{d}!$.
      \item[(ii)] $F_{1,n-1,n}^{(1)}(0)=- \mathbf{d}!(\mathbf{d}^{\mathbf{d}}+\mathbf{d}!)$,
       $F_{2,n-1,n-1}^{(1)}(0)=- \mathbf{d}!\mathbf{d}^{\mathbf{d}}  $.
\end{enumerate}
\end{lemma}
\begin{proof}
Suppose $I=(p_0,\dots,p_n)\in \mathbb{Z}_{\geq 0}^{n+1}$, satisfying $|I|=\sum_i p_i=2$ Then
\begin{eqnarray*}
&&(n-1)\partial_{\tau^1}\partial_{\tau^I}F^{(1)}(0)\\
&=&- \sum_{i=0}^{n}(1-i) p_i \partial_{\tau^{I}}F^{(1)}(0)
- (1-n) M_{n}^1 p_n \partial_{\tau^{I+\mathbf{e}_1- \mathbf{e}_n}}F^{(1)}(0)\\
&&- W_{n-1}^0  p_0 \partial_{\tau^{I+\mathbf{e}_0- \mathbf{e}_{n-1}}}F^{(1)}(0)
-(2-n) M_{n-1}^0 p_{n-1} \partial_{\tau^{I+\mathbf{e}_0- \mathbf{e}_{n-1}}}F^{(1)}(0)
+\partial_{\tau^I}F^{(1)}(0).
\end{eqnarray*}
By the fundamental class axiom, $\partial_{\tau^{I+\mathbf{e}_0- \mathbf{e}_{n-1}}}F^{(1)}(0)=0$. So
\begin{eqnarray}\label{lem-proof-sumZ-F(1)-0}
\partial_{\tau^1}\partial_{\tau^I}F^{(1)}(0)
=\frac{\sum_i i p_i-1}{n-1} \partial_{\tau^{I}}F^{(1)}(0)
+ M_{n}^1 p_n \partial_{\tau^{I+\mathbf{e}_1- \mathbf{e}_n}}F^{(1)}(0).
\end{eqnarray}
(i) Suppose $\sum_i i p_i=n$, and $p_0=0$. Then $p_n=0$. Then (\ref{lem-proof-sumZ-F(1)-0}) and (\ref{eq-F1-leadingTerms-largeFanoIndex}) yields $\partial_{\tau^1}\partial_{\tau^I}F^{(1)}(0)=- \mathbf{d}!$.
(ii) If $p_{n-1}=p_{n}=1$, then  (\ref{lem-proof-sumZ-F(1)-0}) and (\ref{eq-F1-leadingTerms-largeFanoIndex}) yields
\begin{eqnarray*}
\partial_{\tau^1}\partial_{\tau^I}F^{(1)}(0)
&=& 2 \partial_{\tau^{I}}F^{(1)}(0)+ M_{n}^1  \partial_{\tau^{1}}\partial_{\tau^{n-1}}F^{(1)}(0)\\
&=&-2 \mathbf{d}! \mathbf{d}^{\mathbf{d}}+(\mathbf{d}!- \mathbf{d}^{\mathbf{d}})(- \mathbf{d}!)
=- \mathbf{d}!(\mathbf{d}^{\mathbf{d}}+\mathbf{d}!).
\end{eqnarray*}
By (\ref{eq-recursion-F(1)-WDVV-tau}) and Lemma \ref{lem-proof-sumZ-F(0)} we get 
\begin{eqnarray*}
&& \partial_{\tau^{n-1}}\partial_{\tau^{n-1}}\partial_{\tau^2}F^{(1)}(0)\nn\\
&=&-2\big(\partial_{\tau^{n-1}}\partial_{\tau^1} \partial_{\tau^{1}}\partial_{\tau^{n-1}}F^{(0)}(0)\big)\tauG^{n-1,1}\big(\partial_{\tau^{1}}\partial_{\tau^{n-1}}F^{(1)}(0)\big)\nn\\
&&-\big(\partial_{\tau^{n-1}}\partial_{\tau^{n-1}}\partial_{\tau^1} \partial_{\tau^{1}}\partial_{\tau^1}F^{(0)}(0)\big)\tauG^{1,0}\big(\partial_{\tau^{0}}F^{(1)}(0)\big)\nn\\
&&-\big(\partial_{\tau^{n-1}}\partial_{\tau^{n-1}}\partial_{\tau^1} \partial_{\tau^{1}}\partial_{\tau^n}F^{(0)}(0)\big)\tauG^{n,0}\big(\partial_{\tau^{0}}F^{(1)}(0)\big)\nn\\
&&+ 2\big(\partial_{\tau^{n-1}}\partial_{\tau^1} F^{(1)}(0)\big)\big(
\partial_{\tau^{n-1}}\partial_{\tau^{1}}F^{(1)}(0)\big)\\
&=& - \mathbf{d}!\mathbf{d}^{\mathbf{d}}.
\end{eqnarray*}
\end{proof}

\begin{proof}[Proof of Lemma \ref{lem-CubicHypersurfaces-(2,2)-sumZ}]
Suppose $1\leq a,b,c,d\leq n$ and $a+b+c+d=n+1$. Thus $a,b,c,d\leq n-2$ and the sum of any two of $a,b,c,d$ is $\leq n-1$. Then by (\ref{eq-F1-leadingTerms-largeFanoIndex}) we have 
\begin{equation}\label{eq-lem-CubicHypersurfaces-(2,2)-sumZ-1}
      F_{a,b}^{(1)}(0)y_{c+d}+F_{a,c}^{(1)}(0)y_{b+d}
+F_{a,d}^{(1)}(0)y_{b+c}+F_{c,d}^{(1)}(0)y_{a+b}+F_{b,d}^{(1)}(0)y_{a+c}
+F_{b,c}^{(1)}(0)y_{a+d}=0,
\end{equation}
\begin{eqnarray}\label{eq-lem-CubicHypersurfaces-(2,2)-sumZ-2}
&& F_{a+b,c}^{(1)}(0) y_d+F_{a+b,d}^{(1)}(0) y_c+F_{a+c,b}^{(1)}(0) y_d
+F_{a+c,d}^{(1)}(0) y_b\nn\\
&&+F_{a+d,b}^{(1)}(0) y_c+F_{a+d,c}^{(1)}(0) y_b
+F_{b+c,a}^{(1)}(0) y_d+F_{b+c,d}^{(1)}(0) y_a\nn\\
&&+F_{b+d,a}^{(1)}(0) y_c+F_{b+d,c}^{(1)}(0) y_a
+F_{c+d,a}^{(1)}(0) y_b+F_{c+d,b}^{(1)}(0) y_a\nn\\
&=& -3 p(a,b,c,d)\mathbf{d}!y_1,
\end{eqnarray}
and
\begin{equation}\label{eq-lem-CubicHypersurfaces-(2,2)-sumZ-3}
      F_{a,b}^{(1)}(0)F_{c,d}^{(1)}(0)+F_{a,c}^{(1)}(0)F_{b,d}^{(1)}(0)
+F_{a,d}^{(1)}(0)F_{b,c}^{(1)}(0)=0.
\end{equation}
By Lemma \ref{lem-proof-sumZ-F(0)} and (\ref{eq-pairing2}), we have
\begin{eqnarray}\label{eq-lem-CubicHypersurfaces-(2,2)-sumZ-4}
&& \big(F_{a+b,c,d,f}^{(0)}(0)\tauG^{fi}+F_{a+c,b,d,f}^{(0)}(0)\tauG^{fi}+F_{a+d,b,c,f}(0)\tauG^{fi}\nn\\
&&+F_{c+d,a,b,f}^{(0)}(0)\tauG^{fi}+F_{b+d,a,c,f}^{(0)}(0)\tauG^{fi}+F_{b+c,a,d,f}^{(0)}(0)\tauG^{fi}\big)y_{i}\nn\\
&=& \begin{cases}
\mathbf{d}^{\mathbf{d}}\prod_{i=1}^r d_i\times6\times \tauG^{n-1,i}y_{i},& \mbox{if}\ p(a,b,c,d)=0\\
(\mathbf{d}^{\mathbf{d}}\prod_{i=1}^r d_i\times 3+(\mathbf{d}^{\mathbf{d}}- \mathbf{d}!)\prod_{i=1}^r d_i\times 3)\times \tauG^{n-1,i}y_{i},& \mbox{if}\ p(a,b,c,d)=1\\
(\mathbf{d}^{\mathbf{d}}\prod_{i=1}^r d_i+(\mathbf{d}^{\mathbf{d}}- \mathbf{d}!)\prod_{i=1}^r d_i\times4+(\mathbf{d}^{\mathbf{d}}- 2\mathbf{d}!)\prod_{i=1}^r d_i)\times \tauG^{n-1,i}y_{i},& \mbox{if}\ p(a,b,c,d)=2\\
((\mathbf{d}^{\mathbf{d}}- \mathbf{d}!)\prod_{i=1}^r d_i\times3+(\mathbf{d}^{\mathbf{d}}- 2\mathbf{d}!)\prod_{i=1}^r d_i\times 3)\times \tauG^{n-1,i}y_{i},& \mbox{if}\ p(a,b,c,d)=3
\end{cases}\nn\\
&=& \big(6 \mathbf{d}^{\mathbf{d}}-3 \mathbf{d}! p(a,b,c,d)\big)y_1,
\end{eqnarray}
and
\begin{eqnarray}\label{eq-lem-CubicHypersurfaces-(2,2)-sumZ-5}
&& F_{b,c,d,e}^{(0)}(0)\tauG^{ej}y_{j+a}=\begin{cases}
F_{b,c,d,a}^{(0)}(0)\tauG^{a,j}y_{j+a},& \mbox{if}\ a>1\\
F_{b,c,d,1}^{(0)}(0)\tauG^{1,j}y_{j+1}+F_{b,c,d,n}(0)\tauG^{n,j}y_{j+1},& \mbox{if}\ a=1
\end{cases}\nn\\
&=& \begin{cases}
\frac{1}{\prod_{i=1}^r d_i}F_{a,b,c,d}^{(0)}(0)y_{n},& \mbox{if}\ a>1\\
\frac{1}{\prod_{i=1}^r d_i}\big(F_{b,c,d,1}^{(0)}(0)y_{n}- \mathbf{d}^{\mathbf{d}}F_{b,c,d,1}^{(0)}y_1+F_{b,c,d,n}^{(0)}(0)y_{1}\big),& \mbox{if}\ a=1
\end{cases}\nn\\
&=& \begin{cases}
0,& \mbox{if}\ a>1\\
\mathbf{d}! y_1,& \mbox{if}\ a=1
\end{cases}.
\end{eqnarray}
By Lemma \ref{lem-proof-sumZ-F(1)} (i),
\begin{equation}\label{eq-lem-CubicHypersurfaces-(2,2)-sumZ-6}
      F_{a+b,c,d}^{(1)}(0)+F_{a+c,b,d}^{(1)}(0)
+F_{a+d,b,c}^{(1)}(0)+F_{c+d,a,b}^{(1)}(0)
+F_{b+d,a,c}^{(1)}(0)+F_{b+c,a,d}^{(1)}(0)=-6 \mathbf{d}!.
\end{equation}
By (\ref{eq-F1-leadingTerms-largeFanoIndex}),  (\ref{eq-pairing2}) and Lemma \ref{lem-proof-sumZ-F(0)} we have
\begin{eqnarray}\label{eq-lem-CubicHypersurfaces-(2,2)-sumZ-7}
&& F_{b,c,d,e}^{(0)}(0)\tauG^{e,i}F_{i,a}^{(1)}(0)
=\begin{cases}
F_{b,c,d,a}^{(0)}(0)\tauG^{a,i}F_{i,a}^{(1)}(0),& \mbox{if}\ a>1\\
F_{b,c,d,1}^{(0)}(0)\tauG^{1,i}F_{i,1}^{(1)}(0)+F_{b,c,d,n}^{(0)}(0)\tauG^{n,i}F_{i,1}^{(1)}(0),& \mbox{if}\ a=1
\end{cases}\nn\\
&=&\begin{cases}
\frac{1}{\prod_{i=1}^r d_i}F_{b,c,d,a}^{(0)}(0)F_{n-a,a}^{(1)}(0),& \mbox{if}\ a>1\\
\frac{1}{\prod_{i=1}^r d_i}\big(F_{b,c,d,1}^{(0)}(0)F_{n-1,1}^{(1)}(0)- \mathbf{d}^{\mathbf{d}} F_{b,c,d,1}^{(0)}(0)F_{0,1}^{(1)}(0)
+F_{b,c,d,n}^{(0)}(0)F_{0,1}^{(1)}(0)\big),& \mbox{if}\ a=1
\end{cases}\nn\\
&=& 0.
\end{eqnarray}
Applying (\ref{eq-lem-CubicHypersurfaces-(2,2)-sumZ-1})-(\ref{eq-lem-CubicHypersurfaces-(2,2)-sumZ-7}) to (\ref{eq-def-zabcd}) we get
\begin{eqnarray*}
w_{a,b,c,d}
&=& -\frac{1}{6}\big(6 \mathbf{d}^{\mathbf{d}}-3 \mathbf{d}! p(a,b,c,d)\big)y_1
+\frac{1}{2}\times p(a,b,c,d) \mathbf{d}!y_1\nn\\
&&+\frac{2}{3}N\times 3 \mathbf{d}!p(a,b,c,d)y_1-\frac{1}{3}N (-6 \mathbf{d}!) y\nn\\
&=& \big(2 Np(a,b,c,d)\mathbf{d}!+p(a,b,c,d)\mathbf{d}!- \mathbf{d}^{\mathbf{d}}\big)y_1+2N \mathbf{d}!y.
\end{eqnarray*}
It follows that
\begin{eqnarray}\label{eq-lem-CubicHypersurfaces-(2,2)-sumZ-1stPart}
\sum_{i=1}^{n-2}w_{i,n-1-i,1,1}
=\big(4(n-1)\mathbf{d}!N+(2\mathbf{d}!- \mathbf{d}^{\mathbf{d}})n+2 \mathbf{d}^{\mathbf{d}}-2 \mathbf{d}!\big)y_1+2(n-2)N \mathbf{d}!y.
\end{eqnarray}
Moreover (\ref{eq-F1-leadingTerms-largeFanoIndex}),  (\ref{eq-pairing2}), Lemma \ref{lem-proof-sumZ-F(0)} and \ref{lem-proof-sumZ-F(1)} yield
\begin{eqnarray}\label{eq-lem-CubicHypersurfaces-(2,2)-sumZ-2ndPart}
&&w_{n-1,n-1,1,1}\nn\\
&=&8N F_{n-1,1}^{(1)}(0)y_n
-\frac{1}{6}\big(F_{2n-2,1,1,n-1}^{(0)}(0)+4F_{n,n-1,1,n-1}^{(0)}(0)+F_{2,n-1,n-1,n-1}^{(0)}(0)\big)\tauG^{n-1,1}y_1\nn\\
&&+\big(2F_{n-1,1,1,n-1}^{(0)}(0)\tauG^{n-1,1}y_n+F_{n-1,n-1,1,1}^{(0)}(0)\tauG^{1,0}y_{1}
+F_{n-1,n-1,1,n}^{(0)}(0)\tauG^{n,0}y_{1}\big)\nn\\
&&-\frac{2}{3}N\times 6F_{n,n-1}^{(1)}(0)y_1
+\Big(\frac{8}{3}N(3N-2)F_{n-1,1}^{(1)}(0)F_{n-1,1}^{(1)}(0)\nn\\
&&-\frac{1}{3}N\big(F_{2n-2,1,1}^{(1)}(0)+4F_{n,n-1,1}^{(1)}(0)+F_{2,n-1,n-1}^{(1)}(0)\big)\nn\\
&&+2N\big(F_{n-1,1,1,n-1}^{(0)}(0)\tauG^{n-1,1}F_{1,n-1}^{(1)}(0)+F_{n-1,n-1,1,1}^{(0)}(0)\tauG^{1,n-1}F_{n-1,1}^{(1)}(0)\big)\Big)y\nn\\
&=& (4\mathbf{d}! N-\mathbf{d}^{\mathbf{d}}+2 \mathbf{d}!)(-2y_n+\mathbf{d}^{\mathbf{d}}y_1+2N \mathbf{d}! y).
\end{eqnarray}
Hence (\ref{eq-CubicHypersurfaces-(2,2)-sumZ}) follows from (\ref{eq-lem-CubicHypersurfaces-(2,2)-sumZ-1stPart}) and (\ref{eq-lem-CubicHypersurfaces-(2,2)-sumZ-2ndPart}).
\end{proof}

\end{appendices}

\textsc{School of Mathematics, Sun Yat-sen University, Guangzhou 510275, P.R. China}

 \emph{E-mail address:}  huxw06@gmail.com

\end{document}